\documentclass{article}
\usepackage{graphicx}
\usepackage[T1]{fontenc}
\usepackage[english]{babel}
\usepackage{graphicx}
\usepackage{subfigure}
\usepackage{psfrag}
\usepackage{amssymb, amsmath, amsfonts}
\usepackage{dsfont}
\usepackage{vmargin}
\usepackage{makeidx}
\usepackage{color}

\everymath{\displaystyle}

\newcommand{\SDP}{Simple_Double_Puits.eps}
\newcommand{\SPINOQUARTICUN}{SpinoQuartic1.eps}
\newcommand{\SPINOQUARTICDEUX}{SpinoQuartic2.eps}
\newcommand{\SPINOQUARTICTROIS}{SpinoQuartic3.eps}
\newcommand{\SPINOQUARTICQUATRE}{SpinoQuartic4.eps}
\newcommand{\SPINOQUARTICCINQ}{SpinoQuartic5.eps}
\newcommand{\SPINOQUARTICSIX}{SpinoQuartic6.eps}
\newcommand{\SPINOLOGUN}{SpinoLog1.eps}
\newcommand{\SPINOLOGDEUX}{SpinoLog2.eps}
\newcommand{\SPINOLOGTROIS}{SpinoLog3.eps}
\newcommand{\SPINOLOGQUATRE}{SpinoLog4.eps}
\newcommand{\SPINOLOGCINQ}{SpinoLog5.eps}
\newcommand{\SPINOLOGSIX}{SpinoLog6.eps}
\newcommand{\CROIX}{Croix.eps}
\newcommand{\BULLE}{Bulle.eps}
\newcommand{\CROIXLOG}{Croix.eps}
\newcommand{\BULLELOG}{BulleLog.eps}
\newcommand{\Points}{Points.eps}
\newcommand{\ConvergencePoints}{ConvergencePoints.eps}
\newcommand{\Interface}{Interface.eps}
\newcommand{\FitCurveExemple}{FitCurveM35Q3.eps}
\newcommand{\LongueurInterface}{LongueurInterface.eps}
\newcommand{\Enerlog}{EnergiesLog.eps}
\newcommand{\ErrLDEUXlog}{ErreursL2PolyLog.eps}
\newcommand{\EnerlogDEUXD}{EnergiesLog2D.eps}
\newcommand{\ErrLDEUXlogDEUXD}{ErreursL2PolyLog2D.eps}
\newcommand{\SpinoAAA}{SpinodaleT1.eps}
\newcommand{\SpinoAAB}{SpinodaleT2.eps}
\newcommand{\SpinoAAC}{SpinodaleT3.eps}
\newcommand{\SpinoAAD}{SpinodaleT4.eps}
\newcommand{\SpinoAAE}{SpinodaleT5.eps}
\newcommand{\SpinoAAF}{SpinodaleT6.eps}
\newcommand{\MOVINGENERGY}{MovingEnergy.eps}
\newcommand{\FitCurveMAQAZZ}{FitCurveM18Q1.eps}
\newcommand{\FitCurveMBQAZZ}{FitCurveM36Q1.eps}
\newcommand{\FitCurveMCQAZZ}{FitCurveM72Q1.eps}
\newcommand{\FitCurveMAQBZZ}{FitCurveM9Q2.eps}
\newcommand{\FitCurveMBQBZZ}{FitCurveM18Q2.eps}
\newcommand{\FitCurveMCQBZZ}{FitCurveM36Q2.eps}
\newcommand{\FitCurveMAQCZZ}{FitCurveM6Q3.eps}
\newcommand{\FitCurveMBQCZZ}{FitCurveM12Q3.eps}
\newcommand{\FitCurveMCQCZZ}{FitCurveM24Q3.eps}
\newcommand{\SlopesQWMeshW}{SlopesQWMeshW.eps}
\newcommand{\SlopesQWLogMeshW}{SlopesQWLogMeshW.eps}
\newcommand{\SlopesDLA}{SlopesDL90.eps}
\newcommand{\SlopesDLB}{SlopesDL180.eps}
\newcommand{\SlopesDLC}{SlopesDL300.eps}
\newcommand{\ErrEnergiesQZIV}{ErrEnergiesQ1234.eps}
\newcommand{\ErrEnergiesQZIVLog}{ErrEnergiesQ1234Log.eps}
\newcommand{\ErrEnergiesDLWW}{ErrEnergiesDL90180300630.eps}
\newcommand{\ErrEnergiesQZVDEUXD}{ErrEnergies2DQ12345.eps}
\newcommand{\ErrEnergiesQZVDEUXDLog}{ErrEnergies2DQ12345Log.eps}
\newcommand{\ErreursQWMeshW}{ErreursL2Q12345.eps}
\newcommand{\ErreursQWLogMeshW}{ErreursL2Q12345Log.eps}
\newcommand{\Bosses}{DessinsInterfaces.eps}
\newcommand{\BANDE}{Bande-Blanc.eps}
\newcommand{\CERCLE}{Cercle-Blanc.eps}
\newcommand{\MODEONZE}{Mode11-Blanc.eps}
\newcommand{\HAUTMODE}{Mode4114-Blanc.eps}
\newcommand{\BifurcationsRectangleOrdre}{BifurcationsRectangleOrdre.eps}
\newcommand{\ConvergenceCstRectangle}{ConvergenceCstRectangle.eps}
\newcommand{\BifurcationsSegmentOrdre}{BifurcationsSegmentOrdre.eps}
\newcommand{\ConvergenceCstSegment}{ConvergenceCstSegment.eps}
\newcommand{\EigenEllipse}{EigenEllipse.ps} 
\newcommand{\SingularEllipse}{SingularEllipse.eps}
\newcommand{\EigenTrapezoid}{EigenTrapezoid1D.eps} 
\newcommand{\BifurcationsTrapezeOrdre}{BifurcationsTrapezeOrdre.eps}
\newcommand{\CercleMD}{CercleMD.eps}
\newcommand{\BandeH}{BandeH.eps}
\newcommand{\CercleBG}{CercleBG.eps}
\newcommand{\CercleHG}{CercleHG.eps}
\newcommand{\EnergiesTrapezoid}{EnergiesTrapezoid.eps}

\usepackage{hyperref}
\hypersetup{
	colorlinks=true,
	breaklinks=true,
	urlcolor=blue,
	linkcolor=red,
	bookmarksopen=true,
	pdftitle = High order finite element calculations for the deterministic Cahn-Hilliard equation,
	pdfauthor = {Ludovic Gouden\`ege, Daniel Martin, Gr\'egory Vial}
	pdfsubject = High order finite element calculations for the deterministic Cahn-Hilliard equation,
	pdfkeywords = {Cahn-Hilliard, partial differential equations, bi-laplacian, {\it p}-version of the finite element method, spinodal decomposition}
}

\numberwithin{equation}{section} 

\newcommand{\R}{\mathbb{R}} 
\newcommand{\N}{\mathbb{N}} 
\newcommand{\PP}{\mathbb{P}} 
\newcommand{\dd}{\mathrm{d}} 

\newcommand{\sL}{\mathrm{L}}
\newcommand{\sH}{\mathrm{H}}

\newcommand{\BD}{\begin{displaymath}}
\newcommand{\ED}{\end{displaymath}}
\newcommand{\BEA}{\begin{eqnarray}}
\newcommand{\EEA}{\end{eqnarray}}
\newcommand{\BEAS}{\begin{eqnarray*}}
\newcommand{\EEAS}{\end{eqnarray*}}
\newcommand{\BE}{\begin{equation}}
\newcommand{\EE}{\end{equation}}
\newcommand{\BES}{\begin{equation*}}
\newcommand{\EES}{\end{equation*}}

\title{High order finite element calculations for the deterministic Cahn-Hilliard equation}
\author{Ludovic Gouden\`ege\footnotemark[1], Daniel Martin\footnotemark[2], Gr\'egory Vial\footnotemark[1]}
\date{}

\begin{document}

\maketitle
\renewcommand{\thefootnote}{\*}
\footnotetext[0]{\!\!\!\!\!\!\!\!\!\!\!AMS 2000 subject classifications. 65M60, 35K55, 35K65, 35K45, 82C26\\
{\em Key words and phrases} : Cahn-Hilliard, partial differential equations, bi-laplacian, $p$-version of the finite element method, spinodal decomposition, bifurcations.\\}
\renewcommand{\thefootnote}{\fnsymbol{footnote}}
\footnotetext[1]{IRMAR, ENS Cachan Bretagne, CNRS, UEB, 35170 Bruz, France, firstname.name@bretagne.ens-cachan.fr}
\footnotetext[2]{IRMAR, Universit\'e de Rennes 1, Campus de Beaulieu, 35000 Rennes, France, daniel.martin@univ-rennes1.fr}

\begin{abstract}
In this work, we propose a numerical method based on high degree continuous nodal elements for the Cahn-Hilliard evolution. The use of the $p$-version of the finite element method proves to be very efficient and favorably compares with other existing strategies ($\mathcal C^1$ elements, adaptive mesh refinement, multigrid resolution, etc). Beyond the classical benchmarks, a numerical study has been carried out to investigate the influence of a polynomial approximation of the logarithmic free energy and the bifurcations near the first eigenvalue of the Laplace operator.
\end{abstract}

\section*{Introduction}\label{S:0}

We consider an isothermal binary alloy of two species $A$ and $B$, and denote by $u\in [-1,1]$ the ratio between the two components. By thermodynamic arguments, and under a mass conservation property, Cahn and Hilliard described a fourth-order model for the evolution of an isotropic system of nonuniform composition or density. They introduced a free energy density $\bar{f}$ to define a chemical potential, and use it in the classical transport equation (see \cite{C}, \cite{CH1} and \cite{CH2}). The total free energy $\mathcal{F}$ of the binary alloy is a volume integral on $\Omega$ of this free energy density (bulk free energy):
\BE\label{Eq:0.1}
\mathcal{F} := \int_{\Omega} \bar{f}(u,\nabla u, \nabla^2u, \dots) \  \dd V.
\EE
They assumed $\bar{f}$ to be a function of $u$ and its spatial derivatives. 
A truncated Taylor expansion of $\bar{f}$ has thus the following general form:
\BE\label{Eq:0.2}
\bar{f}(u) \sim f(u) + L \cdot \nabla u + K_{1}\otimes \nabla^2u+\nabla u \cdot K_{2} \cdot \nabla u,
\EE
where $\nabla$ is the Nabla operator.
By symmetry arguments, they showed that $L = \vec{0}$ and $K_{1}$ and $K_{2}$ are homothetic operators. Moreover they used Neumann boundary 
condition to cancel the term in $\nabla^2u$ which yields 
\BE\label{Eq:0.3}
\mathcal{F} := \int_{\Omega} \left( f(u) + \kappa |\nabla u|^2 \right)\dd V,
\EE
where $\kappa$ is a parameter (often denoted $\varepsilon^2/2$) which is referred to as the {\em gradient coefficient}.\\
Then, the chemical potential $w$ is defined by:
\BE\label{Eq:0.4}
w := f'(u) - 2\kappa \Delta u.
\EE
$\Delta$ is the Laplace operator. If we denote by $J$ the flux and by $\mathcal{M}(u)$ the mobility, the classical Fick law provide the following equations:
\BE\label{Eq:0.5}
\partial_{t}u = -\nabla \cdot J \text{ and } J = -\mathcal{M}(u) \nabla w.
\EE
Finally, the Cahn-Hilliard equation takes the following general form:
\BE\label{Eq:0.6}
\left\{\begin{array}{ll}
\partial_{t} u = \nabla\cdot\left[\mathcal{M}(u) \nabla w\right],&\text{ on } \Omega\subset \R^d,\\
\\
w = \psi(u)-\varepsilon^2\Delta u ,&\text{ on } \Omega\subset \R^d,\\
\\
\nabla u \cdot \nu = 0 = \nabla w \cdot \nu, &\text{ on } \partial\Omega,\\
\end{array}\right.
\EE
where $t$ denotes the time variable, $\varepsilon$ ($=\sqrt{2\kappa}$) is a measure of the interfacial thickness, $\psi$ ($=f'$) is a nonlinear term, $\mathcal{M}$ is the mobility function, $\nu$ is the outward pointing unit normal on the boundary $\partial\Omega$.
It is well known that the Cahn-Hilliard equation is a gradient flow in $\sH^{-1}$ with
 Lyapunov energy functional $\mathcal{F}$.\par\smallskip
For a regular uniform alloy, the free energy $f$ is explicitly given by: 
\BE\label{Eq:0.7}
f : u \mapsto N_{m} k_{B} T_{c}\frac{1-u^2}{2} + N_{m}k_{B}T\left[\frac{1+u}{2}\ln\left(\frac{1+u}{2}\right)+\frac{1-u}{2}\ln\left(\frac{1-u}{2}\right)\right],
\EE
where $k_{B}$ is the Boltzmann constant, $N_{m}$ a molecular density, $T$ the temperature and $T_{c}> T$ the critical temperature. Thus the nonlinear term $\psi$ is:
\BE\label{Eq:0.8}
\psi:=f' : u \mapsto -N_{m} k_{B} T_{c}u + \frac{N_{m}k_{B}T}{2}\ln\left(\frac{1+u}{1-u}\right),
\EE
which is singular at $u=\pm1$. These singularities give rise to the first difficulty in a numerical study, so this function $\psi$ is often replaced by the derivative of the classic quartic double-well potential, where $f$ takes the following form:
\BE\label{Eq:0.9}
f : u \mapsto \frac{1}{4}\left(1-u^2\right)^2,
\EE
with derivative:
\BE\label{Eq:0.10}
\psi : u \mapsto u^3-u.
\EE
The Cahn-Hilliard equation has been extensively studied in the case where $\psi$ is replaced by a polynomial function (see \cite{CH1}, \cite{LANGER} and \cite{MR763473}). Furthermore, this model has been used successfully for describing phase separation phenomena, see for example the survey \cite{MR1657208}, and the references therein, or other recent results on spinodal decomposition and nucleation in \cite{MR1232163, MR2342011, MR1214868, MR1637817, MR1753703, MR1712442, MR1763320, MR2048517}. Recently, Ma and Wang have studied the stationary solutions of the Cahn-Hiliard equation (see \cite{MAWANG}).
The case of non smooth $\psi$ has been the object of much less research (see \cite{MR1123143} and \cite{MR1327930}).

Other frequent simplifications are often made. The mobility $\mathcal M$ is often assumed to be constant and the physical parameters are set to $1$ - as we have done above in \eqref{Eq:0.9}. For a more physically relevant choice of mobility, we mention~\cite{MR1300532} where the following form is proposed $\mathcal{M}(u)=\max\{0,1-u^2\}$.  Among the physical parameters, $\varepsilon$ has a peculiar role since it may lead to different asymptotic behaviors and equilibria (see \cite{MR1950337} and section \ref{S:3}). The study of evolution with $\varepsilon \rightarrow 0$ is of great importance: in particular a constant mobility leads to a Mullins-Sekerka evolution (nonlocal coupling) whereas a degenerate mobility leads to a purely local geometric motion (see \cite{MR1742748}). Furthermore, when the interface thickness is of the order of a nanometer, an artificially large parameter $\varepsilon$ is often used to regularize the numerical problem. When a fine resolution is out of reach, a change in the height of the barrier between wells in the free energy density, coupled with a change on $\varepsilon$, allows simulations with larger length scales (see \cite{MR2464502} for details).
\par\medskip


The evolution of the solution of~\eqref{Eq:0.6} can essentially be split into two stages. The first one is the \emph{spinodal decomposition} described in section \ref{S:2} where the two species quickly separate from each other. In longer time, the evolution is slower, and the solution tends to reduce its interfacial energy. These two evolutions require different methods for an efficient global simulation. In the beginning, a very small time step and a precise grid resolution allow efficient computation. But this is not appropriate to get long-time behaviors. So an adaptative time accurate or/and an adaptative mesh can improve the efficiency of the algorithms. However, in the long-time evolution, the interfaces have to be precisely captured so that a global adaptative mesh cannot be used. In the literature, many technical ideas have been studied: adaptive refinement of the time-stepping or of the mesh, $\mathcal{C}^1$ elements (see \cite{MR2464502}), multigrid resolution (see \cite{MR2183612}).\\

We propose here an alternative method using high degree $\mathcal{C}^0$ lagrangian nodal finite elements under a constant mobility $\mathcal{M}\equiv1$. The use of $p$-version (increasing polynomial degree, see \cite{MR615529}) instead of $h$-version (decrease mesh-step) has proved to be efficient for propagation \cite{MR2084226,MR1353516,MR1445739}, corner singularities \cite{MR947469}, or oscillating problems \cite{MR2340008}. The numerical results obtained here with the finite element library \textsc{M\'elina}~\cite{Melina} show that this method is suitable in the Cahn-Hilliard framework as well.

Our paper is organized as follows: in section \ref{S:1}, we shortly describe the discretization (in both time and space) including the nonlinear solver and the high degree finite elements we used. Section \ref{S:2} and section \ref{S:3} are respectively devoted to the numerical results for the one-dimensional and the two dimensional problem. We investigate the performance of our method through different quantitative and qualitative aspects of the Cahn-Hilliard equation: comparison to explicit profile-solution in 1D (see section \ref{S:2}), spinodal decomposition (see section \ref{S:2}), discussion about polynomial approximations of the logarithmic potential (see section \ref{S:2}), impact of the temperature and the parameter $\varepsilon$ (see section \ref{S:2} and \ref{S:3}), long-time behavior and asymptotic stable states (see section \ref{S:3}). The numerical results are compared with existing ones in the literature, validating our approach. 



\section{Discretization}\label{S:1}
\subsection{Space-Time schemes}\label{S:1.1}

We start with the description of the time discretisation. Given a large integer $N$, a time step $\tau$, and an initial data $(w_{0},u_{0})$, we denote by $(w_{n},u_{n})_{n \leq N}$ the sequence of approximations at uniformly spaced
times $t_n=n\tau$. The backward Euler scheme is given by:
\BE\label{Eq:1.1}
\left\{\begin{array}{ll}
\frac{u_{n+1}-u_{n}}{\tau} = \Delta w_{n+1},\\
\\
w_{n+1} = \psi(u_{n+1})-\varepsilon^2\Delta u_{n+1}.
\end{array}\right.
\EE
A Crank-Nicolson scheme could easily be implemented but our experiences show that it gives results quite similar to the ones we shall show in the sequel.
The schemes are immediately generalized to our case. We denote by $\langle\cdot,\cdot\rangle$ the scalar product in $\sL^2(\Omega)$. We use the standard Sobolev space $\sH^1(\Omega)$ 
equipped with the seminorm 
\BD
|h|_{1} =  \|\nabla h\|_{\sL^2},
\ED
and with the norm
\BD
\|h\|_{1} =  \left(|h|_{1}^2 + \|h\|^2_{\sL^2}\right)^{1/2}.
\ED
The weak form of the equation \eqref{Eq:1.1} reads:
\BE\label{Eq:1.3}
\left\{\begin{array}{ll}
\langle u_{n+1}-u_{n},\chi\rangle = -\tau\langle \mathcal{M}(u_{n+1}) \nabla w_{n+1},\nabla\chi\rangle, \text{ for all } \chi \in X_{1},\\
\\
\langle w_{n+1},\xi\rangle = \langle \psi(u_{n+1}),\xi\rangle+\langle\varepsilon^2\nabla u_{n+1},\nabla\xi\rangle, \text{ for all } \xi \in X_{2},
\end{array}\right.
\EE
where $X_{1}$ and $X_{2}$ are the spaces of test functions ($\sH^1(\Omega)$ for example). We discretise in space by continuous finite elements. Given a polygonal domain $\Omega$, for a small parameter $h>0$, we partition $\Omega$ into a set $\mathcal{T}^h$ of disjoint open elements $K$ such that $h=\max_{K \in \mathcal{T}^h} (\mathrm{diam}(K))$ and $\mathop{\bigcup}_{K \in \mathcal{T}^h}\overline{K} = \overline\Omega$.
Thus, we define the finite element space
\BE\label{Eq:1.4}
V^h = \left\{\chi \in \mathcal{C}(\bar{\Omega}) : \chi\big|_{K} \in \mathbb{P} \text{ for all } K \in \mathcal{T}^h\right\},
\EE
where $\mathbb{P}$ is a space of polynomial functions, see section \ref{S:1.3}. We denote by $(\varphi_{j})_{j\in J}$ the standard basis of nodal functions. Thus, for $u$ and $v \in \mathcal{C}(\overline{\Omega})$, we define the \emph{lumped scalar product}  by:
\BD
\langle u , v \rangle^h := \sum_{i,j} \langle u,\varphi_{i}\rangle\langle v,\varphi_{j}\rangle\langle \varphi_{i},\varphi_{j}\rangle.
\ED
The scheme \eqref{Eq:1.3} can be rewritten in the fully discrete form, just by replacing the continuous scalar product with the lumped scalar product.

We denote $\mathbf{u}= (u_{j})_{j\in J}$ and $\mathbf{w}= (w_{j})_{j\in J}$, the finite dimensional representation of $u$ and $w$ (we omit here the subscript $n$ of the time scheme). Then we define the matrices $\mathbf{A}$ and $\mathbf{M}$, whose coefficients are given by the following relations:
\BD
\begin{array}{rcll}
[\mathbf{A}]_{ij} &:=& \langle \nabla\varphi_{i}, \nabla\varphi_{j} \rangle,&\text{``stiffness'' matrix}, \text{ for all } i,j \in J,\\
\\
\left[\mathbf{M}\right]_{ij} &:=& \langle \varphi_{i}, \varphi_{j} \rangle,&\text{``mass'' matrix}, \text{ for all } i,j \in J.\\
\end{array}
\ED
For each time-step, given a previous solution $(\mathbf{w}_{n},\mathbf{u}_{n})$, $(\mathbf{w}_{n+1},\mathbf{u}_{n+1})$ is solution of the system
 \BE\label{Eq:1.5}
\left\{\begin{array}{llcl}
\tau\mathbf{A} \mathbf{w}_{n+1} &+ \mathbf{M} \mathbf{u}_{n+1}  &=& \mathbf{M} \mathbf{u}_{n},\qquad\qquad\quad\!\!\\
\\
\mathbf{M}\mathbf{w}_{n+1} &- \varepsilon^2 \mathbf{A}\mathbf{u}_{n+1}- \mathbf{M}\mathbf{\Psi}(\mathbf{u}_{n+1}) &=&0,
\end{array}\right.
\EE
where $\mathbf{\Psi}$ is a pointwise operator (\emph{related to $\psi$}), and with $(\mathbf{w}_{0},\mathbf{u}_{0})$ the finite dimensional representation of the initial data. The system \eqref{Eq:1.5} is clearly block-symmetric. The proof of the convergence of this scheme can be found in \cite{MR1609678}.



\subsection{Nonlinear solver}\label{S:1.2}
At each time step, we use a Newton procedure to solve the implicit nonlinear system~\eqref{Eq:1.5} 
. 
For \eqref{Eq:1.5}, we define the operator $\mathbf{L}$ by:
\BES
\mathbf{L} = \left(\begin{array}{cc}
\tau\mathbf{A} &\mathbf{M} \\
\mathbf{M} &- \varepsilon^2 \mathbf{A}
\end{array}\right).
\EES
Then denote by $\mathbf{S}$ the matrix of the left hand side of the backward Euler scheme,
\BES
\mathbf{S} = \left(\begin{array}{cc}
0 &\mathbf{M} \\
0&0
\end{array}\right).
\EES
Denote also by $\mathbf{G}$ the following operator:
\BES
G (\mathbf{w},\mathbf{u}) :=\left(\begin{array}{cc}
0\\-\mathbf{M}\mathbf{\Psi}(\mathbf{u})
 \end{array}\right).
\EES
Finally denote by $\mathbf{Y}_{n}$ the couple $(\mathbf{w}_{n},\mathbf{u}_{n})$ for each $n \leq N$.
The backward Euler scheme at each time-step satisfies the following formula:
\BE\label{Eq:1.7}
\mathbf{L} \mathbf{Y}_{n+1} + \mathbf{G} (\mathbf{Y}_{n+1}) - \mathbf{S} \mathbf{Y}_{n} = 0.
\EE
The Newton iterates $(\mathbf{Y}_{n}^k:=(\mathbf{w}_{n}^k,\mathbf{u}_{n}^k))_{k\in\N}$ satisfy for each $n\leq N$
\BE\label{Eq:1.8}
\left\{\begin{array}{lcl}
\mathbf{Y}_{n}^0 &=& \mathbf{Y}_{n},\\
\\
\mathbf{Y}_{n}^{k+1} &=& \mathbf{Y}_{n}^{k}  - \left[\mathbf{L}+D_{\mathbf{G}}\left(\mathbf{Y}_{n}^{k}\right)\right]^{-1} \left[\left(\mathbf{L}+\mathbf{G}-\mathbf{S}\right)\left(\mathbf{Y}_{n}^{k}\right)\right], \text{ for all } k \in \N,
\end{array}\right.
\EE
where $D_{\mathbf{G}}\left(\mathbf{Y}_{n}^{k}\right)$ is the differential of $\mathbf{G}$ at point $\mathbf{Y}_{n}^{k}$. Actually, we stop the procedure at $k=k_{n}$ when the residual is small, and define $\mathbf{Y}_{n+1}:= \mathbf{Y}_{n}^{k_{n}}$. System \eqref{Eq:1.8} is an implicit linear system for each Newton-step, handled with a biconjugate gradient method.\\
When the nonlinear term is logarithmic, we should deal with the singularities at $\pm 1$. However, in 
all our computations, the solution stays far from $\pm 1$ so that no special care is needed. This is expected. Indeed, it is known that in the one dimensional case the solution satisfies an $\rm{L}^\infty$ bound which
is strictly less than one (see \cite{MR1182511}). The same result has not been proved in higher dimension but it is probably true.\\
A simple remark shows the mass conservation through the total scheme. Indeed, if we multiply the first component of the second equation in the system \eqref{Eq:1.8} by the vector $\mathbb{I}:=(1,1,...,1)$ which belongs to $V^{h}$, we get for all $k \in \N$:
\BES
\mathbb{I}\,\mathbf{M}\mathbf{u}_{n}^{k} = \mathbb{I}\,\mathbf{M}\mathbf{u}_{n}.
\EES


\subsection{Implementation with high degree finite elements}\label{S:1.3}

The finite element library \textsc{M\'elina}~\cite{Melina} has the feature of providing lagrangian nodal elements with order up to $64$ (the nodes may be chosen as the Gauss-Lobatto points to avoid Runge phenomenon for large degrees). It can thus be used as a $p$-version code -- see~\cite{MR615529} -- or even to implement spectral methods -- see~\cite{MR1470226}. In the following results, we use quadrangular elements for two-dimensional computations, with degree from $1$ to $10$. So we use the notation $Q_{i}$ with $i\in\{1,2,3,4,5,6,7,8,9,10\}$ to describe these elements. We justify this strategy by the fact that the expected solution is smooth but may present a thin interface ; since high degree polynomials are able to capture high frequencies they are well suited in such situations. Some comparisons are shown below between degree $1$ on a refined mesh, and degree $10$ on a coarse mesh, justifying the efficiency of the method (in both terms of accuracy and computational cost). 


\section{Cahn-Hilliard evolution. Polynomial approximation of the logarithm}\label{S:2}

The temperature plays a crucial role in the evolution of the solution. The function $\psi$ defined in \eqref{Eq:0.8} depends on two 
values of the temperature $T$ and $T_{c}$.  
When the temperature $T$ is greater than the critical temperature $T_{c}$, the second derivative of $\psi$ is non-negative, 
thus function $\psi$ is convex and has only one minimum. We say that the function $\psi$ has a single well profile. 
Thus the solution tends to this unique minimum and the alloy exists in a single homogeneous state. 

\begin{figure}[htp]
\centering
\includegraphics[angle=0,width=12cm]{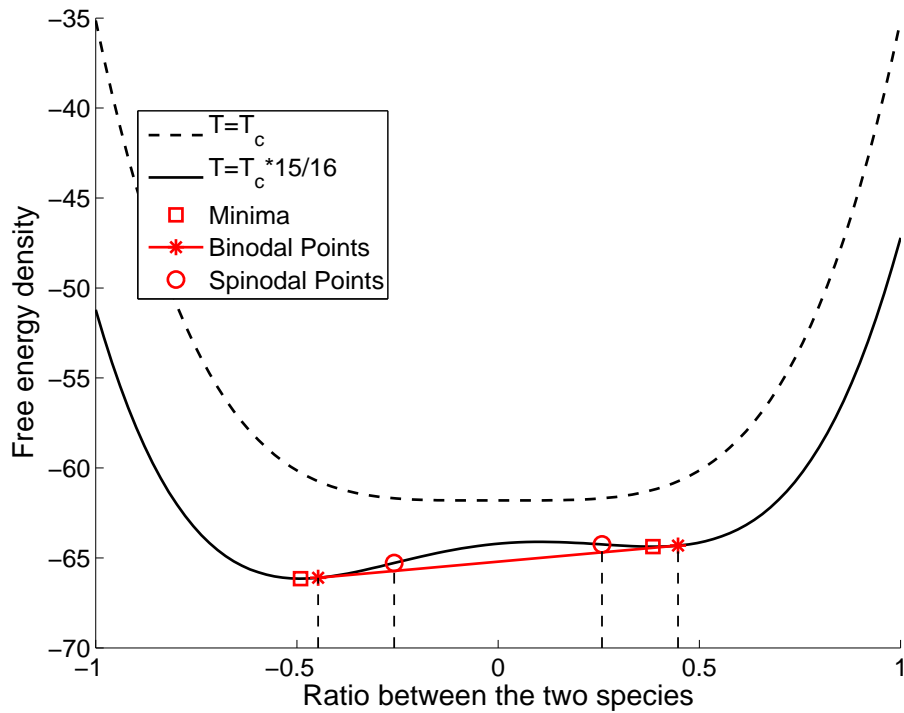}
\caption{Free energy density for two different temperatures.}\label{FIG:A}
\end{figure}

But when the temperature $T$ of the alloy is lowered under the critical temperature $T_{c}$, the function $\psi$ changes from a single well into a double well (see Figure \ref{FIG:A}), and the solution rapidly separates into two phases of nearly homogeneous concentration. This phenomenon is referred to as \emph{spinodal decomposition}. If the initial concentration belongs to the region where the energy density is concave, i.e. between the two \emph{spinodal points} $\sigma_-$ and $\sigma_+$ (see Figure~\ref{FIG:A}), the homogeneous state becomes unstable.

The concentrations of the two regions composing the mixture after a short stabilization have value near the so called \emph{binodal points} $\beta_-$ and $\beta_+$ (see also Figure~\ref{FIG:A}), defined by
\BE\label{Eq:2.6}
f'(\beta_{-}) = f'(\beta_{+}) = \frac{f(\beta_{+})-f(\beta_{-})}{\beta_{+}-\beta_{-}}, \quad \text{ with } \beta_{-} < \beta_{+}.
\EE
If the free energy is symmetric, the binodal points are the minima of each well, but in a more general case they are on a double tangent line (see \cite{MR2464502}). 

\begin{figure}[htp]
\centering
\subfigure[t=0]{
\label{FIG:E1.a}
\includegraphics[angle=0,width=6cm]{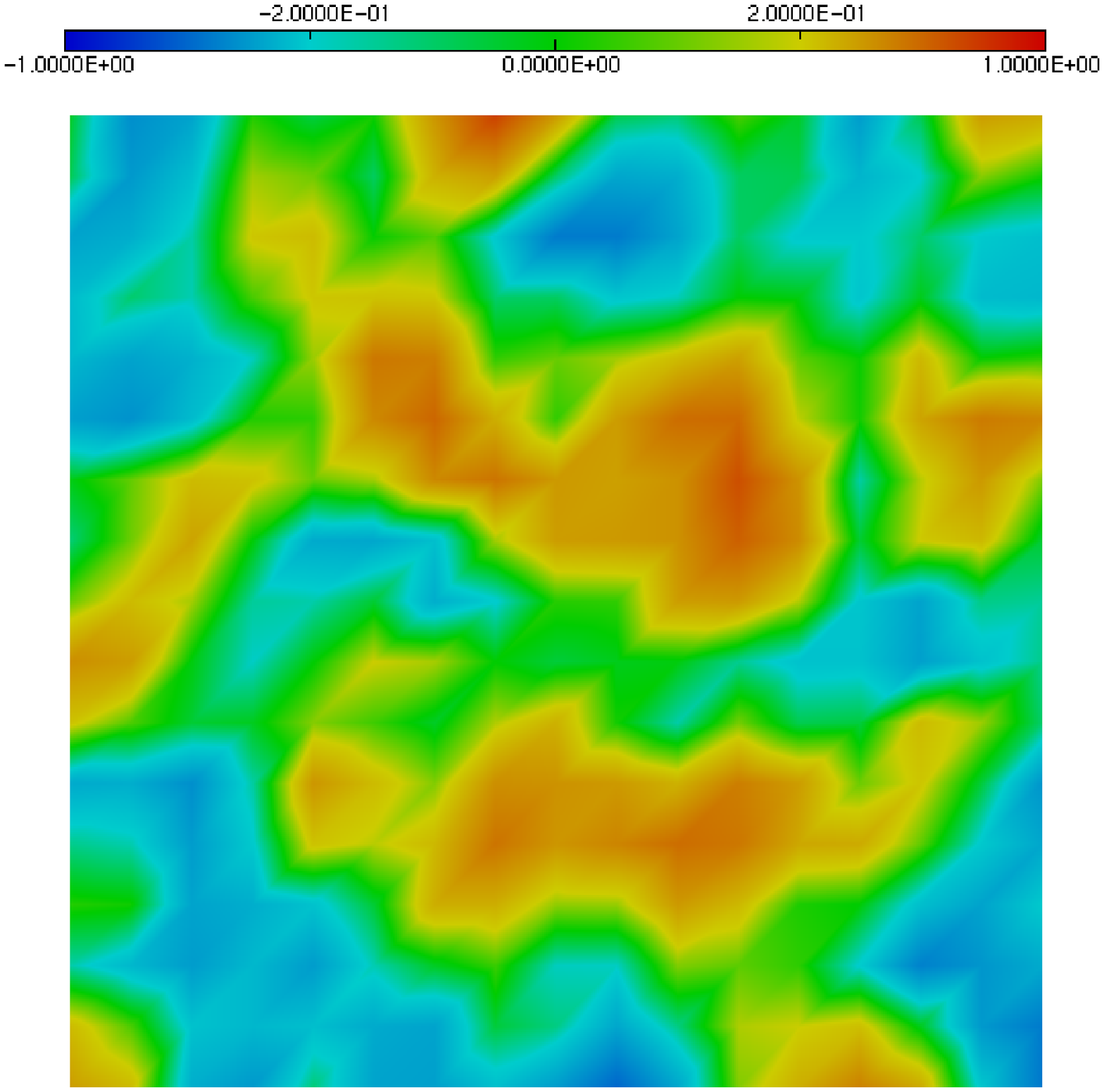}}
\hspace{0.3cm}
\subfigure[t=0.01]{
\label{FIG:E1.b}
\includegraphics[angle=0,width=6cm]{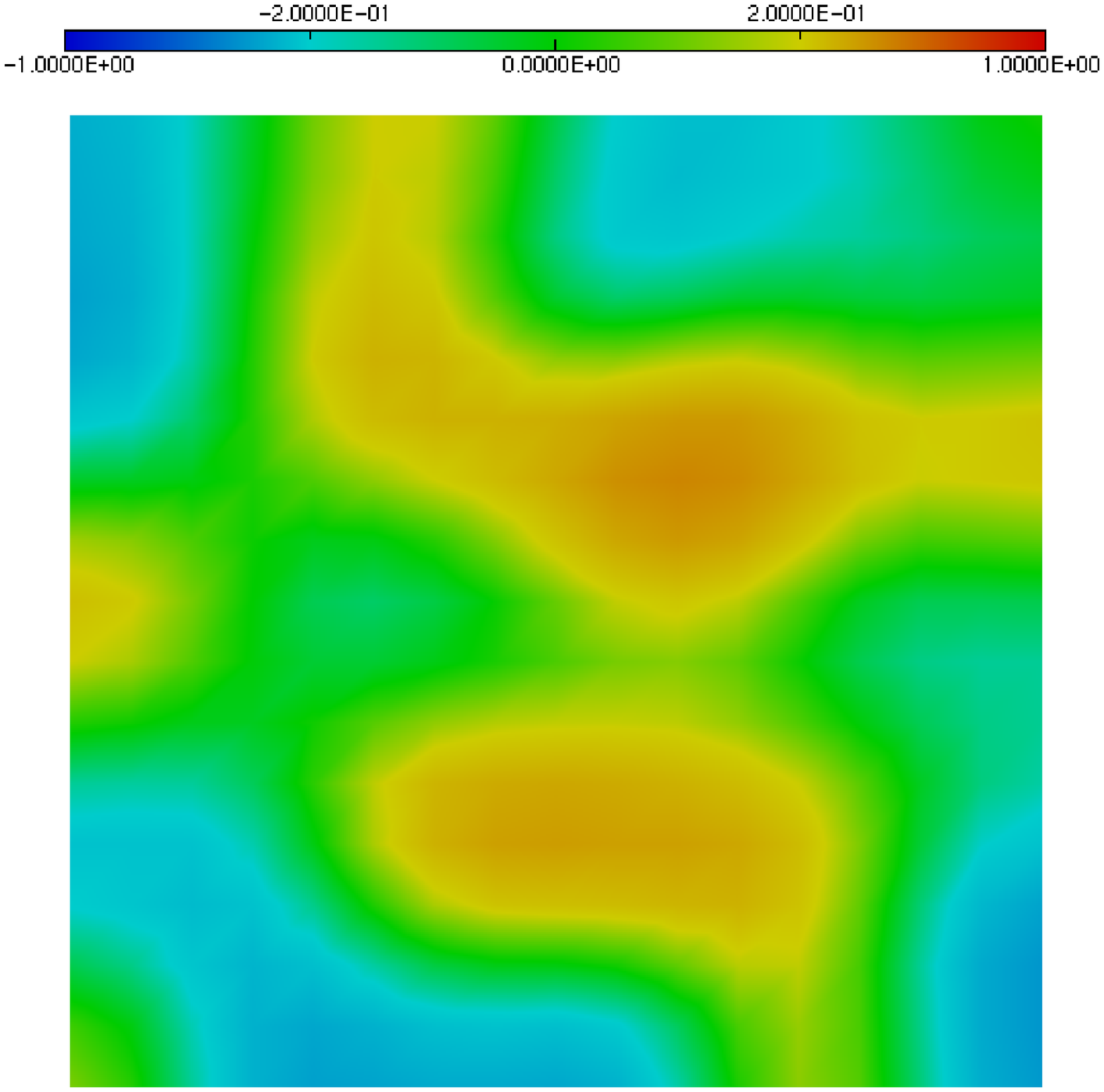}}
\\
\vspace{10pt}
\subfigure[t=0.05]{
\label{FIG:E1.c}
\includegraphics[angle=0,width=6cm]{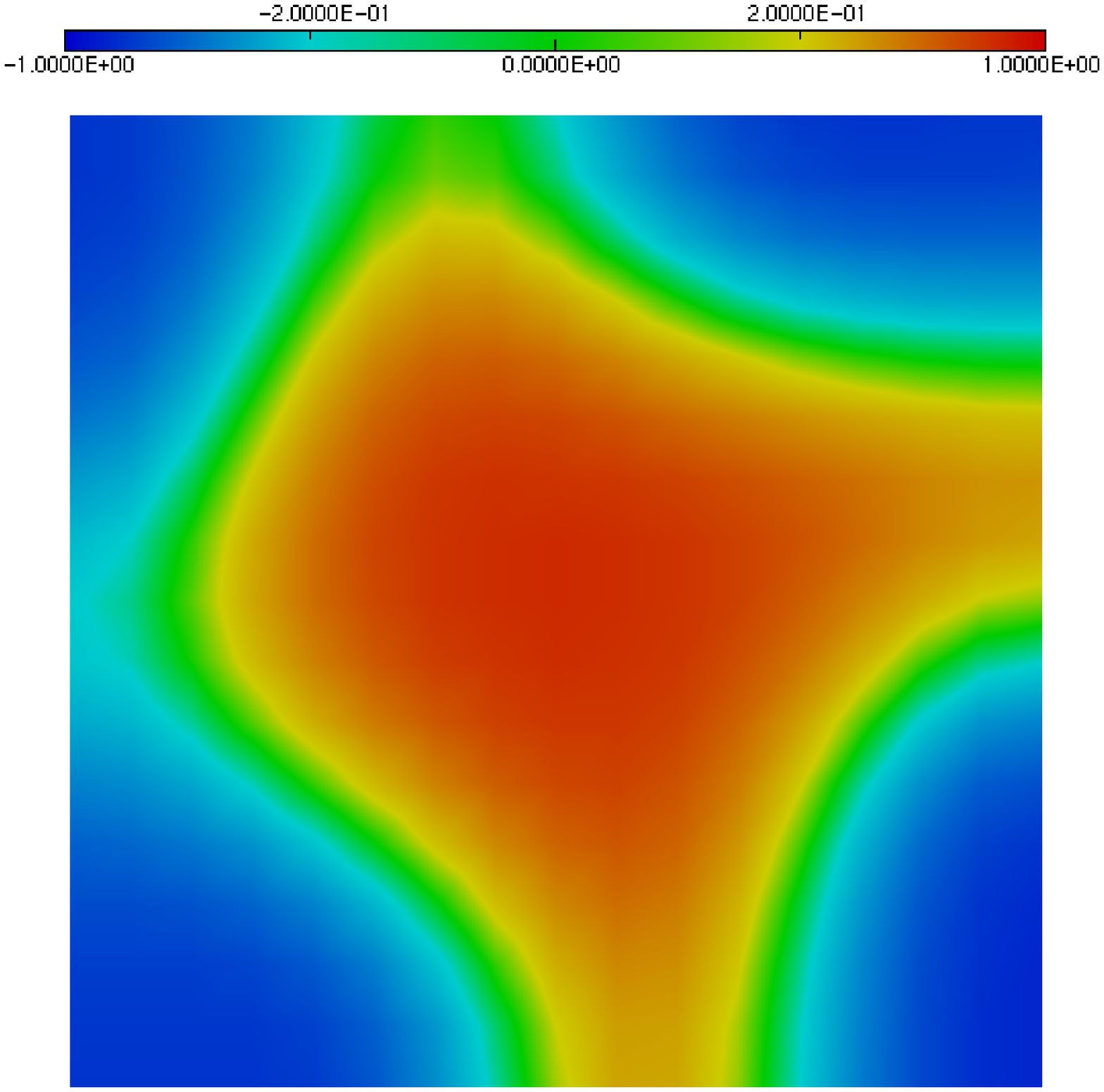}}
\hspace{0.3cm}
\subfigure[t=0.2]{
\label{FIG:E1.d}
\includegraphics[angle=0,width=6cm]{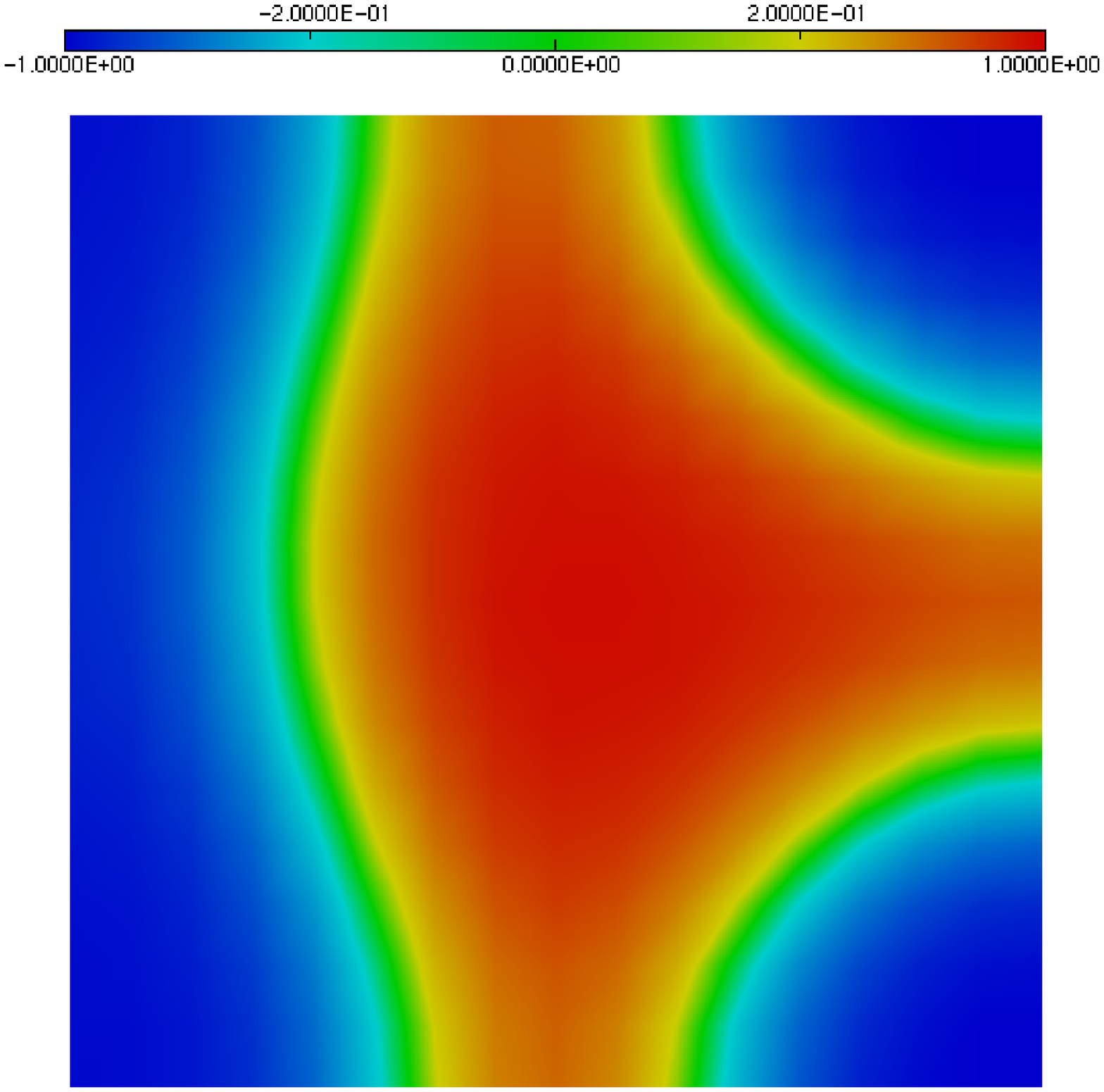}}
\\
\vspace{10pt}
\subfigure[t=0.6]{
\label{FIG:E1.e}
\includegraphics[angle=0,width=6cm]{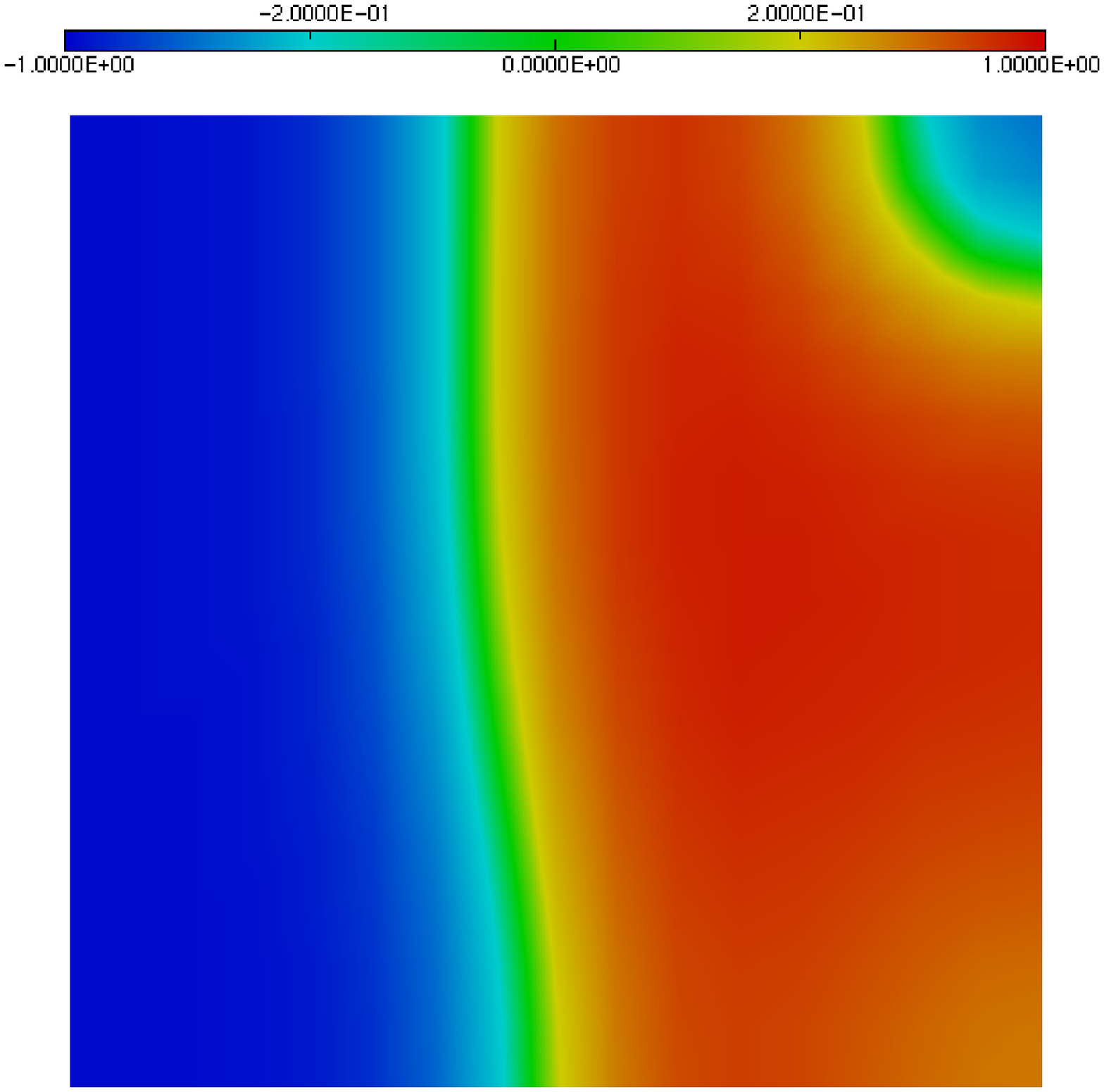}}
\hspace{0.3cm}
\subfigure[t=1]{
\label{FIG:E1.f}
\includegraphics[angle=0,width=6cm]{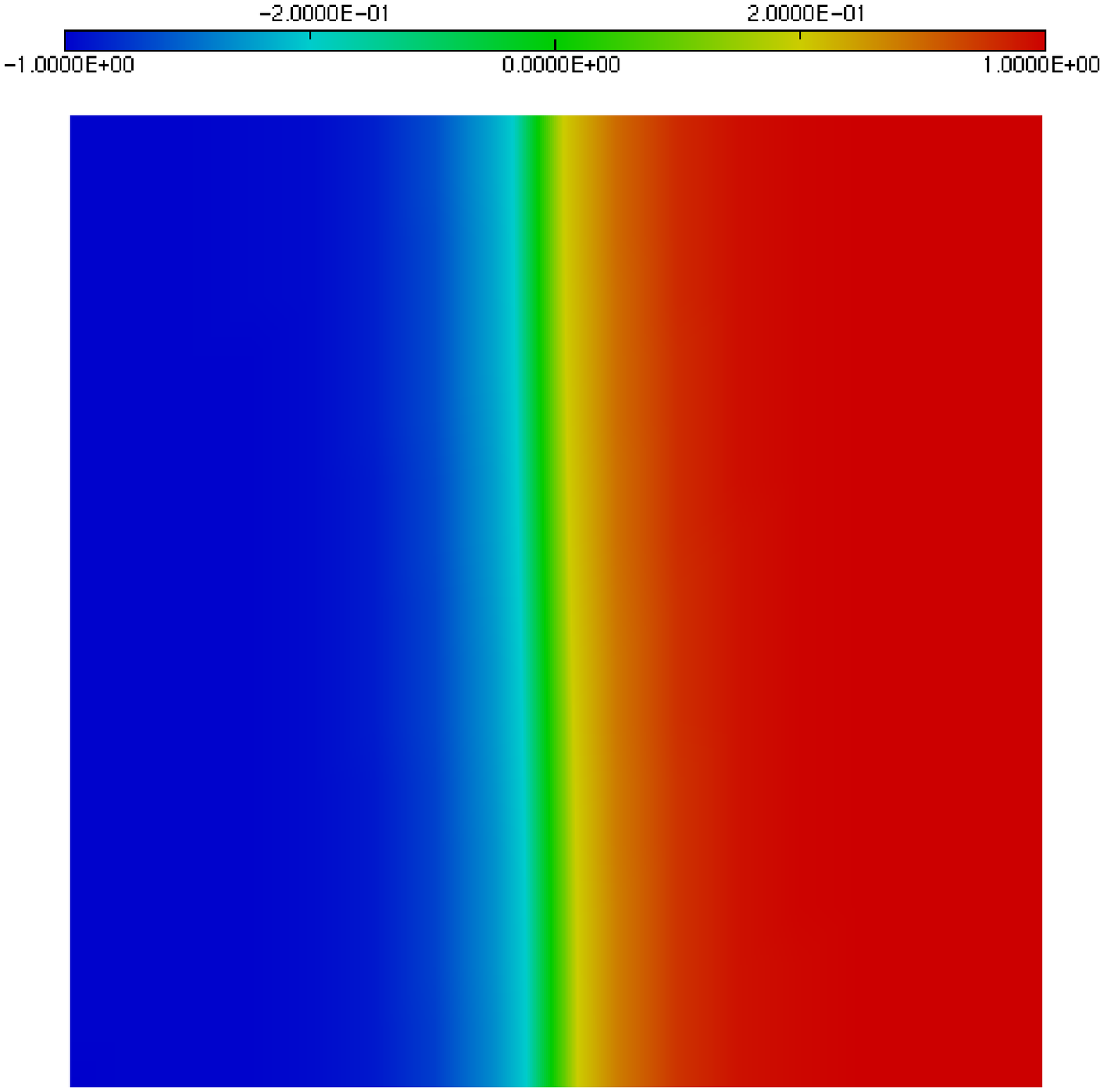}}
\caption{Spinodal decomposition under the classic quartic double-well potential.}
\label{FIG:E1}
\end{figure}

\begin{figure}[htp]
\centering
\subfigure[t=0]{
\label{FIG:E2.a}
\includegraphics[angle=0,width=6cm]{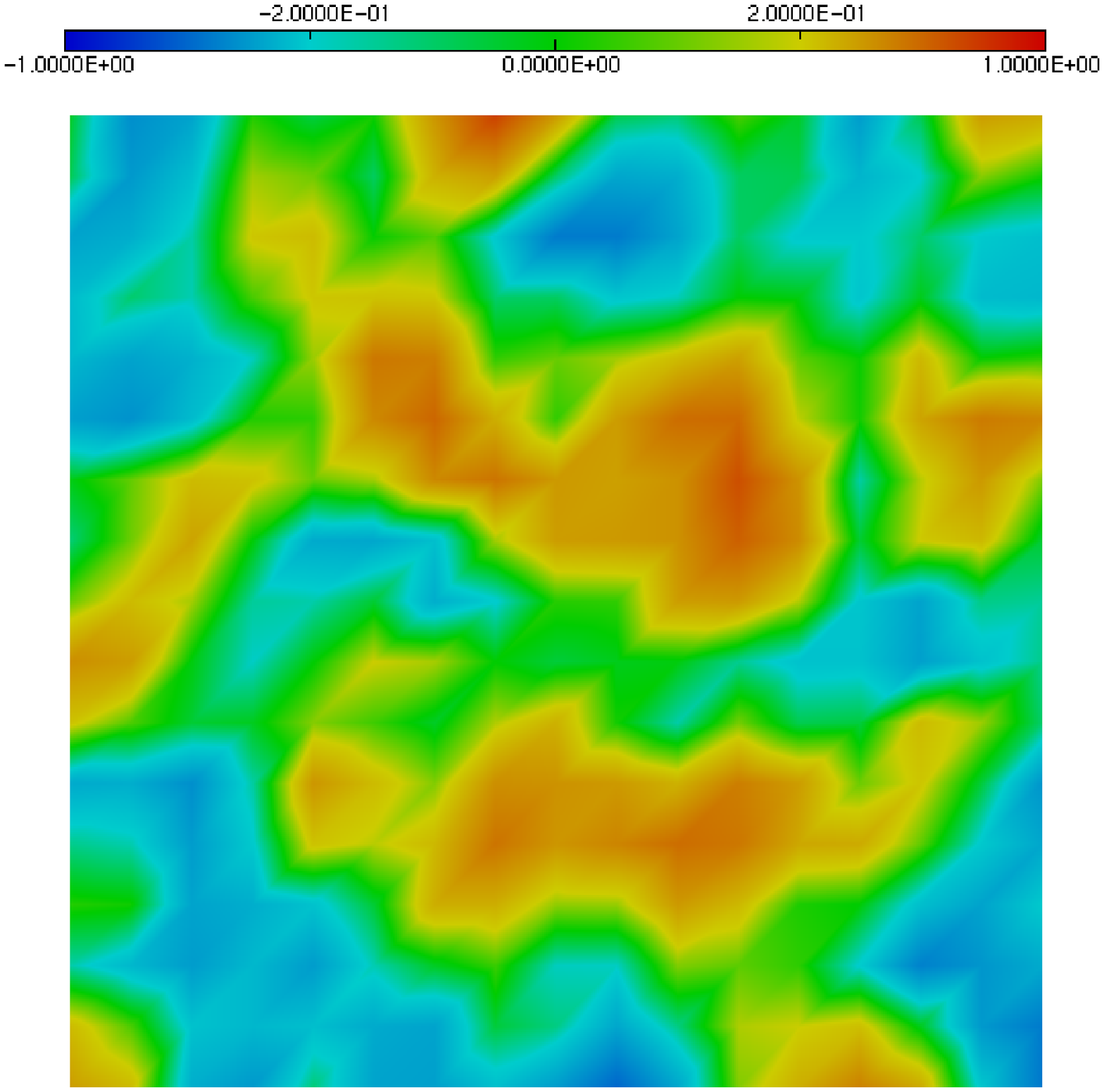}}
\hspace{0.3cm}
\subfigure[t=0.01]{
\label{FIG:E2.b}
\includegraphics[angle=0,width=6cm]{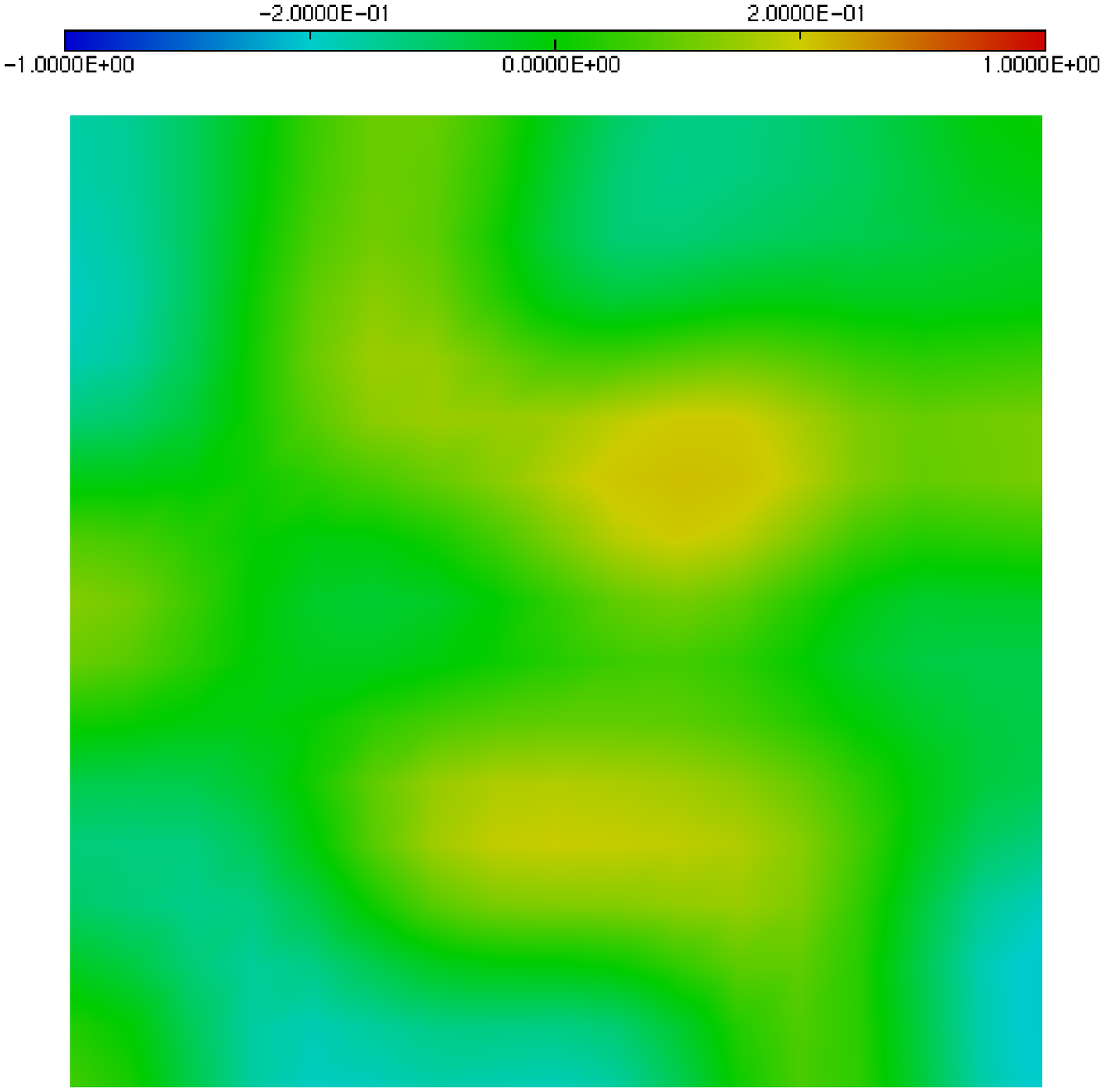}}
\\
\vspace{10pt}
\subfigure[t=0.05]{
\label{FIG:E2.c}
\includegraphics[angle=0,width=6cm]{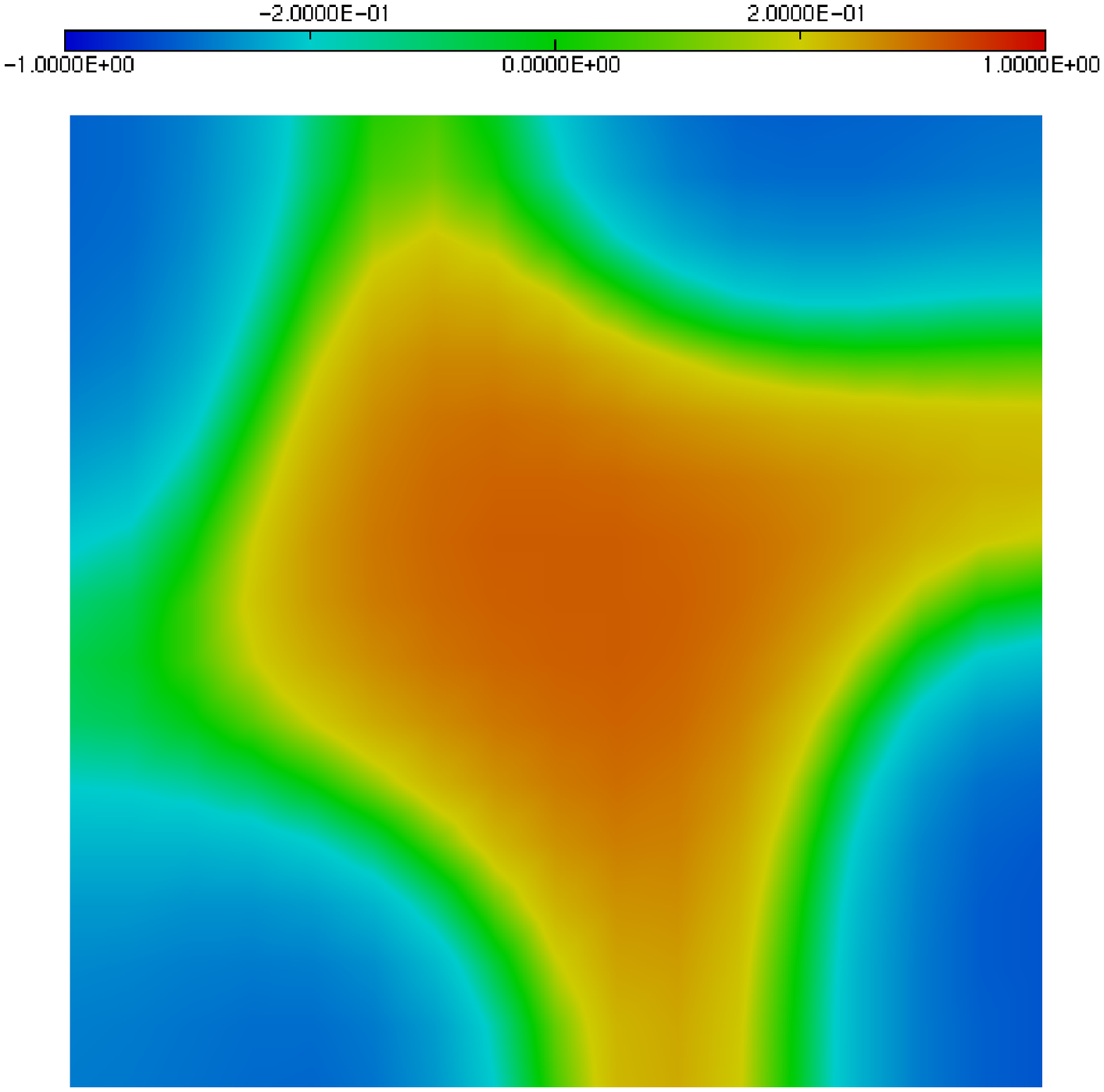}}
\hspace{0.3cm}
\subfigure[t=0.2]{
\label{FIG:E2.d}
\includegraphics[angle=0,width=6cm]{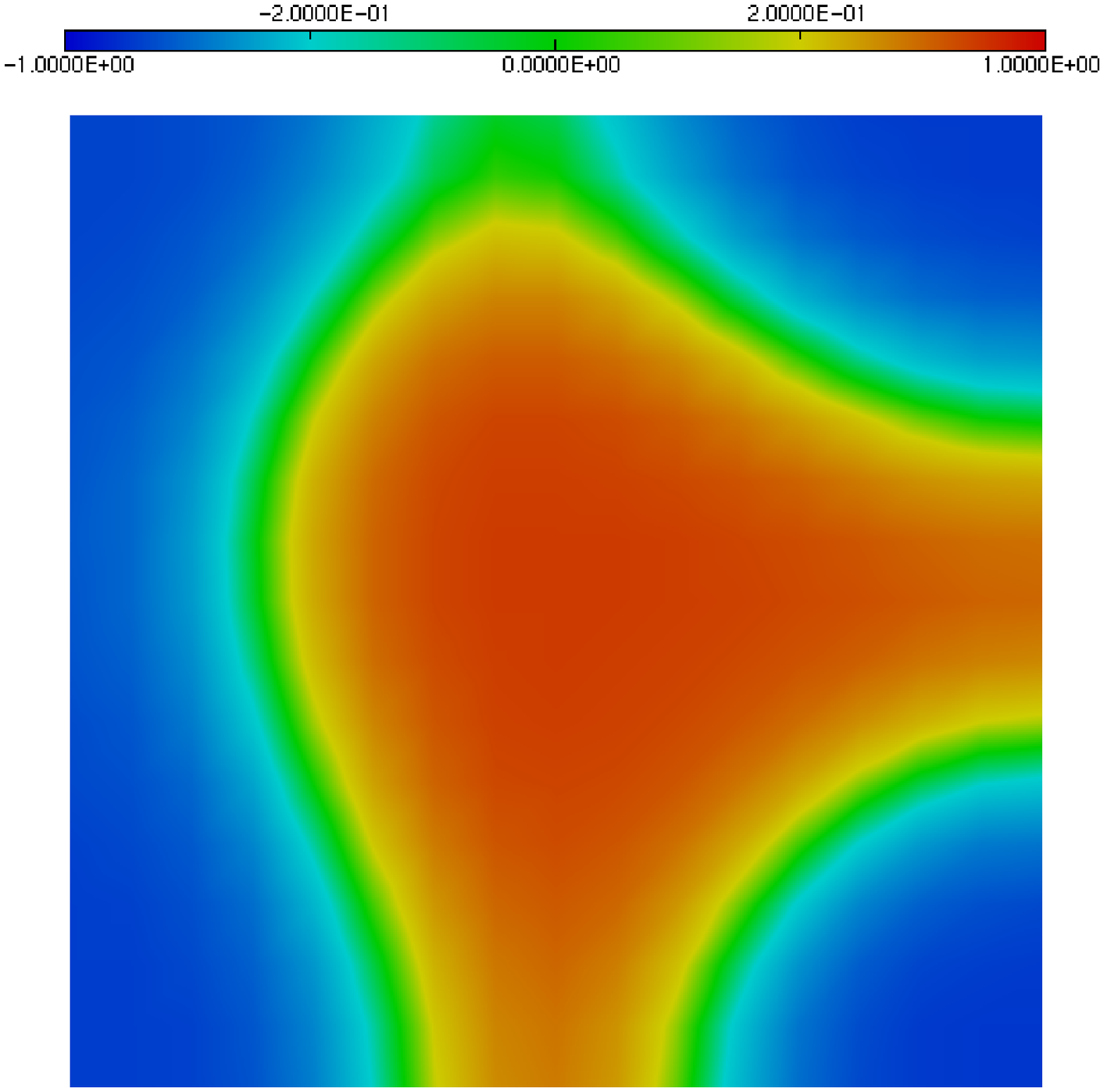}}
\\
\vspace{10pt}
\subfigure[t=0.6]{
\label{FIG:E2.e}
\includegraphics[angle=0,width=6cm]{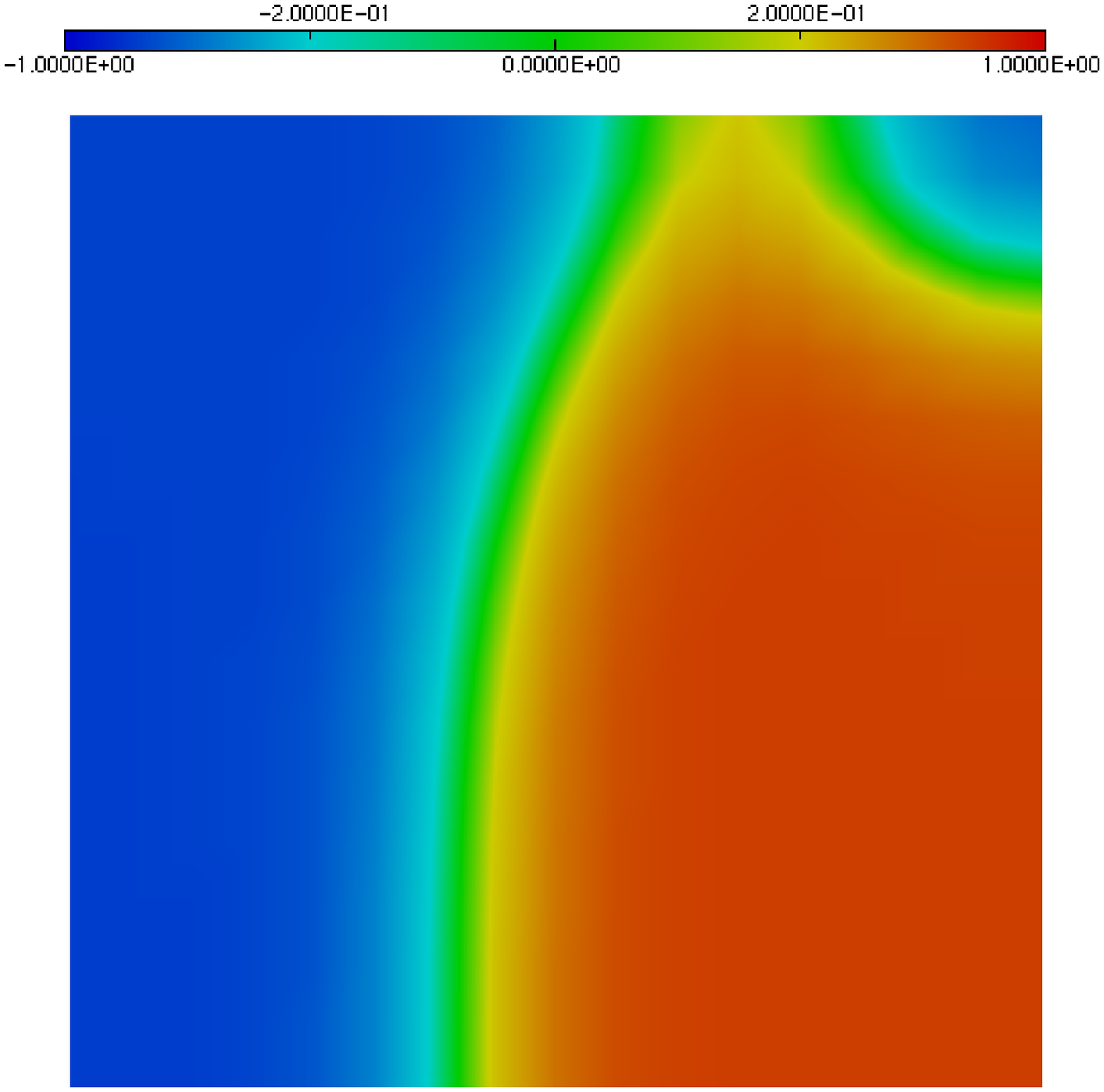}}
\hspace{0.3cm}
\subfigure[t=1]{
\label{FIG:E2.f}
\includegraphics[angle=0,width=6cm]{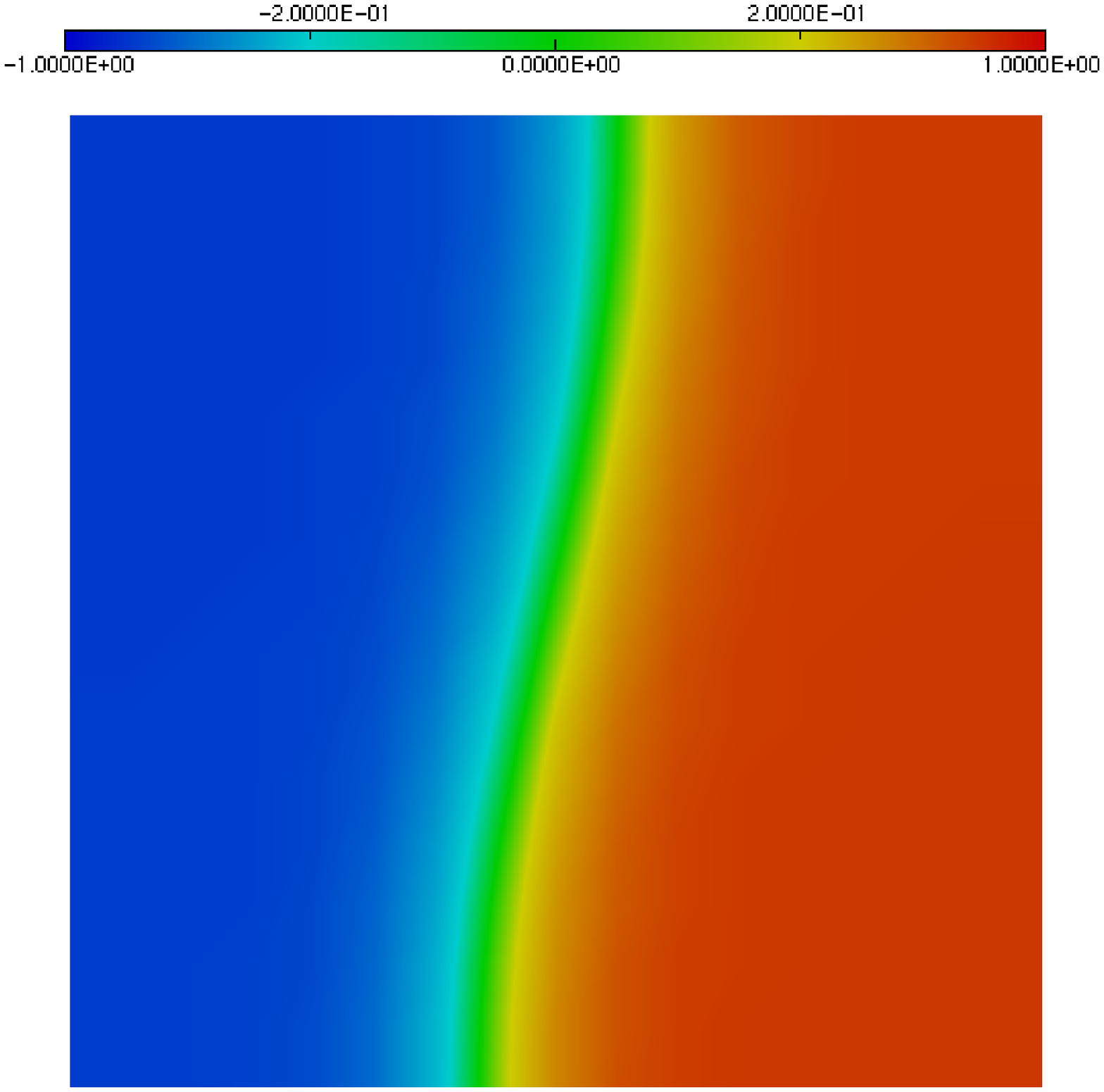}}
\caption{Spinodal decomposition under a logarithmic potential.}
\label{FIG:E2}
\end{figure}

The spinodal decomposition is represented in the first two graphs of Figure~\ref{FIG:E1} or Figure \ref{FIG:E2}. 

In longer time, the separated regions evolve to reduce their interfacial energies. These diffuse interfaces are shortened in an effect
resembling the surface tension on a sharp interface, as the material fronts move to reduce their own curvature (see \cite{MR1401172} and \cite{MR997638}). 
Finally, the solution reaches an equilibrium the location and form of which depend on the total initial concentration (see \cite{MR1950337}). Nevertheless this equilibrium is always a solution with an interface with minimal measure. On Figures \ref{FIG:E1} and \ref{FIG:E2}, this phenomenon is observed on the last four graphs.

Figure \ref{FIG:E1} corresponds to an evolution under the classic quartic double-well potential \eqref{Eq:0.9}
with non scaled coefficients, whereas Figure \ref{FIG:E2} corresponds to an evolution under the logarithmic potential \eqref{Eq:0.7}. They are both simulated on a $12\times12$ mesh under $Q_{1}$ polynomial elements. The $\varepsilon$ parameter is such that $\varepsilon^2 = 0.07$.
We see that the evolutions are quite similar and lead to the same stationary state. On these two evolutions, we can compare the difference of the energies or the $\rm{L}^2$ norm of the difference (see next paragraph).
Note that the polynomial approximation of the logarithm does not change the qualitative behavior. The 
same patterns appear and the long time behavior is very similar. The only notable difference is that with
the logarithmic nonlinearity, the dynamic is slower. This is particularly clear on the graphs (b). The spinodal 
decomposition is almost completed only for the polynomial. Similarly, on graphs (f), we see that at time
$t=1$, the logarithmic evolution has not reached equilibrium yet.
We have observed this in all our tests. 

\begin{figure}[htp]
\centering
\subfigure[Quartic potential - t=0]{
\label{FIG:D.a}
\includegraphics[angle=0,width=6cm]{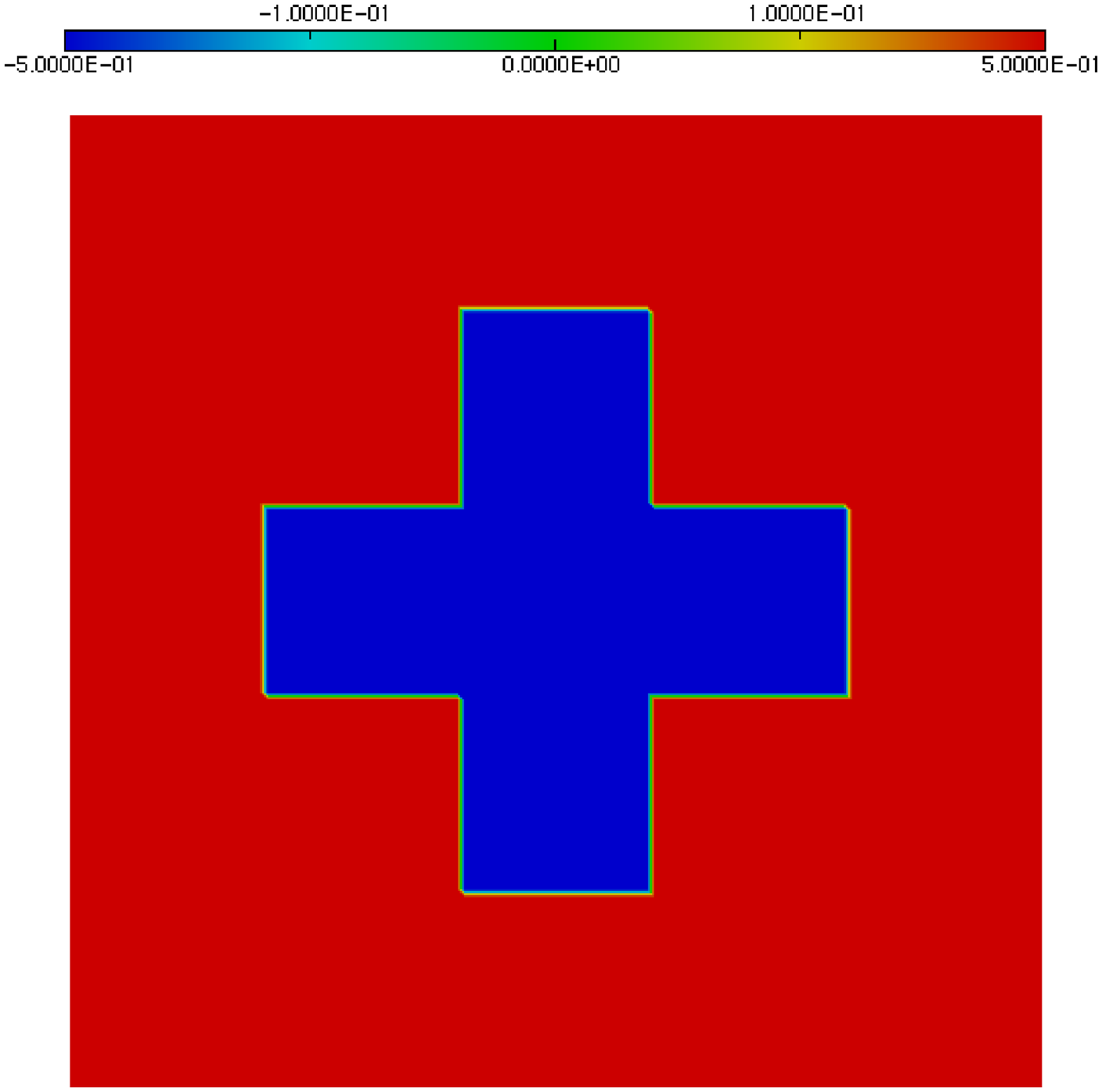}}
\hspace{0.3cm}
\subfigure[t=1]{
\label{FIG:D.b}
\includegraphics[angle=0,width=6cm]{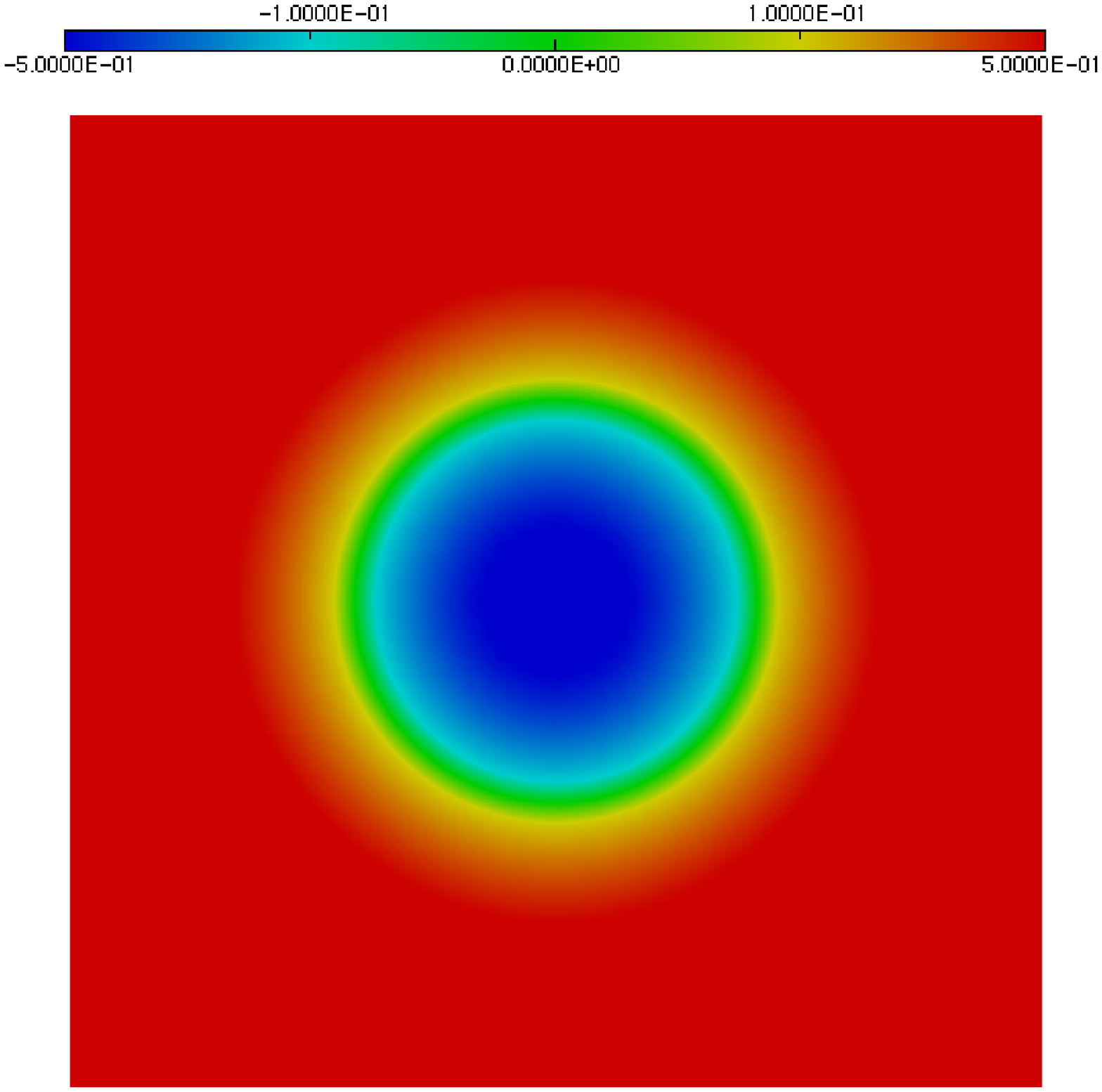}}
\\
\vspace{10pt}
\subfigure[Logarithmic potential - t=0]{
\label{FIG:D.c}
\includegraphics[angle=0,width=6cm]{\CROIXLOG}}
\hspace{0.3cm}
\subfigure[t=1]{
\label{FIG:D.d}
\includegraphics[angle=0,width=6cm]{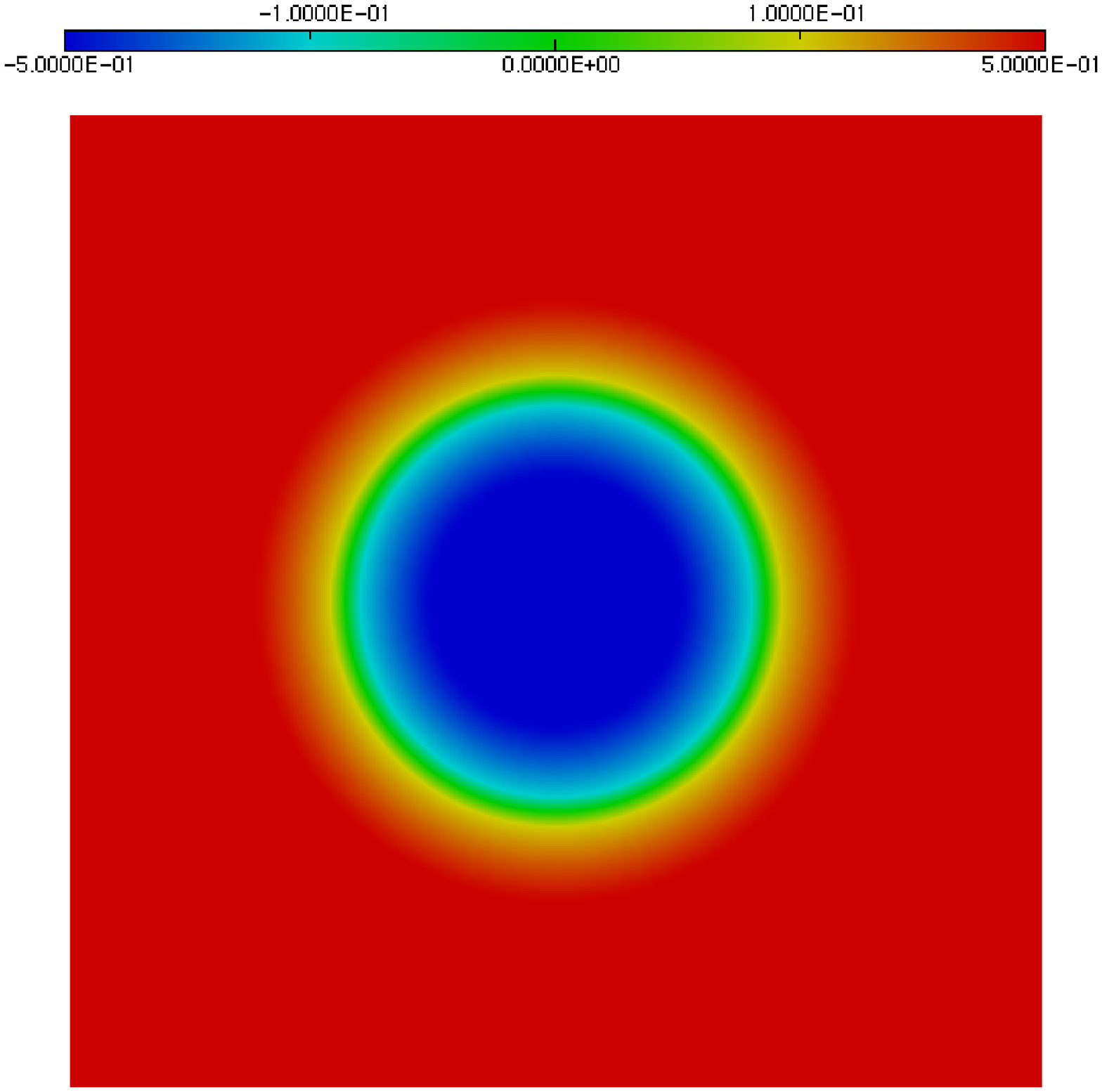}}
\caption{Evolution of a cross-shaped initial condition to a bubble.}
\label{FIG:D}
\end{figure}

The second evolution is often illustrated by the classical benchmark cross. It can be considered as a qualitative validation of the numerical methods. This long time behavior is illustrated in Figure~\ref{FIG:D}. Starting from a cross-shaped initial condition, the interface first diffuses from the arbitrary width of the initial condition to the equilibrium interface width. Next, the solution tries to reduce its interfacial energy and tends to a circular form. In the total free energy \eqref{Eq:0.3}, the term with the free energy function $f$ is responsible to the spinodal decomposition, whereas the gradient term is responsible for the interfacial reduction. This phenomenon has been simulated on a $256 \times 256$ mesh under $Q_{3}$ polynomial elements. Figures \ref{FIG:D} (a), (b) are obtained with
the quartic nonlinearity. We see on Figure \ref{FIG:D} (c) and (d) that again the qualitative behavior is 
very similar with the logarithm.

It is difficult to measure precisely the qualitative difference between the two evolutions. The only physical 
 quantity which can be measured in two dimensions is the energy. A detailed study of this aspect is 
 performed below. Moreover in the one-dimensional case, we are able to measure the interface. We will see that 
the quartic nonlinearity tends to thicken the interface. 

 The replacement of  the logarithmic free energy by the quartic one has been done by many authors in order to
avoid numerical and theoretical difficulties raised by the singular values $\pm1$. More generally, we can discuss the approximation of the logarithm by polynomial functions. We consider the $2n$-th order polynomial Taylor expansion $f_{2n}$:
\BE\label{Eq:2.7}
f_{2n} := u \mapsto \left(T_{c}\left(\frac{1-u^2}{2}\right) + T\left[-\ln(2) +\sum_{p=1}^{n} \frac{u^{2p}}{2p(2p-1)}\right]\right) + K_{2n}.
\EE
It is defined up to an additive constant $K_{2n}$. The constant $K_{2n}$ is apparently arbitrary. However,
it is preferable to choose it in order that the energy  of a solution $u$
\BE
\mathcal{F}_{2n} (u) := \int_{\Omega} \left( f_{2n}(u) + \kappa |\nabla u|^2 \right)\dd V
\EE
is well defined on unbounded domains. Since it is expected that the solution converges to one
of the binodal values, it is natural to choose $K_{2n}$ so that $f_{2n}$ vanishes at those points. We always consider this choice.

We have seen above that the quartic approximation does not seem to change drastically the 
qualitative behaviour, except that the evolution is faster. We now perform a quantitative study
to measure more precisely the effect of the polynomial approximation.

\begin{figure}[htp]
\centering
\begin{psfrags}
\psfrag{Title}{}
\psfrag{Xlabel}{{\Large \!\!\!\!\!\!\!\!\!\!\!\!\!\!\!\!\!\!\!\!\! Polynomial degree}}
\psfrag{Ylabel}{{\Large \!\!\!\!\!\!\!\!\!\!\!\!\!\!\!\!\!\!\!\!\!\!\!\!\!\!\! Error with the exact point}}
\includegraphics[angle=0,width=12cm]{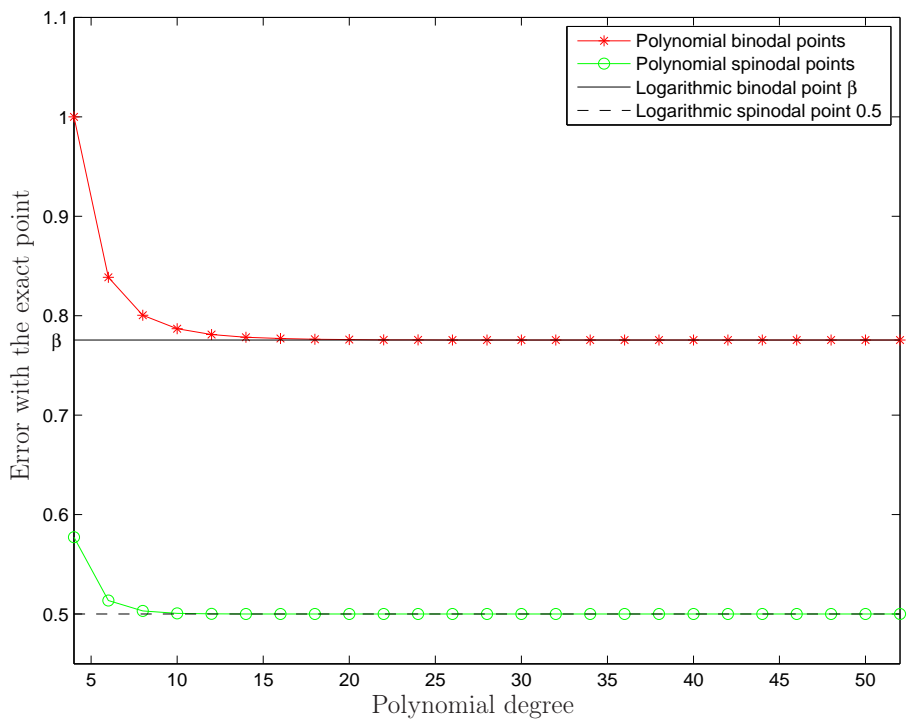}
\caption{Polynomial and logarithmic spinodal and binodal points.}\label{FIG:B}
\end{psfrags}
\end{figure}

\begin{figure}[htp]
\centering
\begin{psfrags}
\psfrag{Title}{}
\psfrag{Xlabel}{{\Large \!\!\!\!\!\!\!\!\!\!\!\!\!\!\!\!\!\!\!\!\! Polynomial degree}}
\psfrag{Ylabel}{{\Large \!\!\!\!\!\!\!\!\!\!\!\!\!\!\!\!\!\!\!\!\!\!\!\!\!\!\!\!\!\!\!\!\!\!\!\!\!\!\!\!\!\!\!\!\!\! Logarithmic error with the exact point}}
\includegraphics[angle=0,width=12cm]{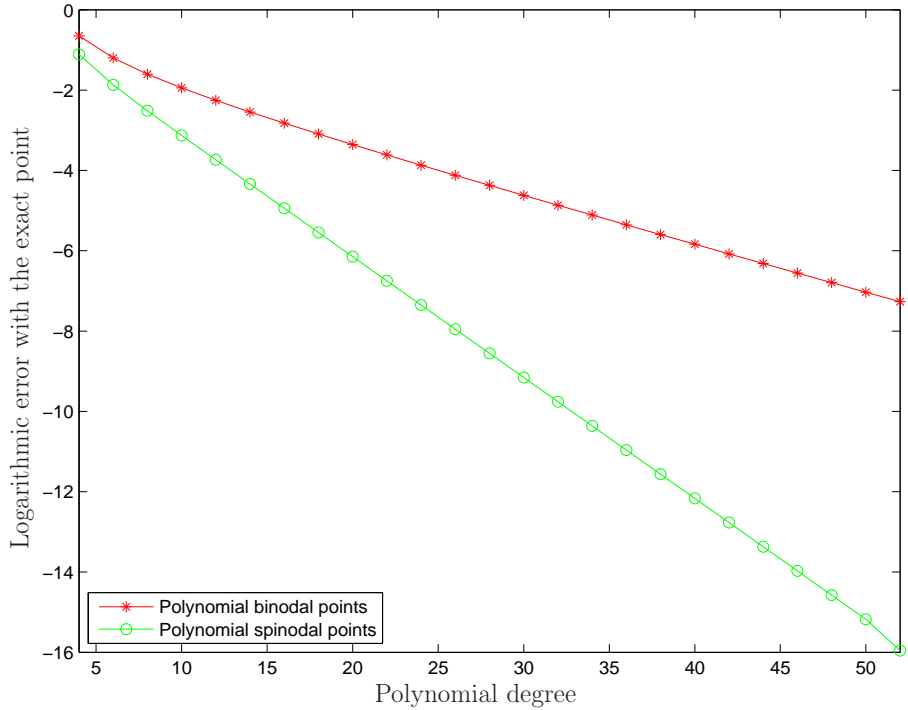}
\caption{Rate of convergence of the polynomial points.}\label{FIG:C}
\end{psfrags}
\end{figure}

The spinodal and binodal points are drawn in Figure \ref{FIG:B} for various $n$. When $n$ increases, the spinodal and binodal points converge to the corresponding values for the logarithmic potential. However, the convergence is rather slow (see Figures \ref{FIG:B} and \ref{FIG:C}).

\begin{figure}[htp]
\centering
\begin{psfrags}
\psfrag{-u+}{{\huge$-u_{+}$}}
\psfrag{u+}{{\huge$u_{+}$}}
\psfrag{Slope mu}{{\LARGE\color{blue}Slope $\mu$}}
\psfrag{l}{{\LARGE\color{blue}$\ell$}}
\includegraphics[angle=0,width=10cm]{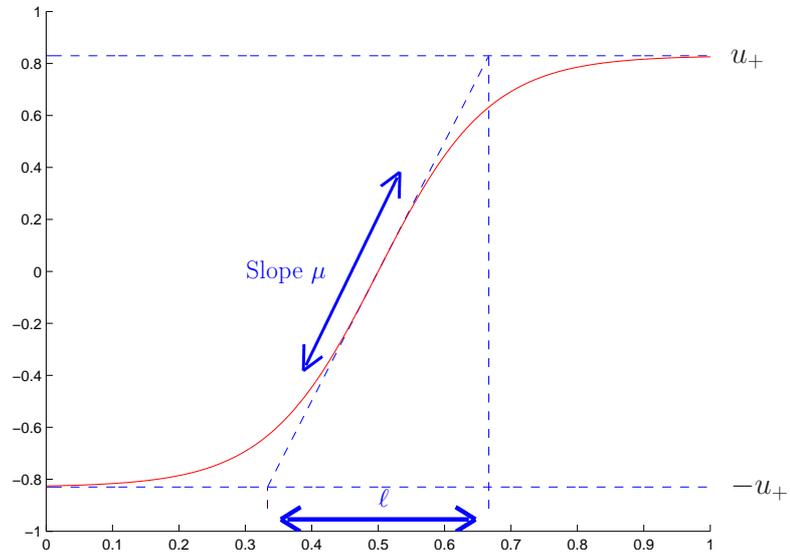} 
\caption{Interface length for the solution $u_{4}$.}\label{FIG:F}
\end{psfrags}
\end{figure}

In the one-dimensional case, it is possible to study the thickness of the interface.
Let us consider the domain $\Omega=\mathbb R$ and the quartic potential
\BE\label{Eq:2.1}
\psi_{4} := u \mapsto - T_{c}u + T\left(u+\frac{u^3}{3}\right),
\EE
which is the derivative of
\BE\label{Eq:2.2}
f_{4} := u \mapsto T_{c}\left(\frac{1-u^2}{2}\right) + T\left[\frac{u^2}{2} +\frac{u^4}{12}\right]+ K_{4}.
\EE
Then a stationary solution of the Cahn-Hilliard equation~\eqref{Eq:0.6} can be explicitly computed (under a constant mobility $\mathcal{M}(u)\equiv1$), see~\cite{MR983721}:
\BE\label{Eq:2.3}
u_{} : x\mapsto u_{+}\tanh\left( x\mu\right),
\EE
where
\BE\label{Eq:2.4}
u_{+} = \sqrt{3\left(\frac{T_c}{T}-1\right)} \quad\text{ and }\quad \mu = \frac{\sqrt{T_c-T}}{\varepsilon\sqrt{2}}.
\EE
It is important to remark that the solution is constrained in $[-u_{+},u_{+}]$. 
We can define a characteristic length $\ell$ (see Figure \ref{FIG:F}), corresponding to the width of the region containing the main variations of a solution $u$ :
\BD
\ell := \frac{|\lim_{x\rightarrow +\infty} u(x)|+|\lim_{x\rightarrow -\infty} u(x)|}{\text{Slope in interface point}},
\ED
where the interface point is the point $x_{0}$ where $u(x_{0})=0$. Thus we can compute explicitly this length and obtain:
\BE\label{Eq:ell}
\ell =\frac{2u_+}{u'_{4}(0)} =\frac{2\varepsilon\sqrt{2}}{\sqrt{T_{c}-T}}.
\EE
Cahn and Hilliard have defined a parameter $\lambda:=\frac{2\varepsilon\sqrt{2}}{\sqrt{T_{c}}}$ in order to characterize the interface length. With this parameter $\lambda$ we obtain the following expression for $\ell$:
\BE\label{Eq:2.5}
\ell = \frac{\lambda}{\sqrt{1-\frac{T}{T_{c}}}}.
\EE
Cahn and Hilliard have shown that in the case of the logarithmic free density the interface length 
is of the same order. This suggests that the quartic double well approximation preserves important features of the solution.

\begin{figure}[htp]
\centering
\includegraphics[angle=0,width=12cm]{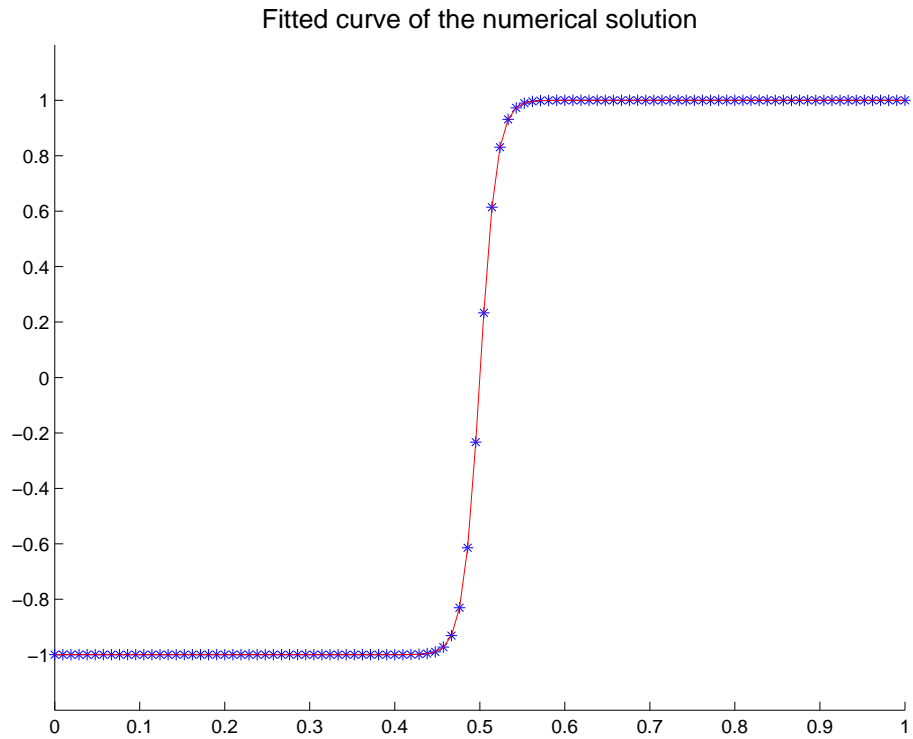}
\caption{Fitted curve on the ``tanh-profile''.}
\label{FIG:G}
\end{figure}

In Figure~\ref{FIG:G}, we present the numerical solution for $\Omega=[0,1]$ (blue stars), and the ``tanh-profile'' whose coefficients $u_+$ and $\mu$ have been fitted to the data. The fitting on $u_{+}$ corresponds to the value of the solution on the boundaries of the domain $\Omega$. And the fitting on $\mu$ corresponds to a least square method between the numerical solution and a "tanh-profile" solution interpolated on the same meshes. The ``tanh-profile'' (defined over $\mathbb R$) may be considered as a good approximation of the solution on $\Omega =[0,1]$ since the interface is very thin.
The numerical solution is computed with 35 $Q_{3}$-elements.

\begin{figure}[htp]
\centering
\includegraphics[angle=0,width=12cm]{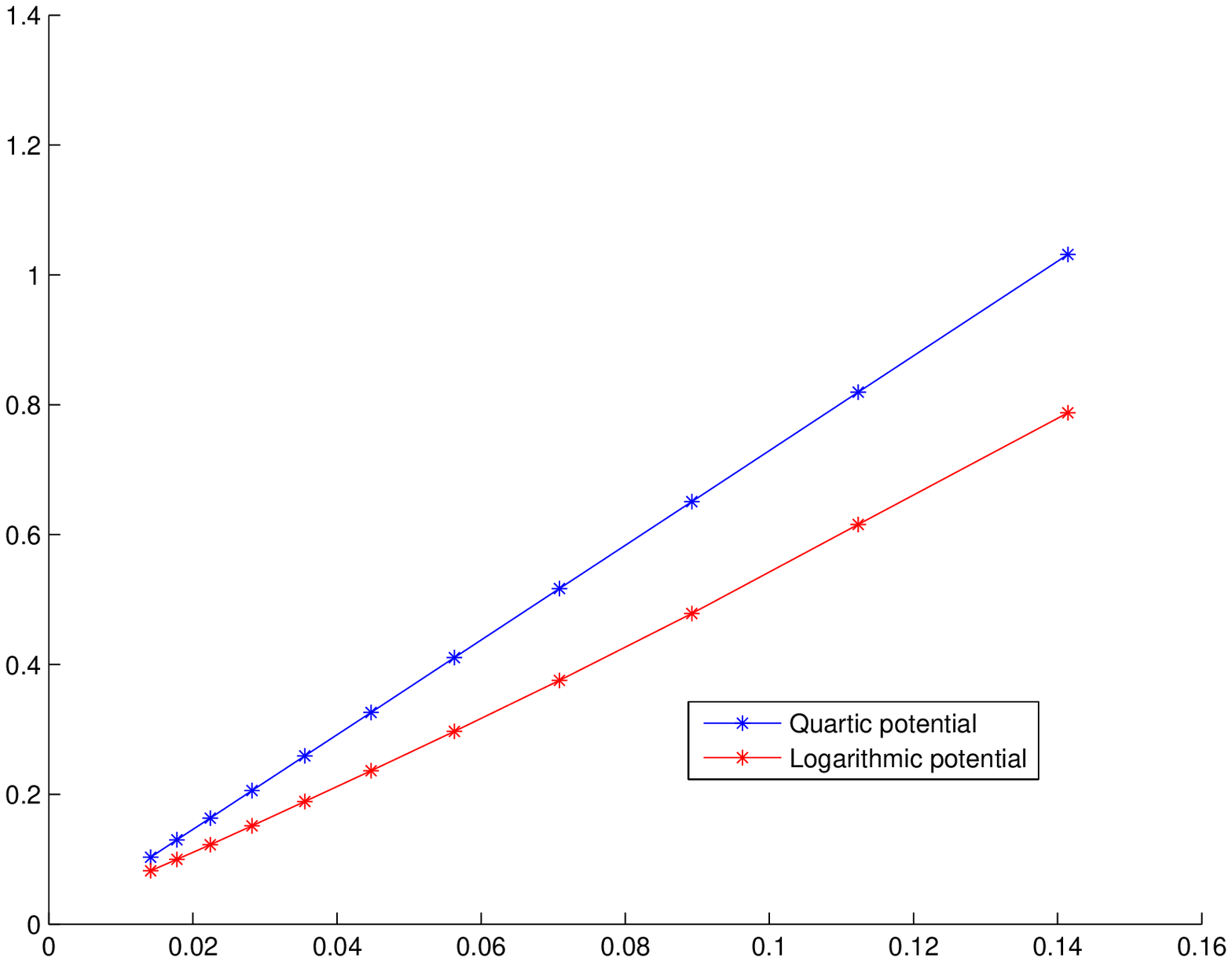}
\caption{Length of the interface for the quartic and logarithmic potentials.}
\label{FIG:GG}
\end{figure}

However, we have measured numerically the interface width in the quartic and logarithmic cases. 
This width is plotted for various $\varepsilon$ on Figure \ref{FIG:GG}. We see that as expected by the formula  \eqref{Eq:ell}, 
it varies linearly with $\varepsilon$. But, for $\varepsilon$ not too
small, the interface width is thinner for the logarithmic equation. The quartic approximation 
introduces a non negligible extra diffusivity.

We can also compare the total free energies. Denote by $u$ the solution of a simulation with the \emph{logarithmic} function $f$ and by 
$(u_{2n})_{n\geq2}$  the family of solutions of the simulations with the polynomial functions $(f_{2n})_{n\geq2}$. For the energy, we take as reference the logarithmic total free energy, and we study
\BE
|\mathcal{F}(u_{2n}) - \mathcal{F}(u)|.
\EE

\begin{figure}[htp]
\centering
\begin{psfrags}
\psfrag{F4}{{\Large$f_{4}$}}
\psfrag{F6}{{\Large$f_{6}$}}
\psfrag{F8}{{\Large$f_{8}$}}
\psfrag{F10}{{\Large$f_{10}$}}
\psfrag{F12}{{\Large$f_{12}$}}
\psfrag{F14}{{\Large$f_{14}$}}
\psfrag{F16}{{\Large$f_{16}$}}
\includegraphics[angle=0,width=12cm]{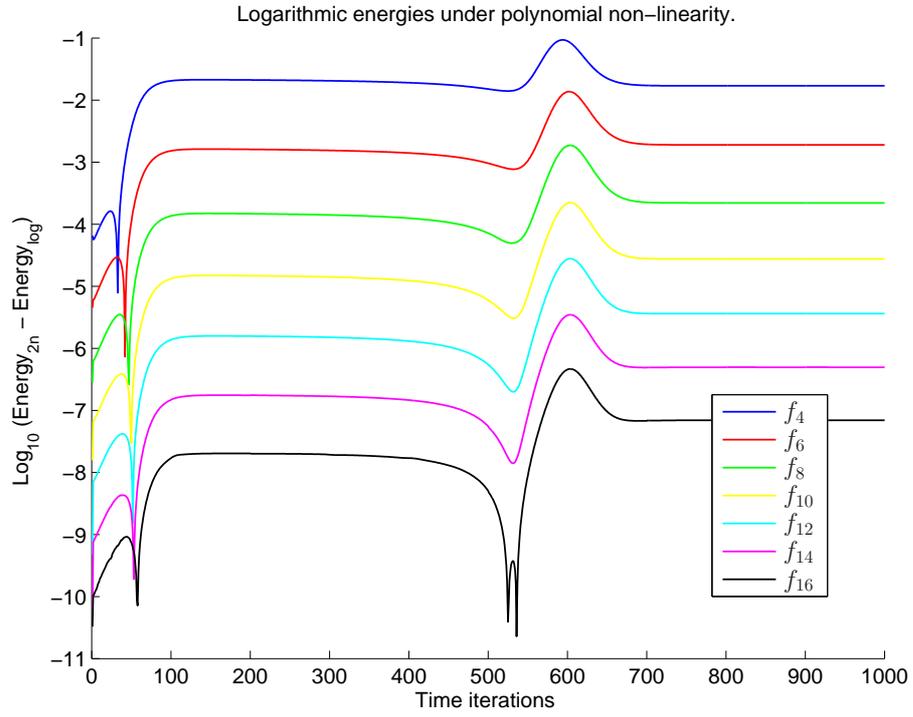}
\end{psfrags}
\caption{Polynomial energies versus logarithmic energy.}\label{FIG:H}
\end{figure}

On Figure~\ref{FIG:H}, the evolution of the logarithm of this 
quantity is plotted during a classical spinodal decomposition in dimension one.

We can see important peaks at the begining and smoother peaks between iterations $500$ and $700$. These peaks appear when the solution has a rapid evolution and when its topological form changes. For instance, these peaks correspond to the changes beetween the fourth and the fifth images of  Figure~\ref{FIG:E1}, and between the fifth and the sixth images. After the iteration $750$, all the solutions are in an asymptotic stable state, and the energies do not change anymore. \par
For a quartic potential ($n=2$ i.e. $f_{4}$ in Figure \ref{FIG:H}), the energy error is significant and 
the polynomial approximation is not good in that respect. 

We could as well have shown the evolution of 
$$
|\mathcal{F}_{2n}(u_{2n}) - \mathcal{F}(u)|
$$
In fact, it is very similar and does not bring new information.

\begin{figure}[htp]
\centering
\begin{psfrags}
\psfrag{F4}{{\Large$f_{4}$}}
\psfrag{F6}{{\Large$f_{6}$}}
\psfrag{F8}{{\Large$f_{8}$}}
\psfrag{F10}{{\Large$f_{10}$}}
\psfrag{F12}{{\Large$f_{12}$}}
\psfrag{F14}{{\Large$f_{14}$}}
\psfrag{F16}{{\Large$f_{16}$}}
\includegraphics[angle=0,width=12cm]{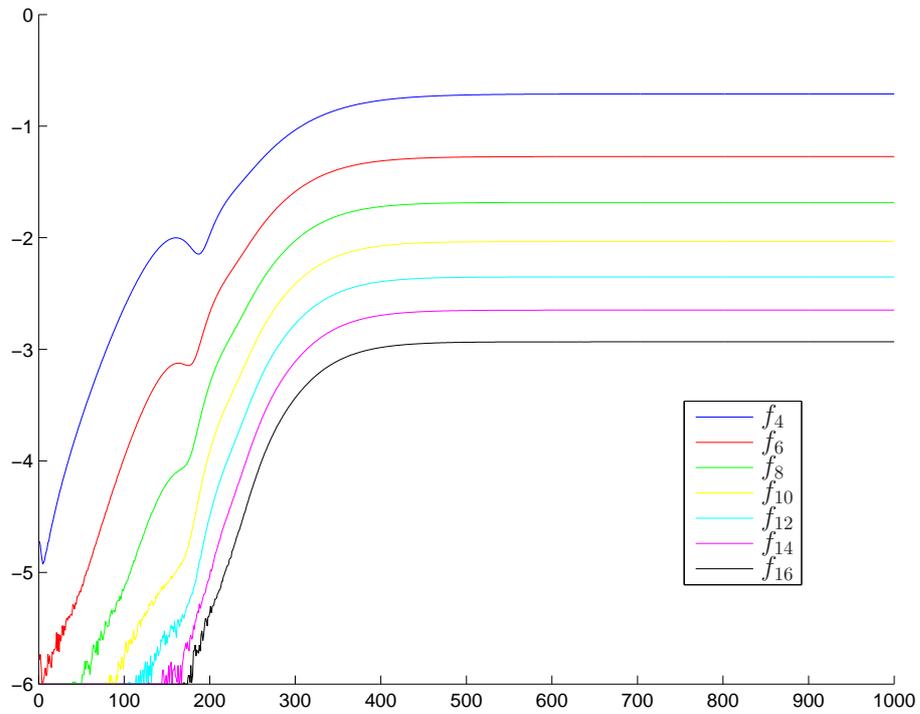}
\end{psfrags}
\caption{$\rm{L}^2$ errors between the polynomial solutions and the logarithmic solution.}\label{FIG:HH}
\end{figure}

On a mathematical point of view, it is interesting to study the error in $L^2$ norm:
$$
\left(\int_\Omega |u_{2n}-u|^2 dx\right)^{1/2}.
$$
We see on Figure \ref{FIG:HH} that for $n=2$, the error is important. It decreases with $n$ but is still 
significant for $n=3$. For $n\ge 6$, it is negligible.

\begin{figure}[htp]
\centering
\begin{psfrags}
\psfrag{Title}{{\Large \!\!\!\!\!\!\!\!\!\!\!\!\!\!\!\!\!\!\!\! Two dimensional case}}
\psfrag{Ylabel}{{\Large \!\!\!\!\!\!\!\!\!\!\!\!\!\!\!\!\!\!\!\!\!\!\!\!\!\!\!\!\!\! Log$_{10}($Energy$_{2n}$ - Energy$_{log})$}}
\psfrag{Xlabel}{{\Large \!\!\!\!\!\!\!\!\!\! Time iterations}}
\psfrag{F4}{{\Large$f_{4}$}}
\psfrag{F6}{{\Large$f_{6}$}}
\psfrag{F8}{{\Large$f_{8}$}}
\psfrag{F10}{{\Large$f_{10}$}}
\psfrag{F12}{{\Large$f_{12}$}}
\psfrag{F14}{{\Large$f_{14}$}}
\psfrag{F16}{{\Large$f_{16}$}}
\includegraphics[angle=0,width=12cm]{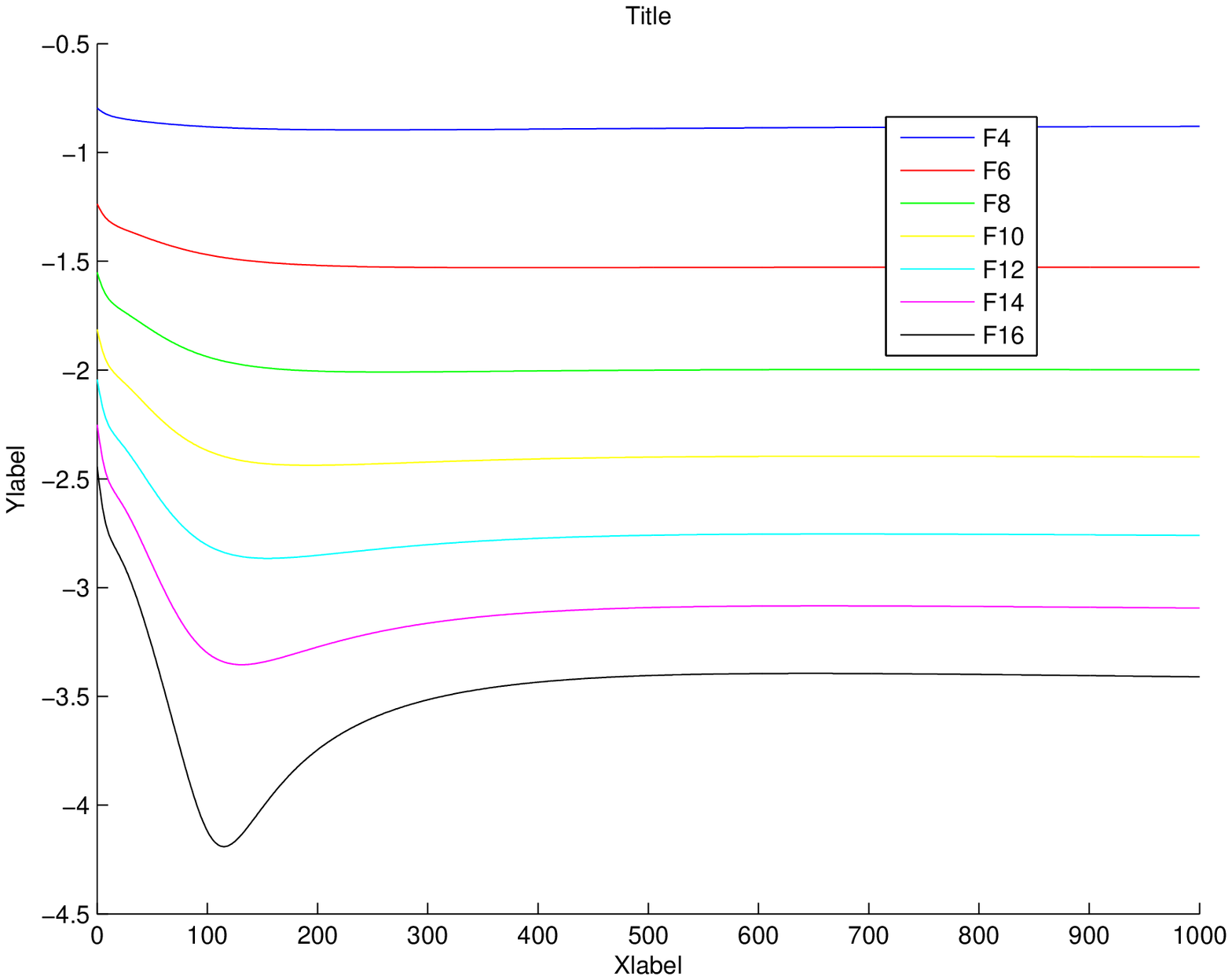}
\end{psfrags}
\caption{Polynomial energies versus logarithmic energy.}\label{FIG:H:2D}
\end{figure}

\begin{figure}[htp]
\centering
\begin{psfrags}
\psfrag{Title}{{\Large \!\!\!\!\!\!\!\!\!\!\!\!\!\!\!\!\!\!\!\! Two dimensional case}}
\psfrag{Ylabel}{}
\psfrag{Xlabel}{{\Large \!\!\!\!\!\!\!\!\!\! Time iterations}}
\psfrag{F4}{{\Large$f_{4}$}}
\psfrag{F6}{{\Large$f_{6}$}}
\psfrag{F8}{{\Large$f_{8}$}}
\psfrag{F10}{{\Large$f_{10}$}}
\psfrag{F12}{{\Large$f_{12}$}}
\psfrag{F14}{{\Large$f_{14}$}}
\psfrag{F16}{{\Large$f_{16}$}}
\includegraphics[angle=0,width=12cm]{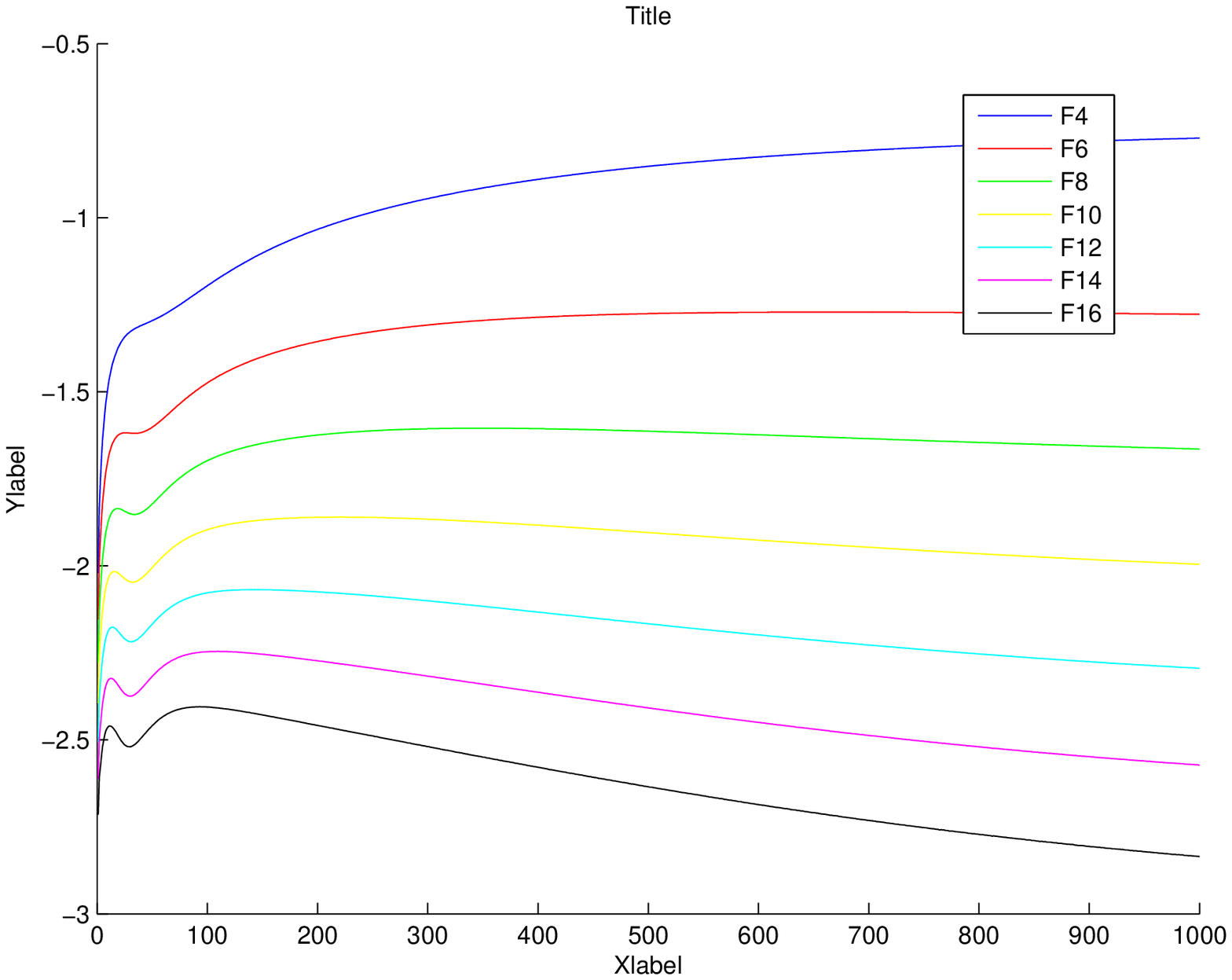}
\end{psfrags}
\caption{$\rm{L}^2$ errors between the polynomial solutions and the logarithmic solution.}\label{FIG:HH:2D}
\end{figure}

Figures \ref{FIG:H:2D} and \ref{FIG:HH:2D} present the same quantities for a two-dimensional 
spinodal decomposition. We observe the same quantitative difference. Note that we clearly see
 that the energy evolution slows down as the degree $n$ grows.


We conclude that the classical quartic approximation of the free energy may be considered as a good 
approximation for qualitative behaviour but it produces a significant error and accelerates the dynamics.
If precision is required, one should consider an approximation with a higher order polynomial.

\section{Validation of the numerical method. Choice of the degree of the elements}\label{S:3}
\begin{figure}[htp]
\centering
\subfigure[t=0]{
\label{FIG:12B.a}
\includegraphics[angle=0,width=4.5cm]{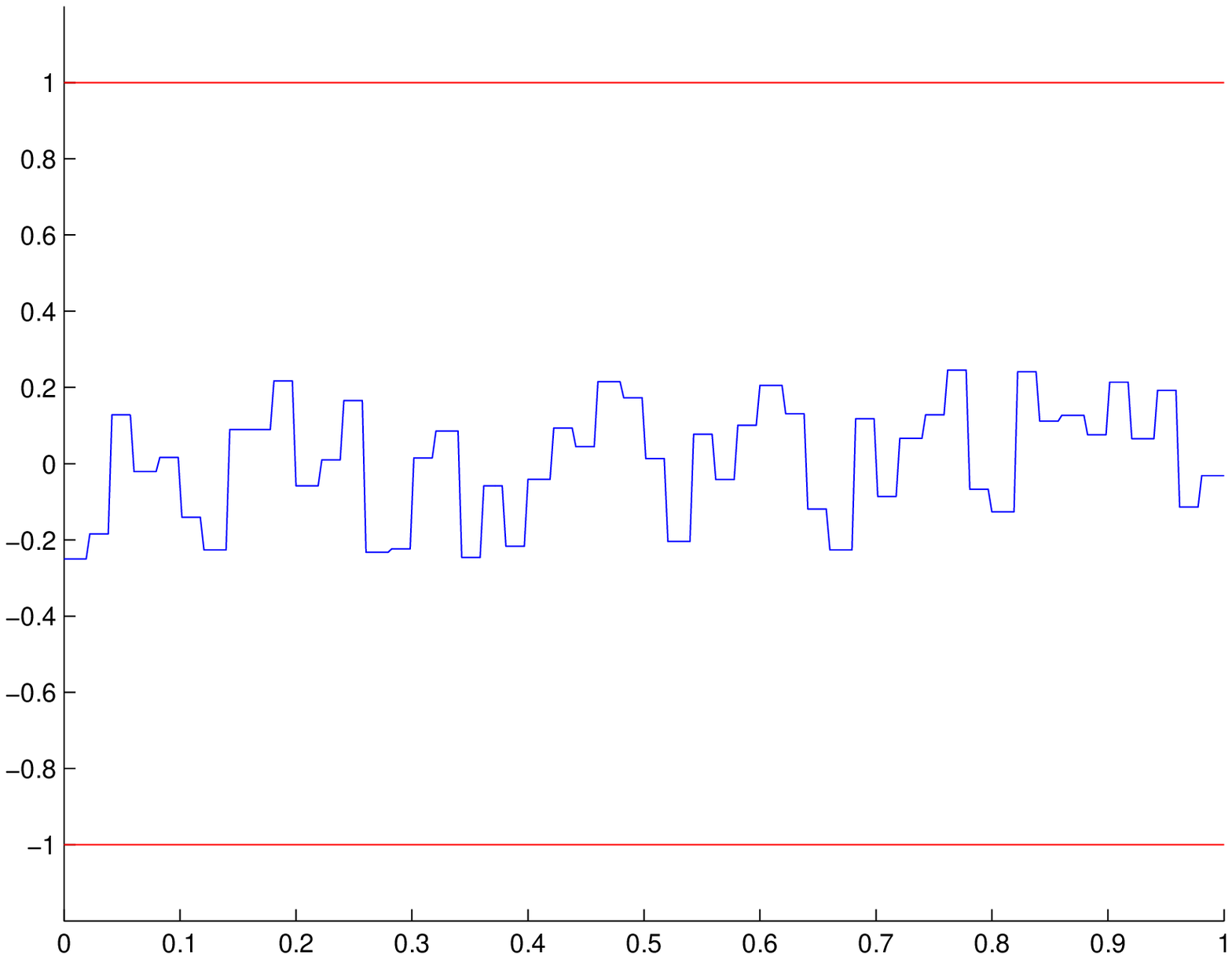}}
\hspace{0.3cm}
\subfigure[t=0.1]{
\label{FIG:12B.b}
\includegraphics[angle=0,width=4.5cm]{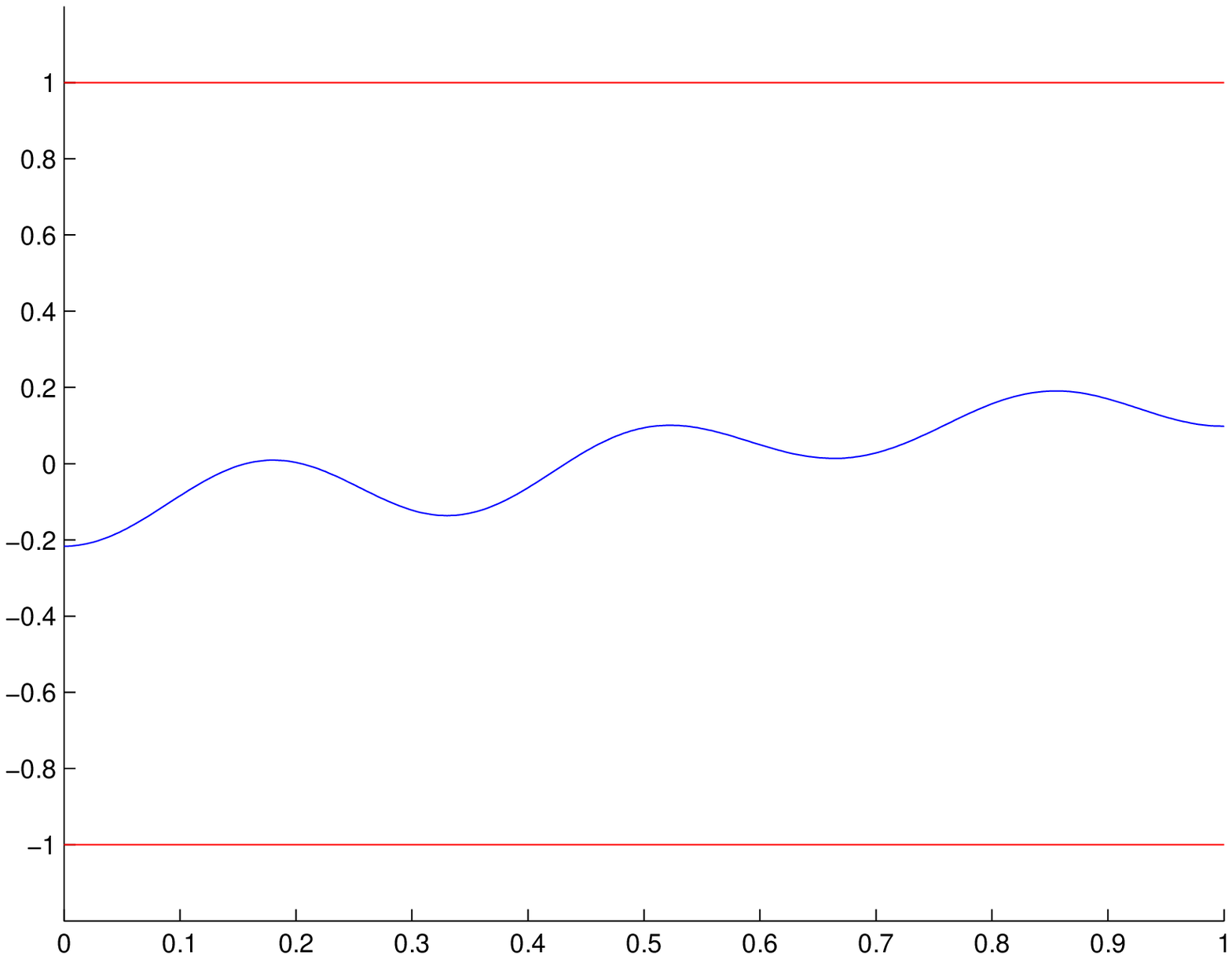}}
\hspace{0.3cm}
\subfigure[t=0.5]{
\label{FIG:12B.c}
\includegraphics[angle=0,width=4.5cm]{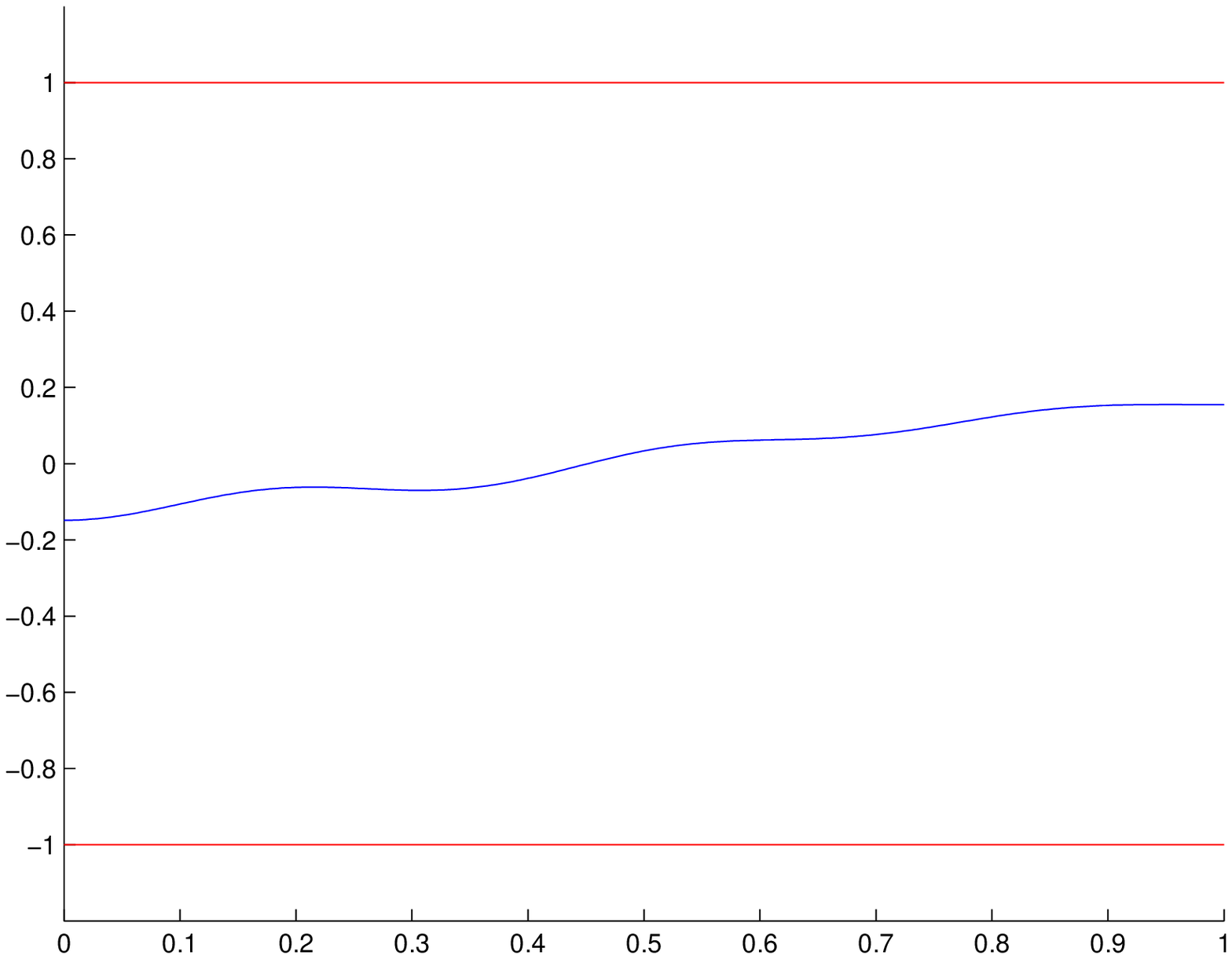}}
\\
\vspace{10pt}
\subfigure[t=10]{
\label{FIG:12B.d}
\includegraphics[angle=0,width=4.5cm]{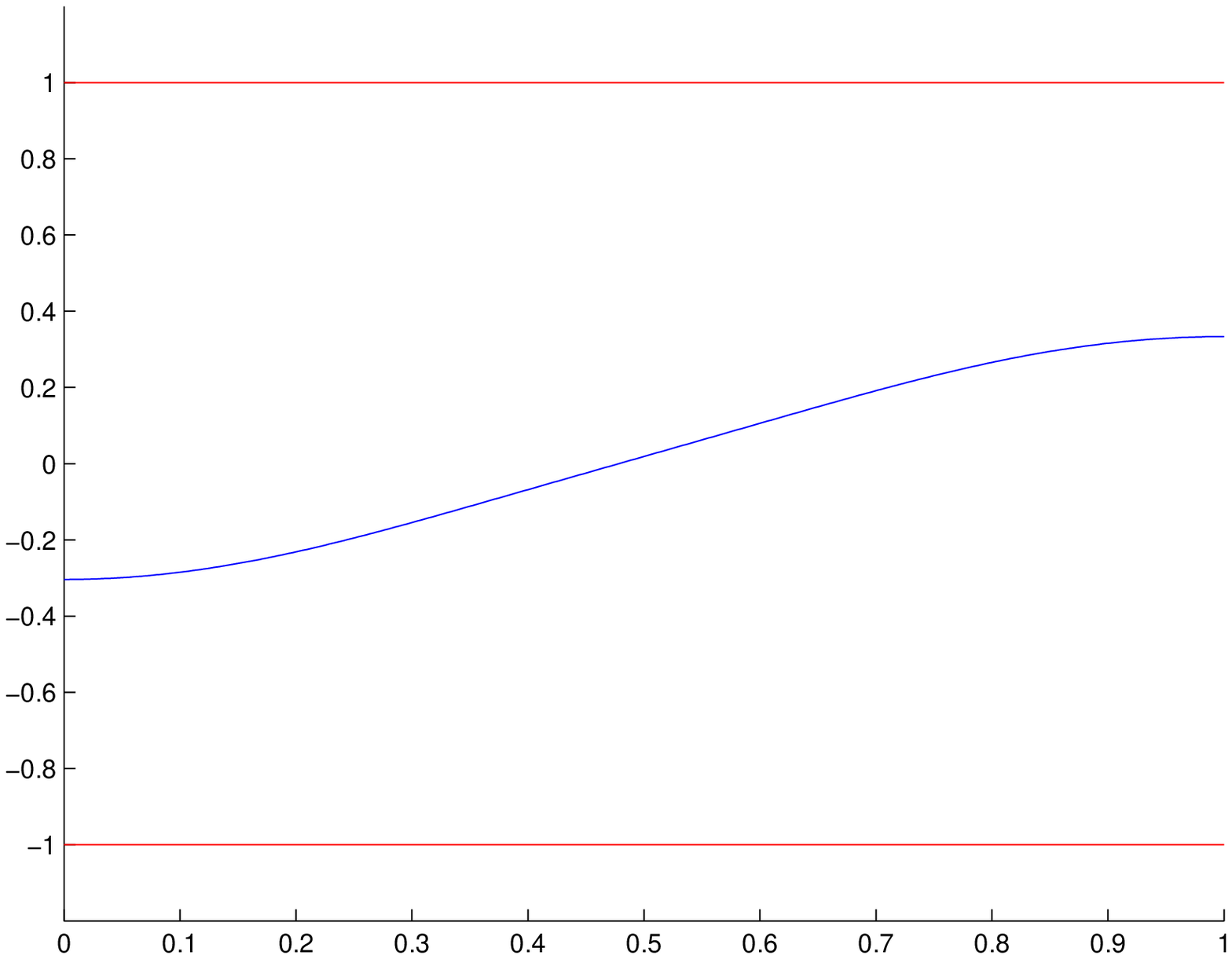}}
\hspace{0.3cm}
\subfigure[t=20]{
\label{FIG:12B.e}
\includegraphics[angle=0,width=4.5cm]{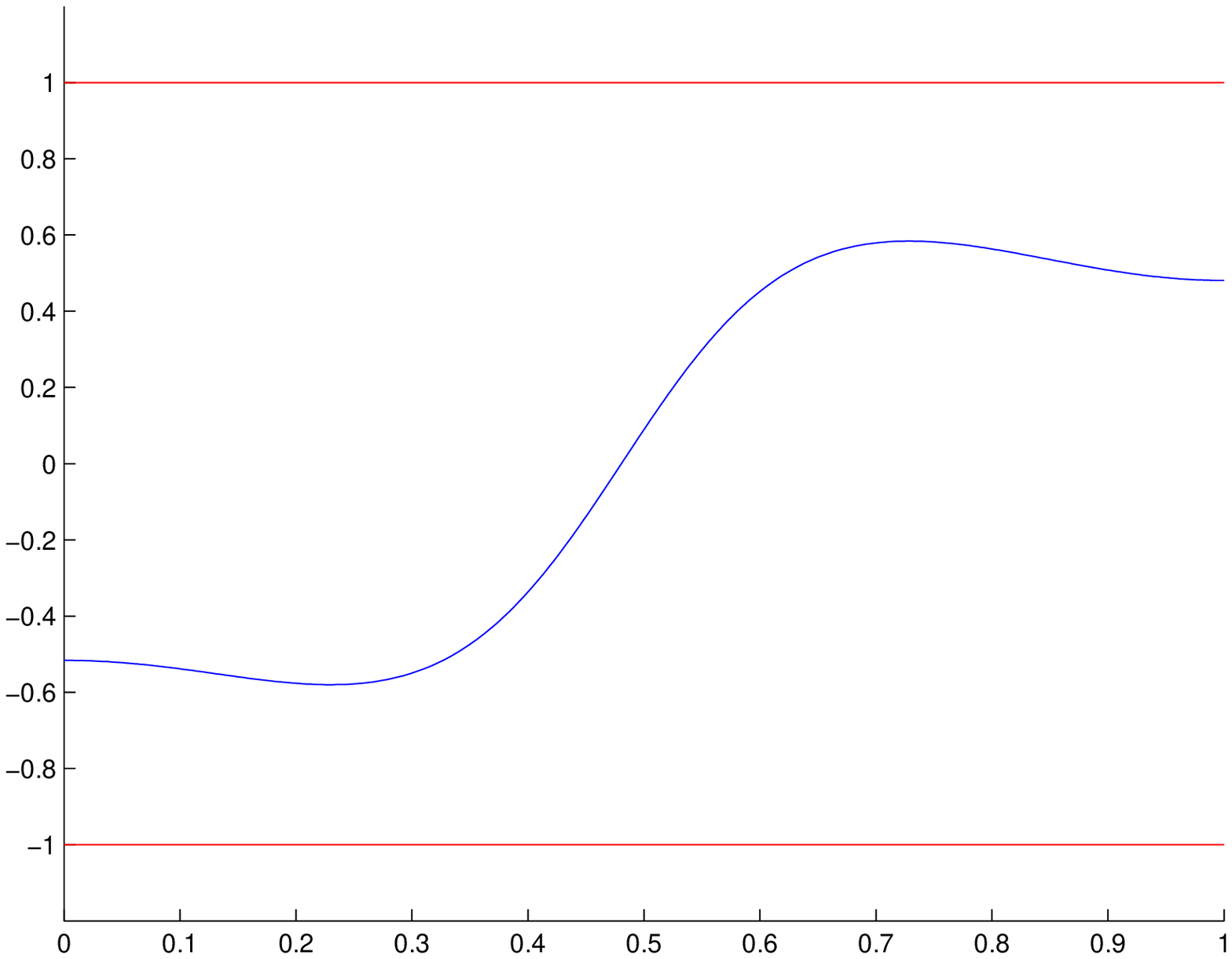}}
\hspace{0.3cm}
\subfigure[t=50]{
\label{FIG:12B.f}
\includegraphics[angle=0,width=4.5cm]{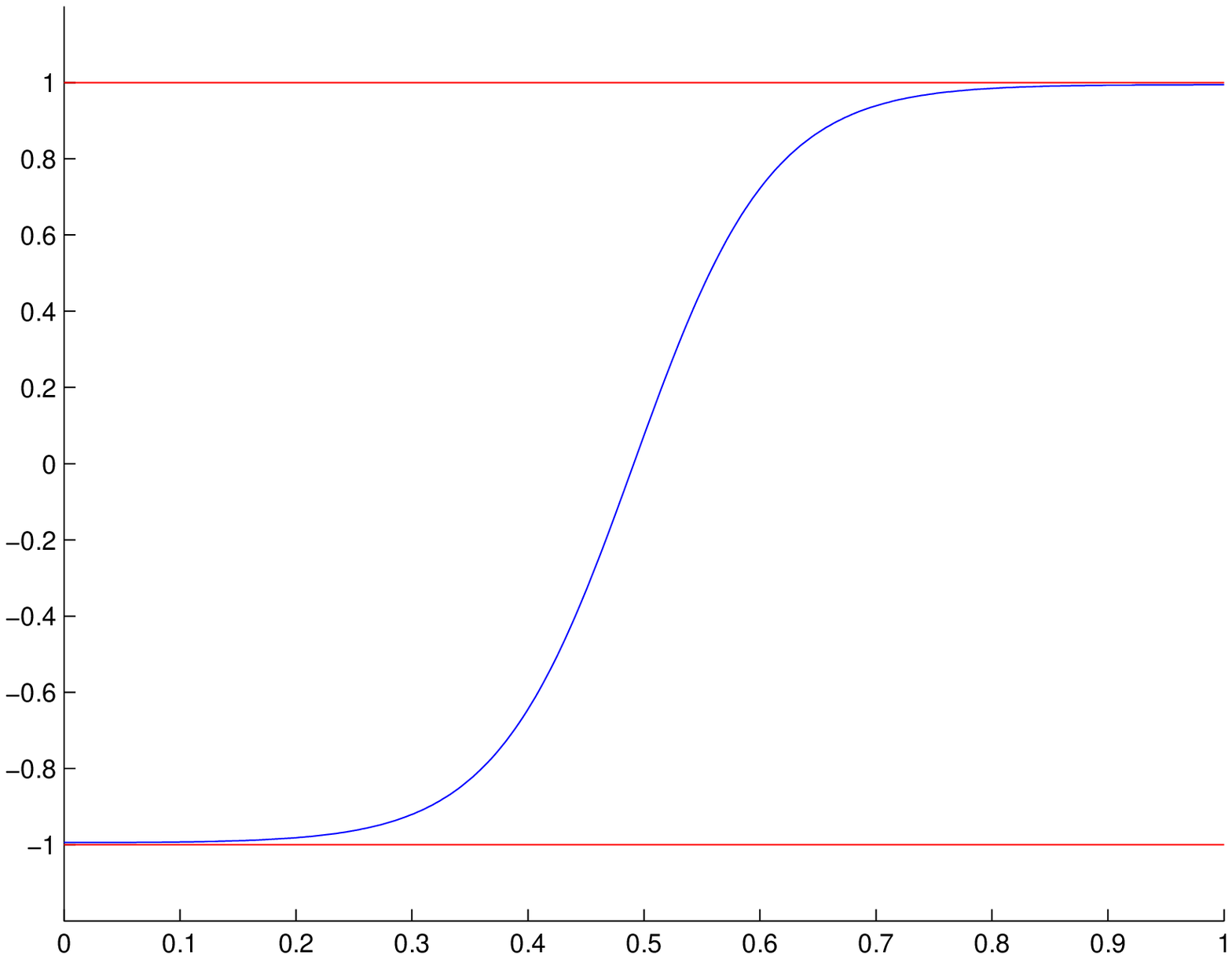}}
\caption{$Q_{1}$ solution on a mesh with 100 elements.}
\label{FIG:12B}
\end{figure}

On Figure \ref{FIG:12B}, we have drawn a numerical solution for different times. 
It is a $Q_{1}$ solution on a mesh with 100 elements under the quartic double-well potential. 
On figure \ref{FIG:12B.f}, the solution has reached its stable state and has binodal values $\pm 1$ on the boundary.

In our first set of tests, we start with the same initial state near the ``tanh profile'' solution. The evolutions are driven by the quartic potential function. We wait for the stabilization of all the solutions and study the error on the energies and on the slopes of the interface.

\begin{figure}[htp]
\centering
\includegraphics[angle=0,width=12cm]{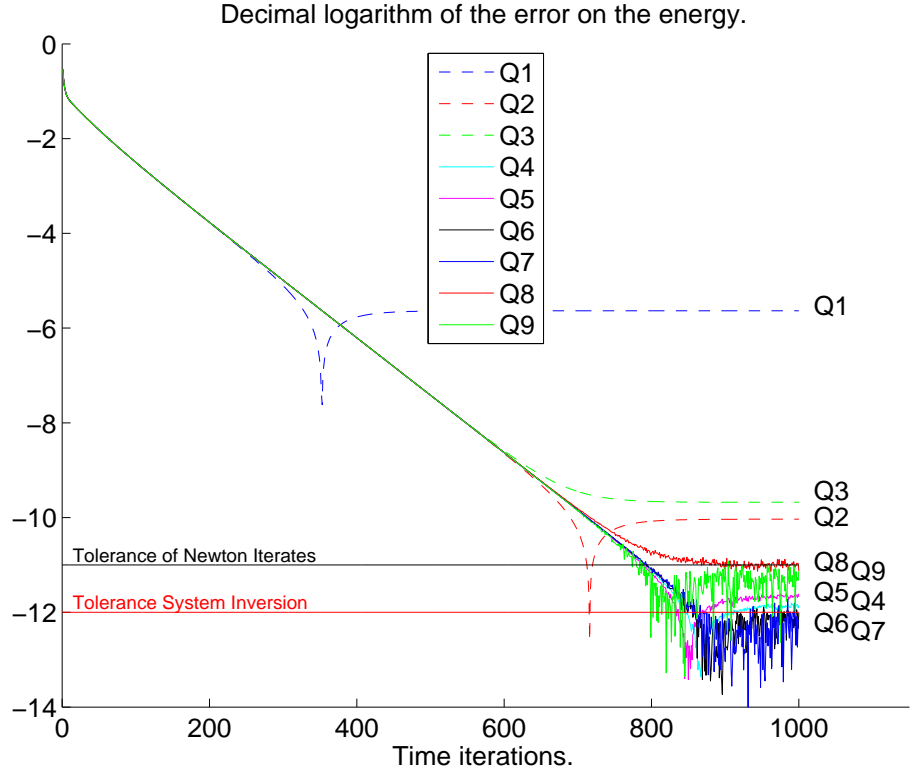}
\caption{Energies during an evolution.}\label{FIG:12}
\end{figure}

Remark that we can explicitly compute the energy of the explicit solution. And since the energy of a numerical simulation is decreasing in time, this energy should converge to the energy of the explicit solution. 
Figure \ref{FIG:12} shows the evolutions of the errors between the numerical energies and the explicit energy according to the degree of the polynomial space $\PP$. Before the 800th iteration in time, the solutions are not stable. They try to minimize their energies. After the 800th iteration, all the solutions are in a stable state. We can see that the evolutions are qualitatively similar at the beginning, but the elements $Q_{1}$, $Q_{2}$ and $Q_{3}$ don't achieve the tolerance zone, whereas the other elements do. However, $Q_{2}$ and $Q_{3}$
give a very good result.

The slope of the interface is an essential physical quantity. So we have compared the errors on the slopes between the numerical solutions and the theoretical solution. Note that these slopes correspond to the values of the derivatives of the numerical solutions at the 
interface and our finite elements have not a $\mathcal{C}^1$ regularity. 

\begin{figure}[htp]
\centering
\subfigure[Mesh 18 - Q1]{
\label{FIG:2.a}
\includegraphics[angle=0,width=4.5cm]{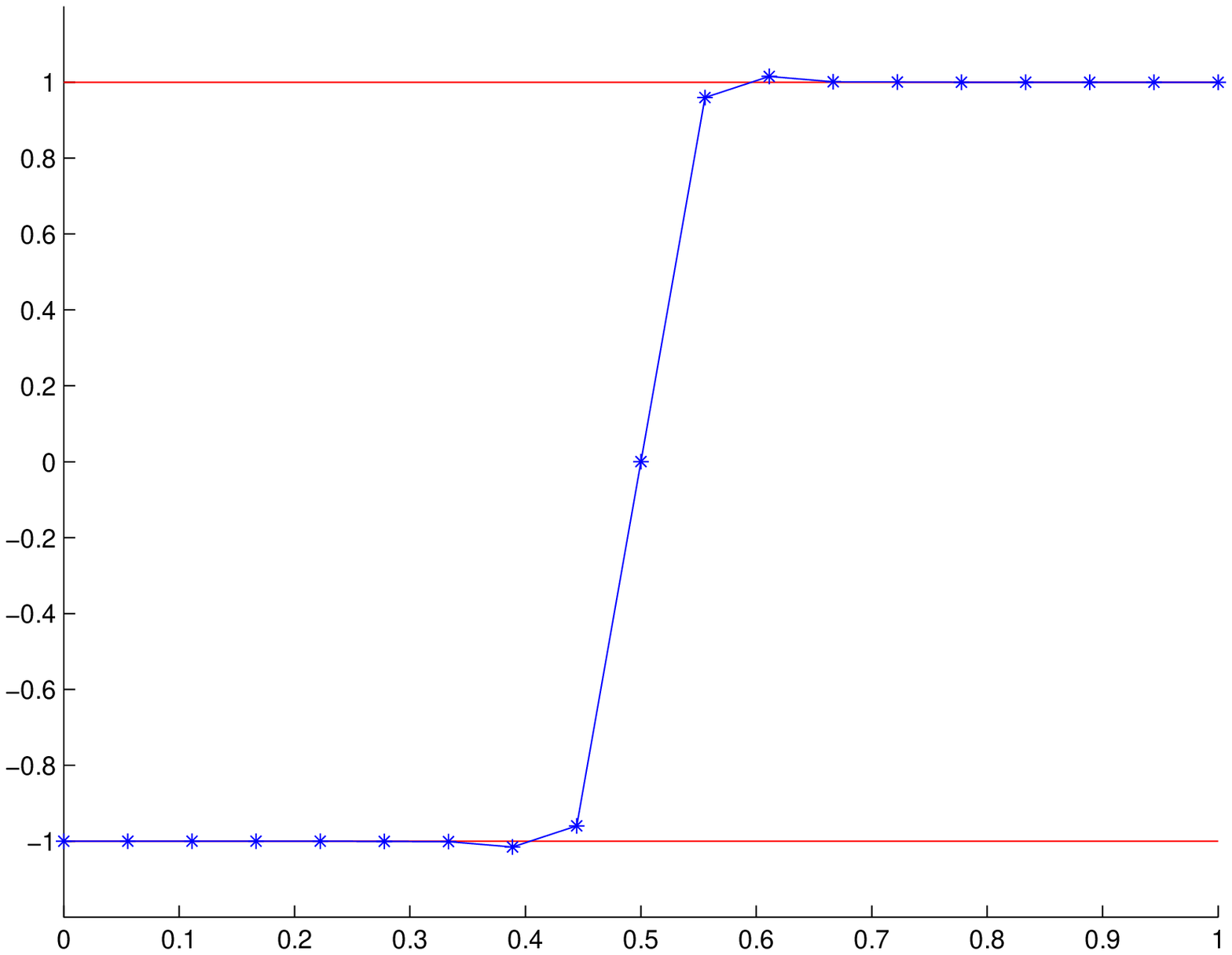}}
\hspace{0.3cm}
\subfigure[Mesh 36 - Q1]{
\label{FIG:2.b}
\includegraphics[angle=0,width=4.5cm]{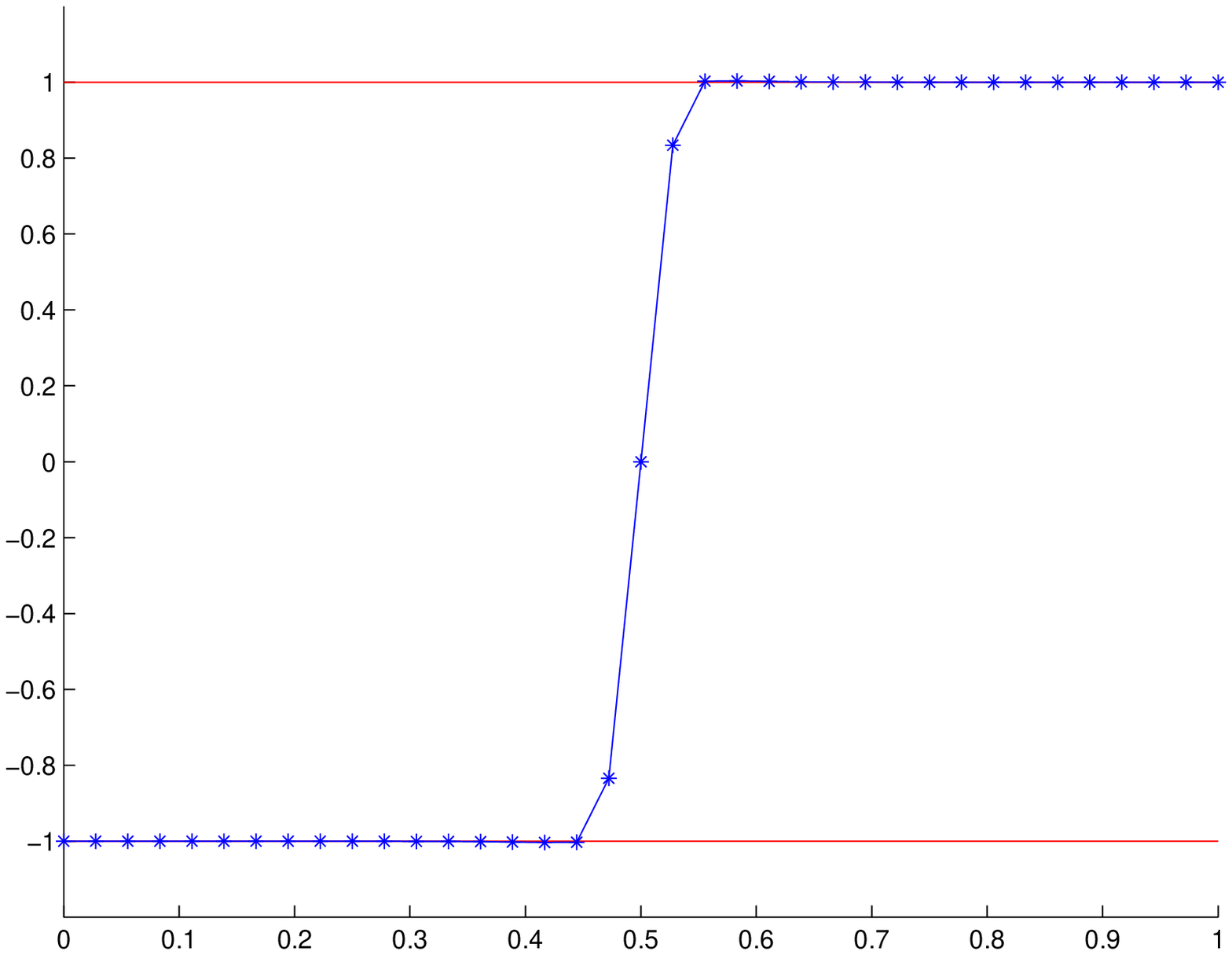}}
\hspace{0.3cm}
\subfigure[Mesh 72 - Q1]{
\label{FIG:2.c}
\includegraphics[angle=0,width=4.5cm]{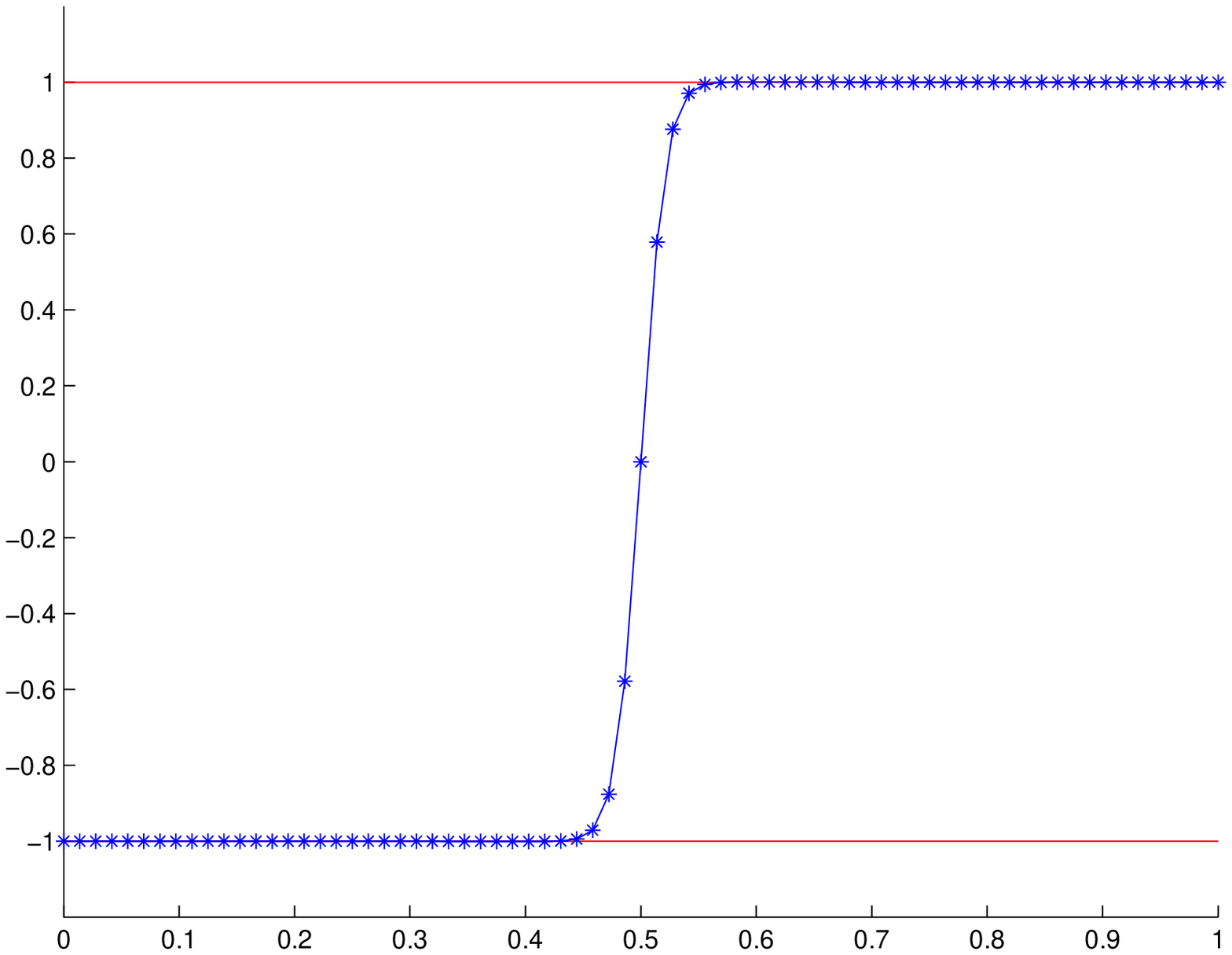}}
\\
\vspace{10pt}
\subfigure[Mesh 9 - Q2]{
\label{FIG:2.d}
\includegraphics[angle=0,width=4.5cm]{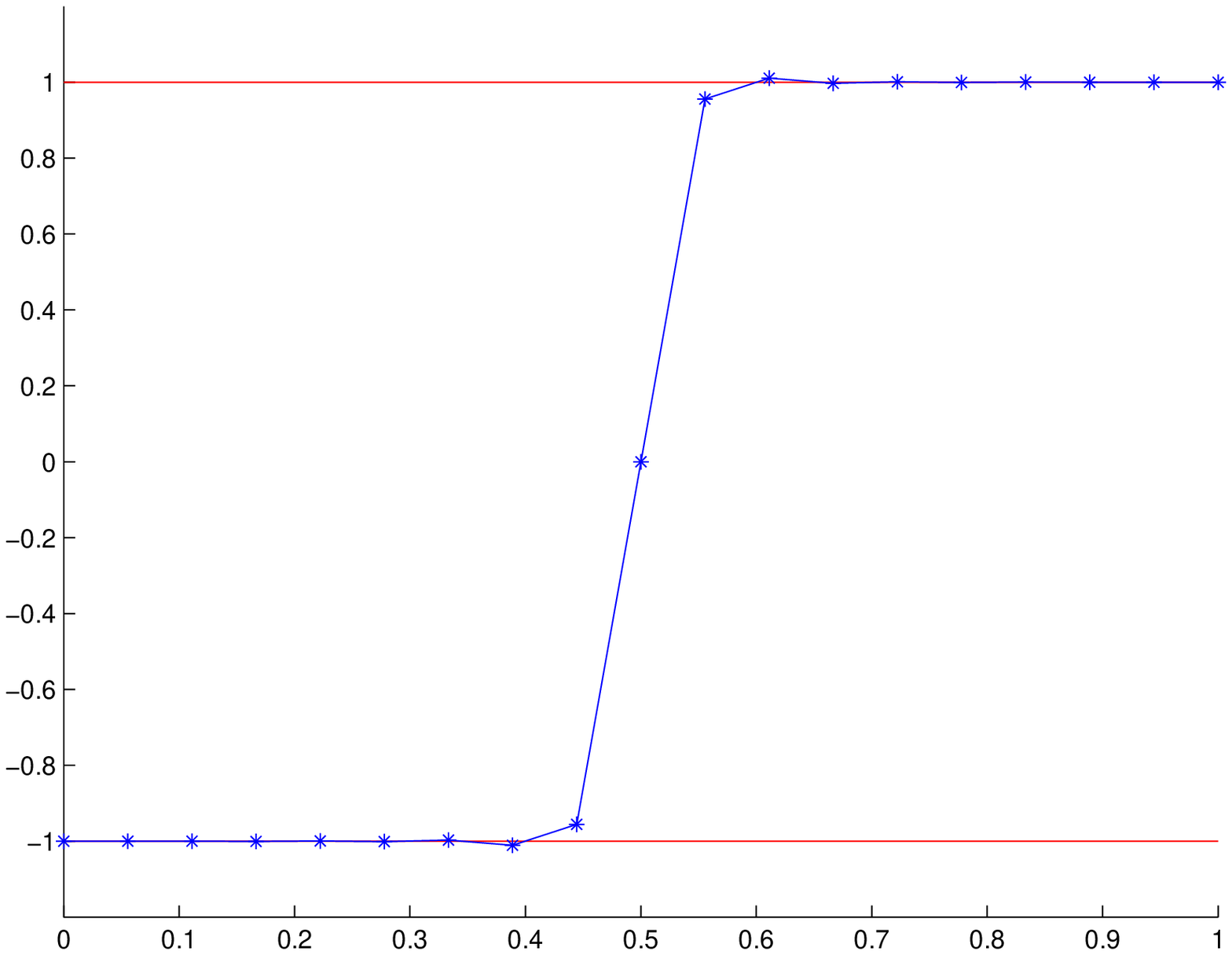}}
\hspace{0.3cm}
\subfigure[Mesh 18 - Q2]{
\label{FIG:2.e}
\includegraphics[angle=0,width=4.5cm]{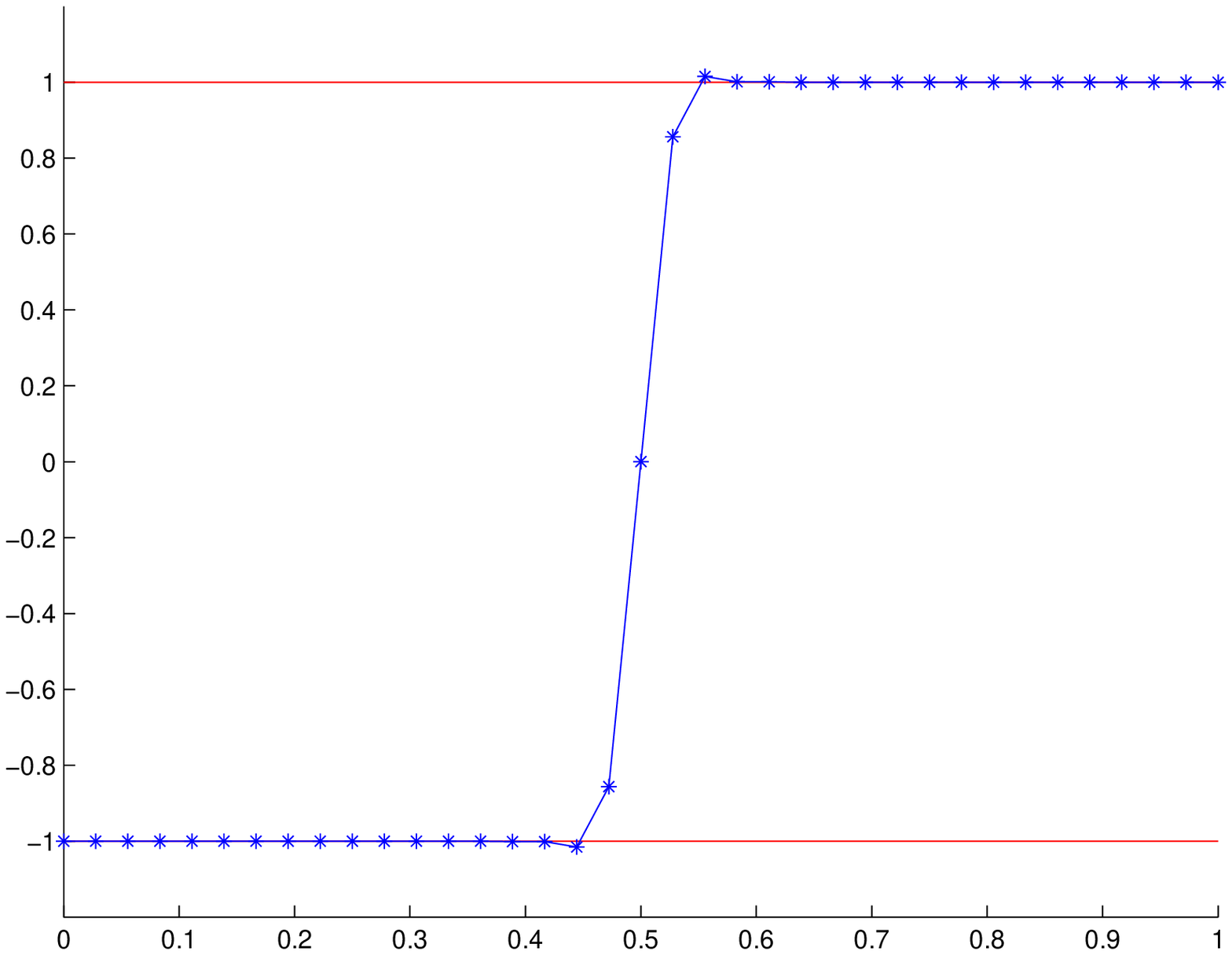}}
\hspace{0.3cm}
\subfigure[Mesh 36 - Q2]{
\label{FIG:2.f}
\includegraphics[angle=0,width=4.5cm]{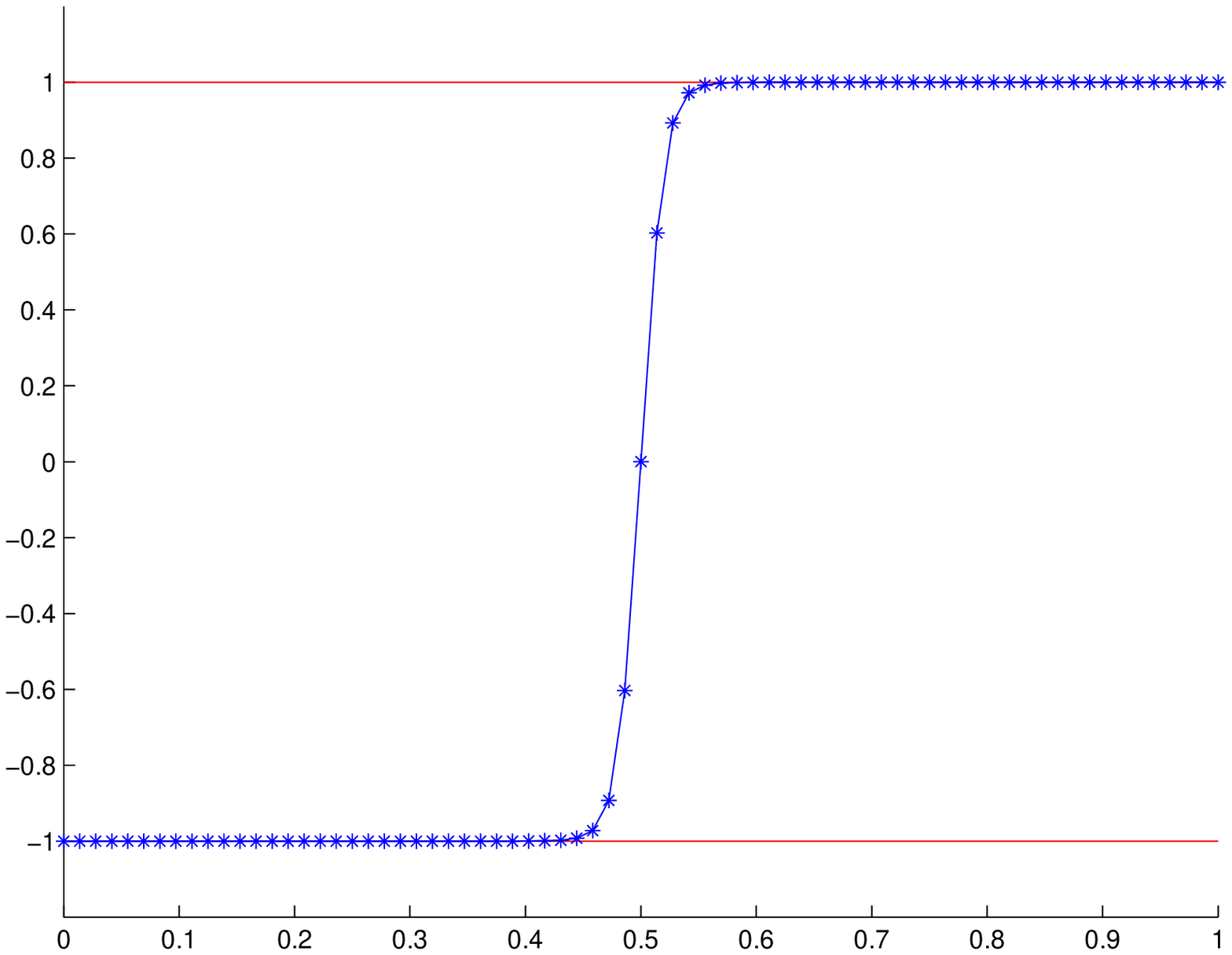}}
\\
\vspace{10pt}
\subfigure[Mesh 6 - Q3]{
\label{FIG:2.g}
\includegraphics[angle=0,width=4.5cm]{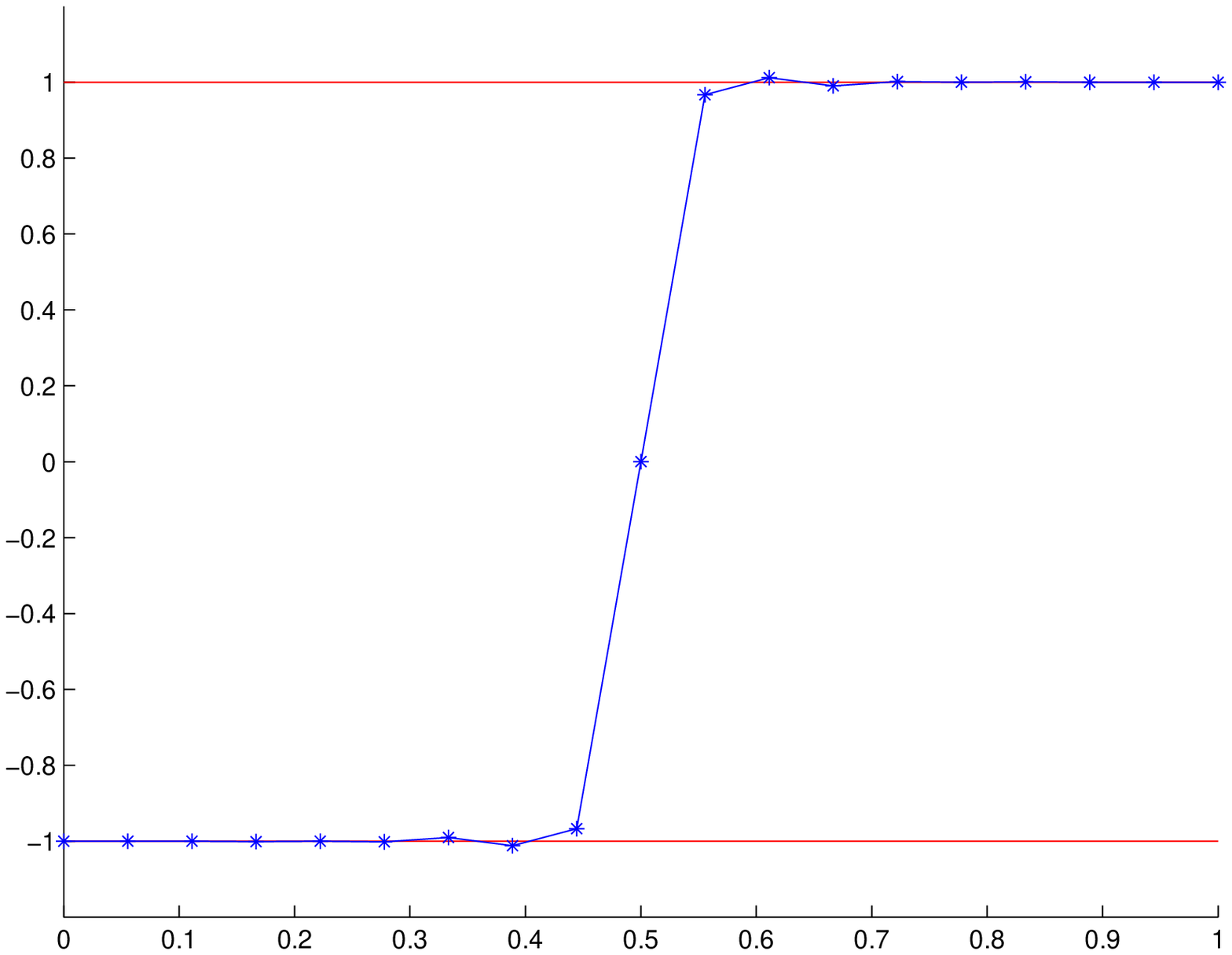}}
\hspace{0.3cm}
\subfigure[Mesh 12 - Q3]{
\label{FIG:2.h}
\includegraphics[angle=0,width=4.5cm]{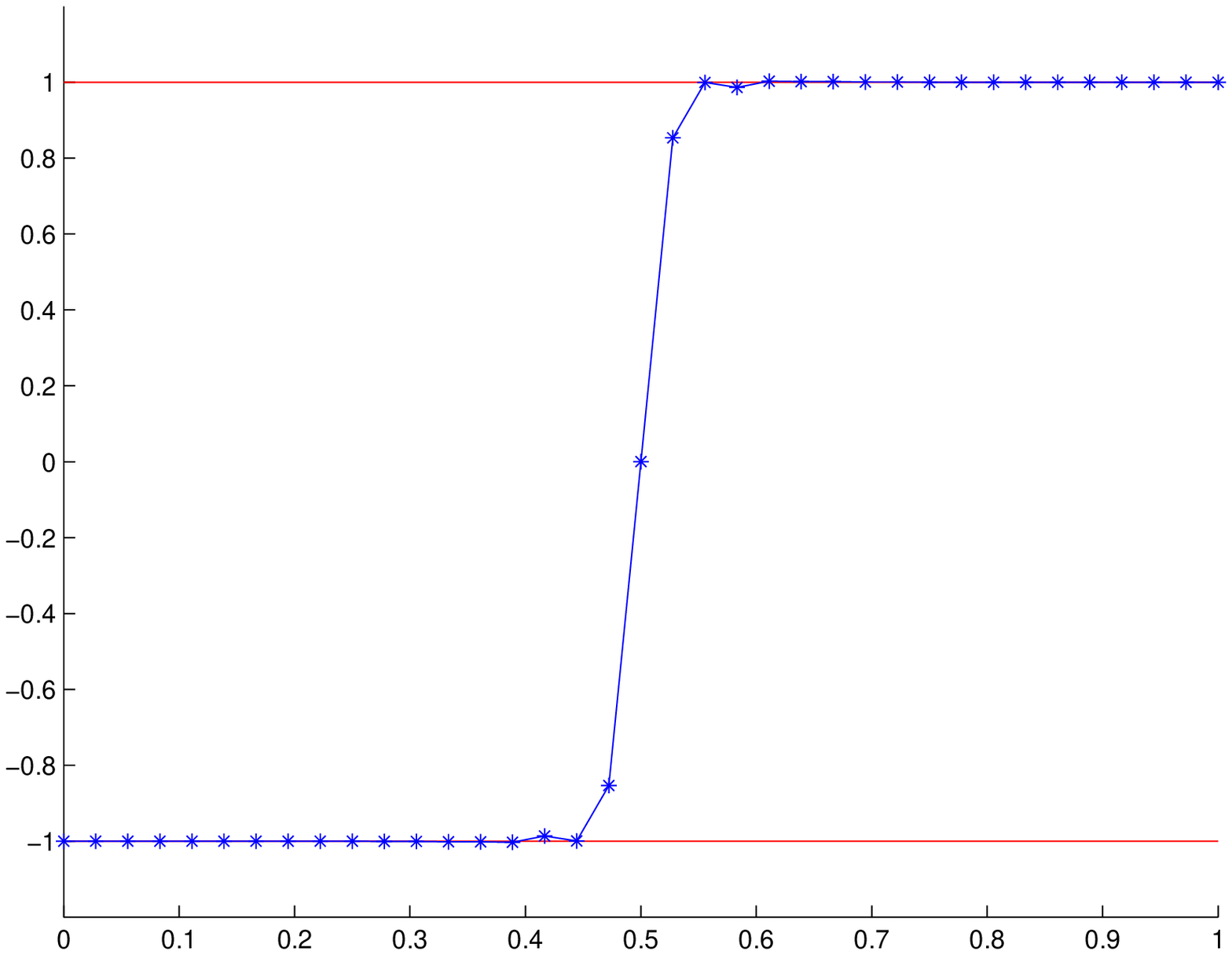}}
\hspace{0.3cm}
\subfigure[Mesh 24 - Q3]{
\label{FIG:2.i}
\includegraphics[angle=0,width=4.5cm]{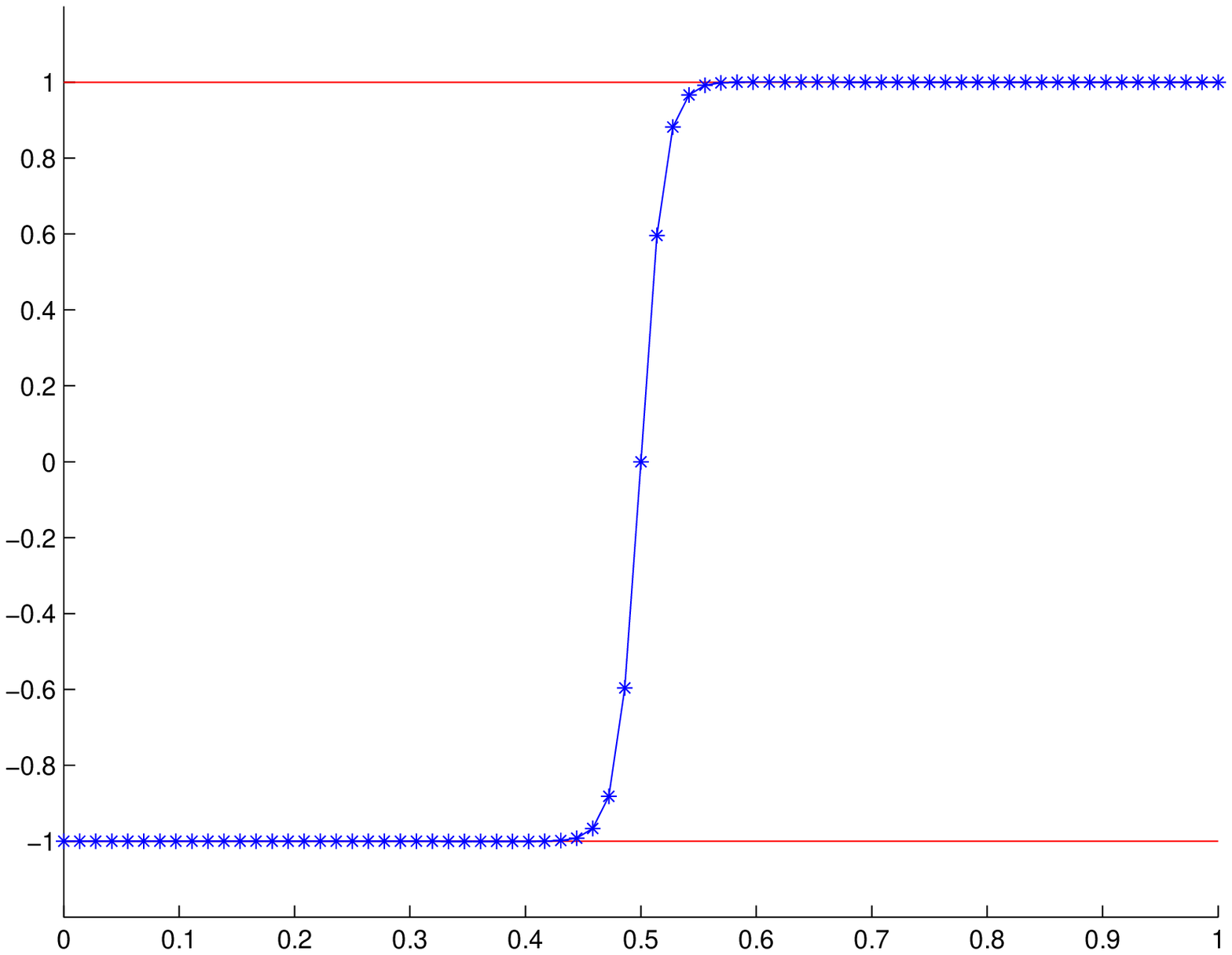}}
\caption{Fitted curves on the ``tanh-profile''.}
\label{FIG:2}
\end{figure}

Under the quartic double-well potential \eqref{Eq:2.2}, we have an explicit slope $\mu$ for the stationary solution. On Figure~\ref{FIG:2}, we present the numerical solution for $\Omega=[0,1]$ (blue stars), and the ``tanh-profile'' whose coefficients $u_+$ and $\mu$ have been fitted to the data. The fitting on $u_{+}$ corresponds to the value of the solution on the boundaries of the domain $\Omega$. And the fitting on $\mu$ corresponds to a least square method between the numerical solution and a ``tanh-profile'' solution interpolated on the same meshes. The ``tanh-profile'' (defined over $\mathbb R$) may be considered as a good approximation of the solution on $\Omega =[0,1]$ since the interface is very thin. 

If we want to compare the solutions between a $Q_{1}$ simulation and a $Q_{10}$ simulation, we need to compare the two simulations under a same complexity which, up to the inversions of the 
linear systems, corresponds to a similar computational cost. In the one dimensional case, the complexity corresponds to the value $Degree \times Number\ of\ elements$. For a $Q_{10}$ simulation, we only need a mesh with 10 times less elements than for a $Q_{1}$ simulation. 

\begin{figure}[htp]
\centering
\subfigure[Versus the number of elements of the mesh]{
\begin{psfrags}
\psfrag{Title}{}
\psfrag{Xlabel}{}
\psfrag{Ylabel}{}
\label{FIG:3ter.a}
\includegraphics[angle=0,width=6cm]{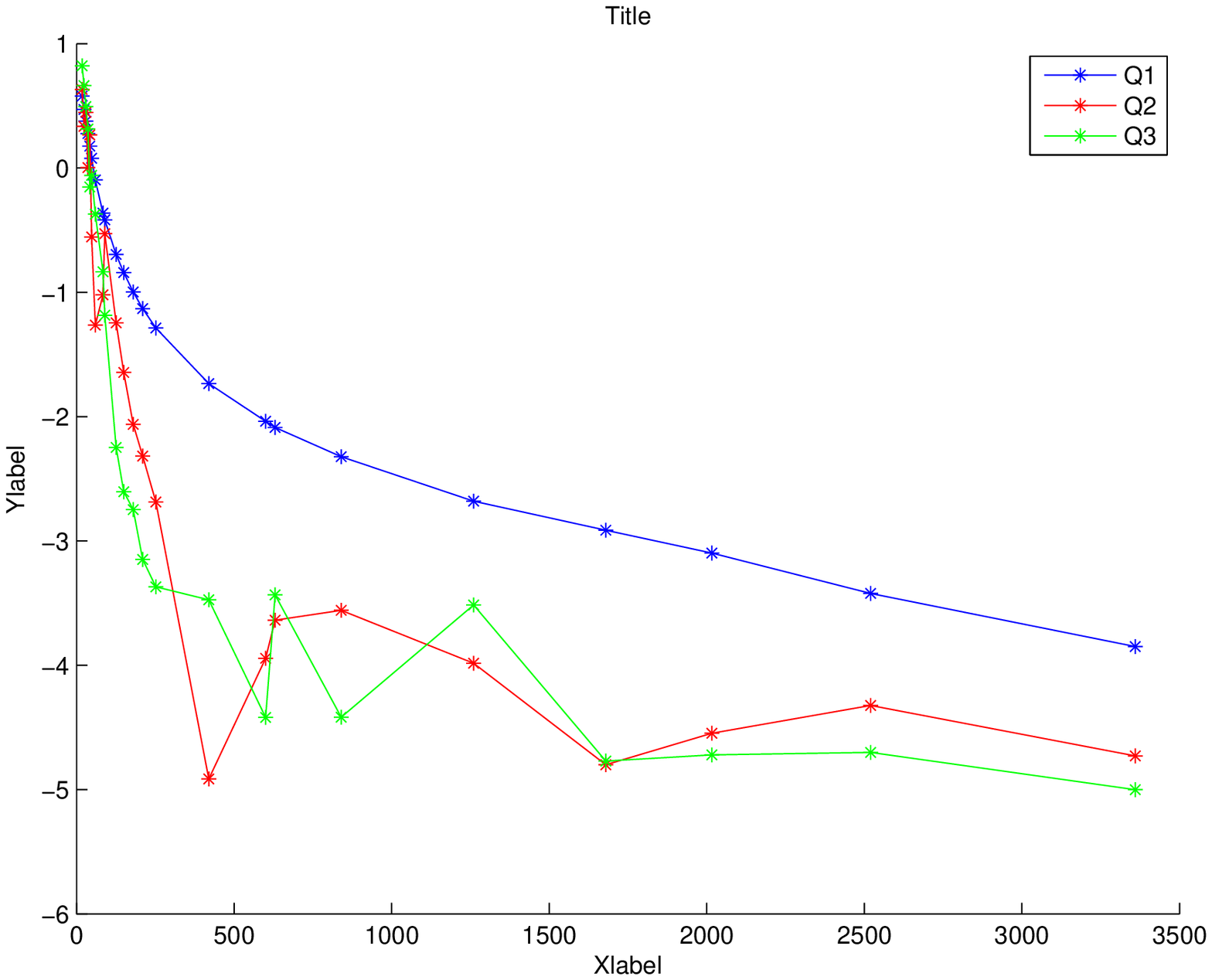}
\end{psfrags}}
\hspace{1cm}
\subfigure[Versus the logarithm of the number of elements of the mesh]{
\begin{psfrags}
\psfrag{Title}{}
\psfrag{Xlabel}{}
\psfrag{Ylabel}{}
\label{FIG:3ter.b}
\includegraphics[angle=0,width=6cm]{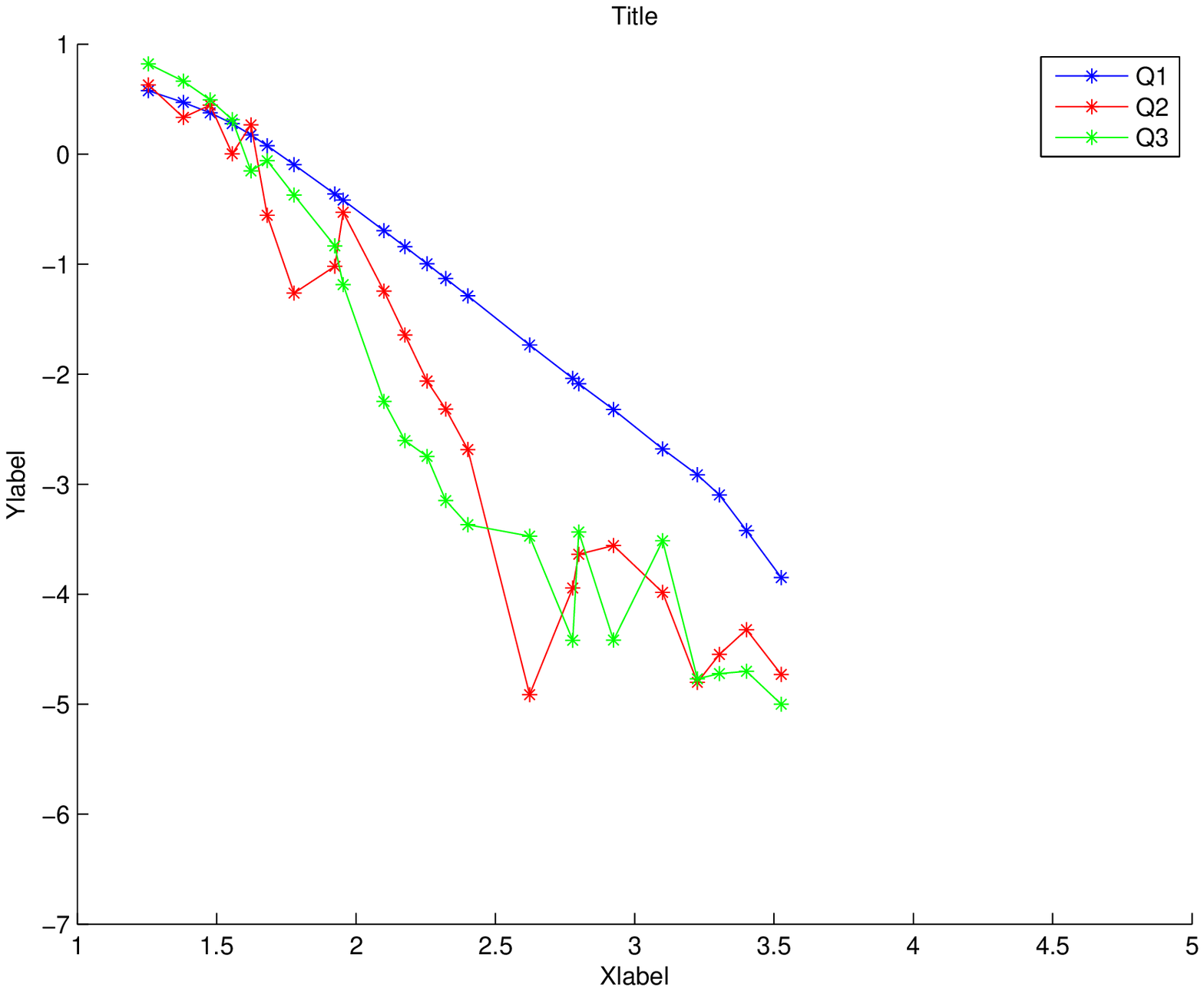}
\end{psfrags}}
\caption{Comparison between the errors on the slope under the same complexity.}
\label{FIG:3ter}
\end{figure}

Figure \ref{FIG:2} represents the numerical solution over mesh grids with three different complexities 18, 36 and 72, and under polynomial functions of degree 1, 2 and 3. For instance, for the elements $Q_{2}$, it corresponds to the mesh grids with 9, 18 and 36 elements. If we increase the number of elements or the degree of the polynomial space $\PP$, then we obtain a better approximation of the slope of the ``tanh-profil'' solution. But for the same complexity, the curves are qualitatively similar. Figures \ref{FIG:3ter.a} and \ref{FIG:3ter.b} show the evolution of this approximation error according to the complexity for $Q_{1}$, $Q_{2}$ and $Q_{3}$ simulations. On Figure \ref{FIG:3ter.b}, we have used a logarithmic scale in order to compare the rate of the convergence.
\clearpage
We obviously conclude that, for elements $Q_{1}$, $Q_{2}$ or $Q_{3}$, a fine mesh allows a better approximation. But the $Q_{2}$ and $Q_{3}$ elements seem to reach faster a saturation. They only need 500 elements in order to reach a $10^{-5}$ precision, whereas the $Q_{1}$ elements need 5000 elements ! Figure \ref{FIG:3ter.b} highlights this better speed on the approximation error of the slope. But $Q_{2}$ and $Q_{3}$ elements seems to have a similar speed before reaching the saturation zone. 

\begin{figure}[htp]
\centering
\subfigure[DoF 90]{
\begin{psfrags}
\psfrag{Title}{}
\psfrag{Xlabel}{}
\psfrag{Ylabel}{}
\psfrag{Q1}{\!\!\huge{\color{blue}$Q_{1}$}}
\psfrag{Q2}{\!\!\!\!\huge{\color{blue}$Q_{2}$}}
\psfrag{Q3}{\!\!\!\!\huge{\color{blue}$Q_{3}$}}
\psfrag{Q4}{\!\!\!\!\huge{\color{blue}$Q_{4}$}}
\psfrag{Q5}{\!\!\!\!\huge{\color{blue}$Q_{5}$}}
\psfrag{Q6}{\!\!\!\!\huge{\color{blue}$Q_{6}$}}
\psfrag{Q7}{\!\!\!\!\huge{\color{blue}$Q_{7}$}}
\psfrag{Q8}{\!\!\!\!\huge{\color{blue}$Q_{8}$}}
\psfrag{Q9}{\!\!\!\!\huge{\color{blue}$Q_{9}$}}
\psfrag{Q10}{\!\!\!\!\huge{\color{blue}$Q_{10}$}}
\label{FIG:4bis.a}
\includegraphics[angle=0,width=4.5cm]{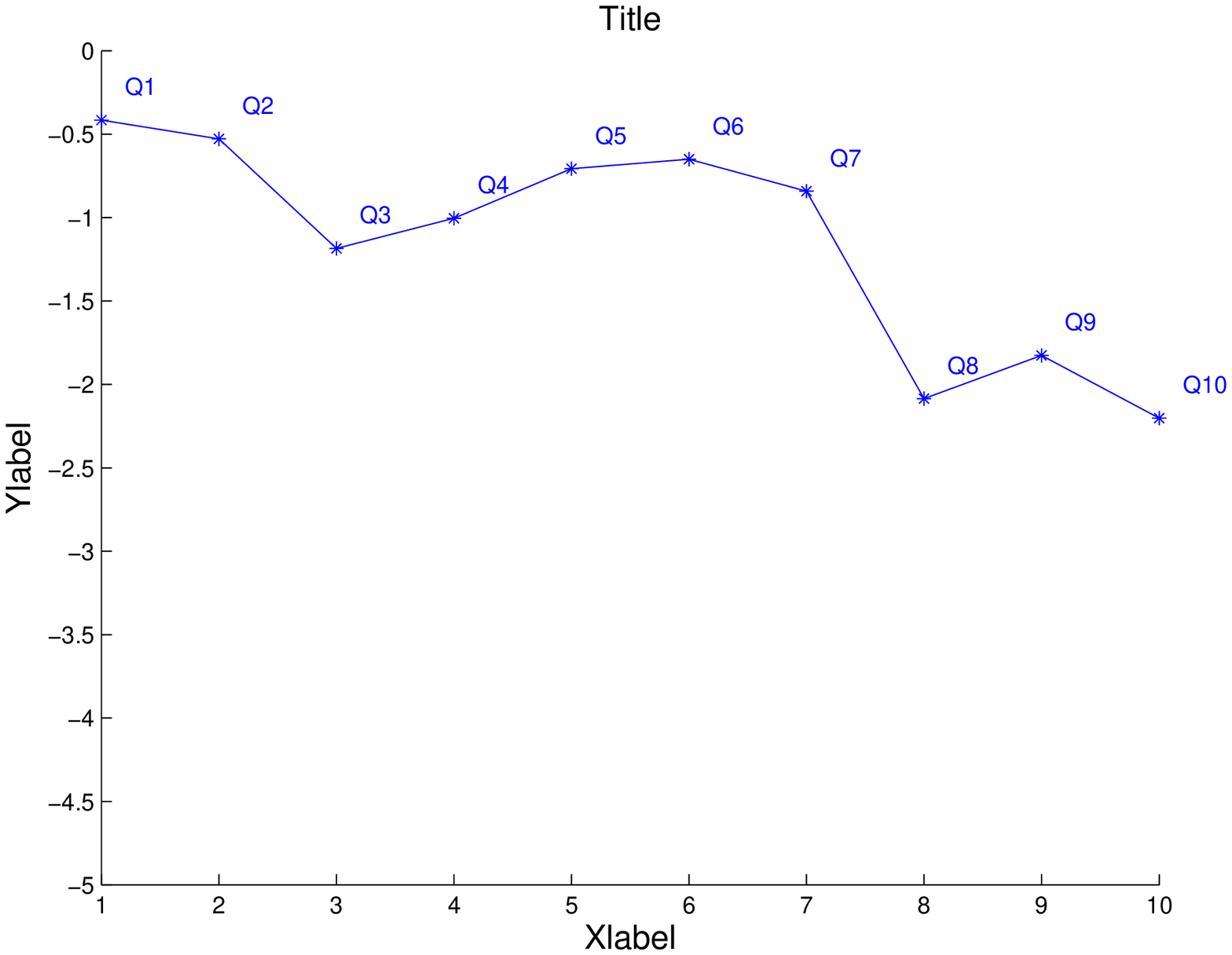}
\end{psfrags}}
\hspace{0.3cm}
\subfigure[DoF 300]{
\begin{psfrags}
\psfrag{Title}{}
\psfrag{Xlabel}{}
\psfrag{Ylabel}{}
\psfrag{Q1}{\!\!\huge{\color{blue}$Q_{1}$}}
\psfrag{Q2}{\!\!\!\!\huge{\color{blue}$Q_{2}$}}
\psfrag{Q3}{\!\!\!\!\huge{\color{blue}$Q_{3}$}}
\psfrag{Q4}{\!\!\!\!\huge{\color{blue}$Q_{4}$}}
\psfrag{Q5}{\!\!\!\!\huge{\color{blue}$Q_{5}$}}
\psfrag{Q6}{\!\!\!\!\huge{\color{blue}$Q_{6}$}}
\psfrag{Q7}{\!\!\!\!\huge{\color{blue}$Q_{7}$}}
\psfrag{Q8}{\!\!\!\!\huge{\color{blue}$Q_{8}$}}
\psfrag{Q9}{\!\!\!\!\huge{\color{blue}$Q_{9}$}}
\psfrag{Q10}{\!\!\!\!\huge{\color{blue}$Q_{10}$}}
\label{FIG:4bis.b}
\includegraphics[angle=0,width=4.5cm]{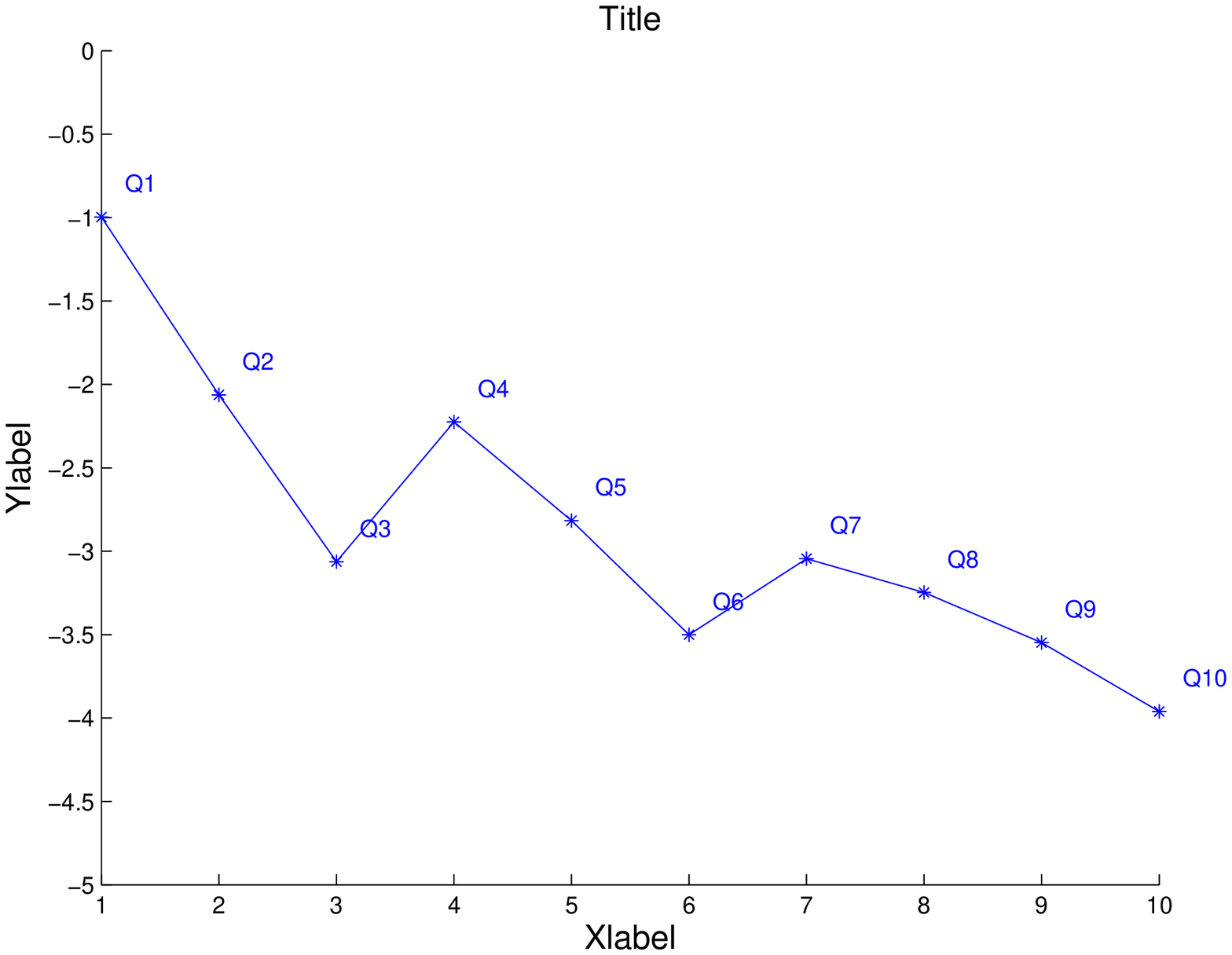}
\end{psfrags}}
\hspace{0.3cm}
\subfigure[DoF 630]{
\begin{psfrags}
\psfrag{Title}{}
\psfrag{Xlabel}{}
\psfrag{Ylabel}{}
\psfrag{Q1}{\!\!\huge{\color{blue}$Q_{1}$}}
\psfrag{Q2}{\!\!\!\!\huge{\color{blue}$Q_{2}$}}
\psfrag{Q3}{\!\!\!\!\huge{\color{blue}$Q_{3}$}}
\psfrag{Q4}{\!\!\!\!\huge{\color{blue}$Q_{4}$}}
\psfrag{Q5}{\!\!\!\!\huge{\color{blue}$Q_{5}$}}
\psfrag{Q6}{\!\!\!\!\huge{\color{blue}$Q_{6}$}}
\psfrag{Q7}{\!\!\!\!\huge{\color{blue}$Q_{7}$}}
\psfrag{Q8}{\!\!\!\!\huge{\color{blue}$Q_{8}$}}
\psfrag{Q9}{\!\!\!\!\huge{\color{blue}$Q_{9}$}}
\psfrag{Q10}{\!\!\!\!\huge{\color{blue}$Q_{10}$}}
\label{FIG:4bis.c}
\includegraphics[angle=0,width=4.5cm]{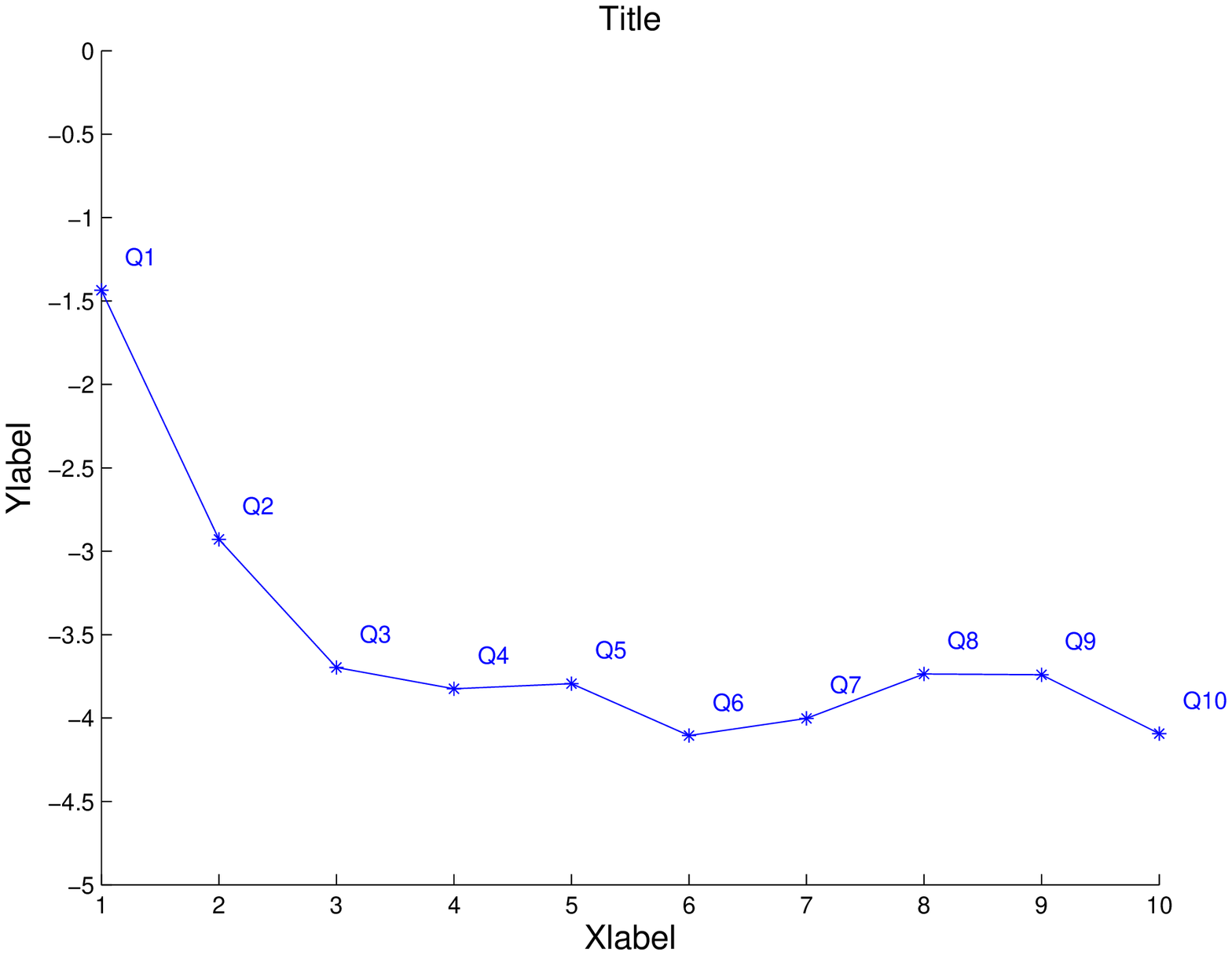}
\end{psfrags}}
\caption{Error on the slope versus the degree of the polynomial space $\PP$ under the same complexity.}
\label{FIG:4bis}
\end{figure}

If we fix the complexity, we can test which degree of the polynomial space $\PP$ can provide the best speed. 
Figure \ref{FIG:4bis} shows this approximation error according to the degree of the polynomial space $\PP$ under the same complexity - quantified by the number of degrees of freedom (DoF) of the finite elements space. 

Under the same complexity, we see on Figure \ref{FIG:4bis} that high degree elements still provide better approximations than $Q_{1}$ elements.
Although very high degree elements always provide better approximations than low degree elements, 
the slopes on Figure \ref{FIG:4bis.c} of the curves for low degrees suggest that $Q_{3}$ elements are a good choice. 
Higher elements increase the computation time for matrix inversion and the gain is not valuable.

\begin{figure}[htp]
\centering
\subfigure[Errors according to the complexity]{
\label{FIG:L}
\includegraphics[angle=0,width=6cm]{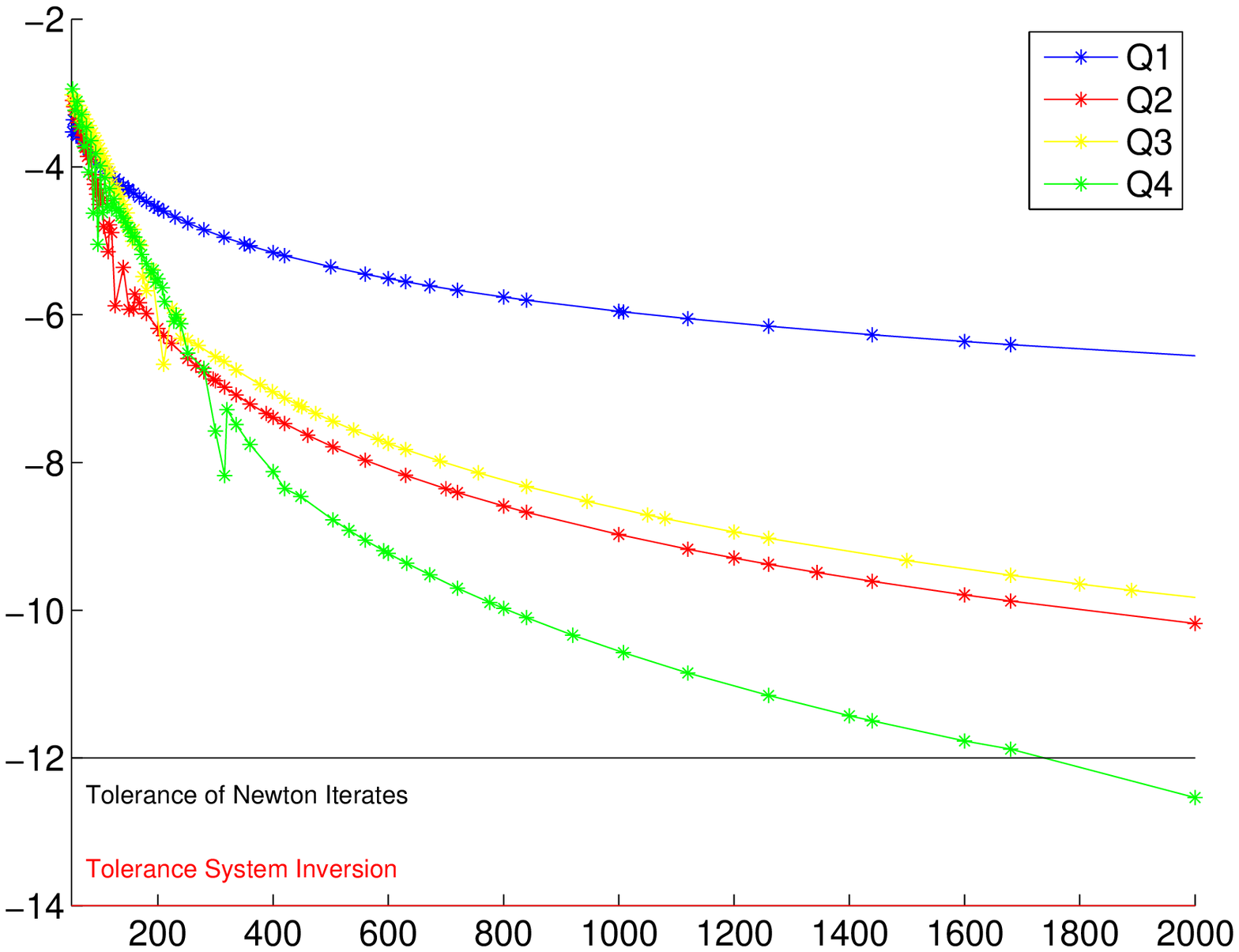}}
\hspace{0.3cm}
\subfigure[Errors according to the logarithm of the complexity]{
\label{FIG:M}
\includegraphics[angle=0,width=6cm]{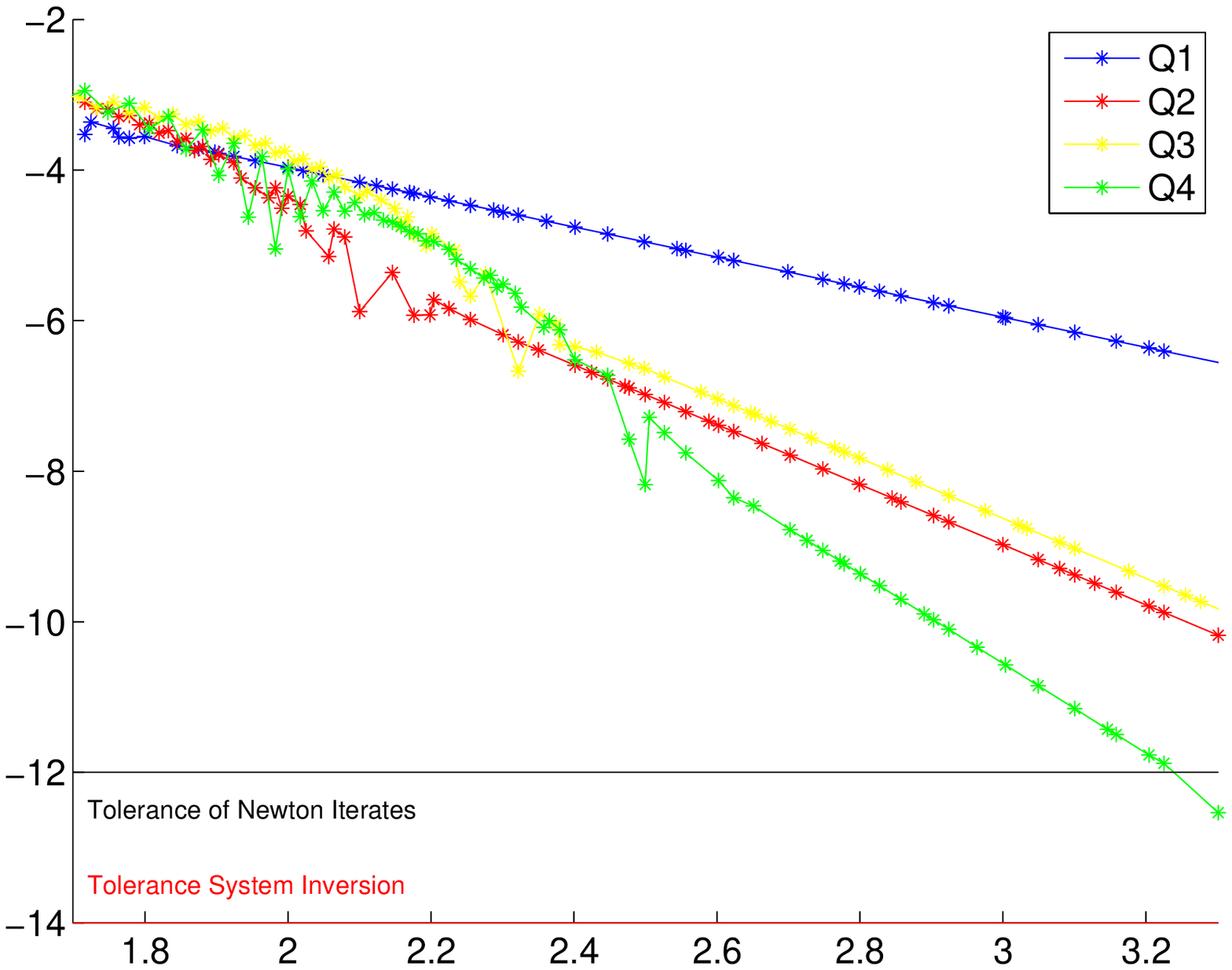}}
\caption{Errors according to the logarithm of the complexity}
\label{FIG:LM}
\end{figure}

Figure \ref{FIG:L} shows the error on the energies according to the complexity under $Q_{1}$, $Q_{2}$, $Q_{3}$ and $Q_{4}$ elements. As for the slopes, we see that the error is decreasing as the number of elements of the mesh is increasing. Whereas the error reaches a $10^{-4}$ precision for the slopes before saturation, the error on the energy reaches the tolerance zone for $Q_{4}$ elements on a mesh with $500$ elements. Figure \ref{FIG:M} shows the logarithm of the error according to the logarithm of the complexity. We see that the evolution is linear for the finest meshes with a good speed. We conclude in particular that we can compute an order of the speed of the convergence. For $Q_{1}$ elements, we find an order $2$, for $Q_{2}$ elements, we find an order $4$, for $Q_{3}$ elements, we find an order $4$ and for $Q_{4}$ elements, we find an order $6$. Note that the error on the energies should be of the order as the $\rm{H}^1$ error. 

\begin{figure}[htp]
\centering
\begin{psfrags}
\psfrag{Title}{}
\psfrag{Xlabel}{}
\psfrag{Ylabel}{}
\psfrag{DL90}{\Large{\color{blue}DoF90}}
\psfrag{DL180}{\Large{\color{red}DoF180}}
\psfrag{DL300}{\Large{\color{green}DoF300}}
\psfrag{DL630}{\Large{\color{black}DoF630}}
\includegraphics[angle=0,width=12cm]{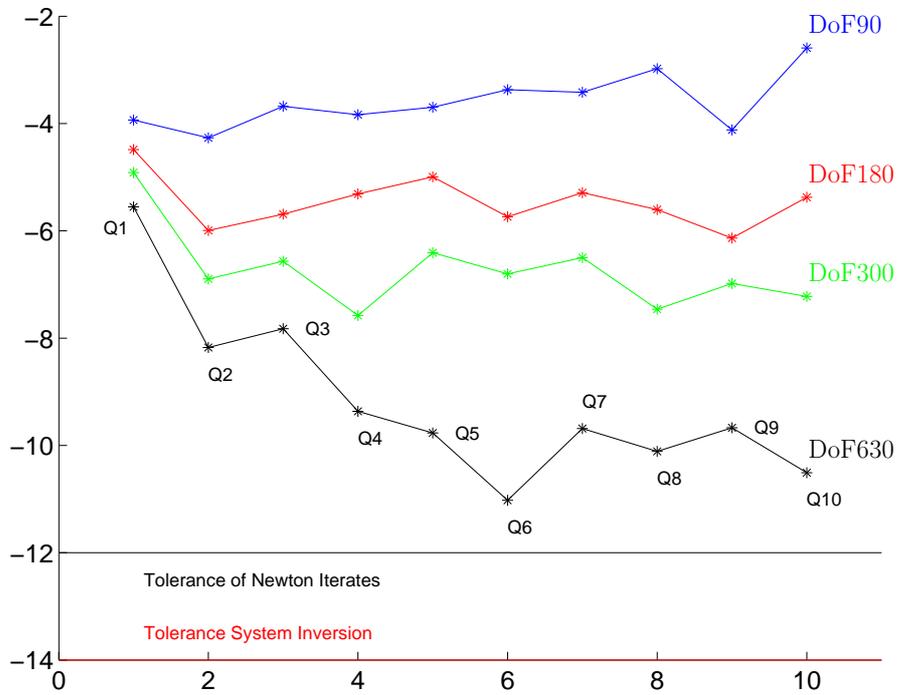}
\caption{Errors on the energy according to the logarithm of the complexity}
\label{FIG:17}
\end{psfrags}
\end{figure}

Now, we fix the complexity and compare the approximation error on the energy according to the degree of the polynomial space $\PP$. For the complexities $90$, $180$, $300$ and $630$, we have drawn the decimal logarithm of the errors on Figure \ref{FIG:17}.

Again, under a same complexity, if we increase the degree of the polynomial space $\PP$, the high degrees can provide better approximation, except on the coarse grids. We can conclude that for a fixed mesh (fine enough), high degrees provide a better approximation. But for each complexity, it seems that we have a saturation because the $Q_{6}$, $Q_{7}$, $Q_{8}$, $Q_{9}$ and $Q_{10}$ elements have almost the same errors. We conclude that we have to use elements with high degrees, but it is not necessary to choose the highest. We have to take into account the computational cost, and the precision of our inverse solver. Indeed, even if the complexity is the same, the finite elements matrices have not the same profil. For instance, the bandwidth of the ``mass'' matrix for $Q_{10}$ elements is much larger than for $Q_{1}$ elements. Figures \ref{FIG:LM} and \ref{FIG:17} indicate that $Q_2$ and $Q_3$ elements are a good compromise to ensure good 
results without increasing the computational cost too much.

\begin{figure}[htp]
\centering
\subfigure[Errors according to the complexity]{
\label{FIG:LL}
\includegraphics[angle=0,width=6cm]{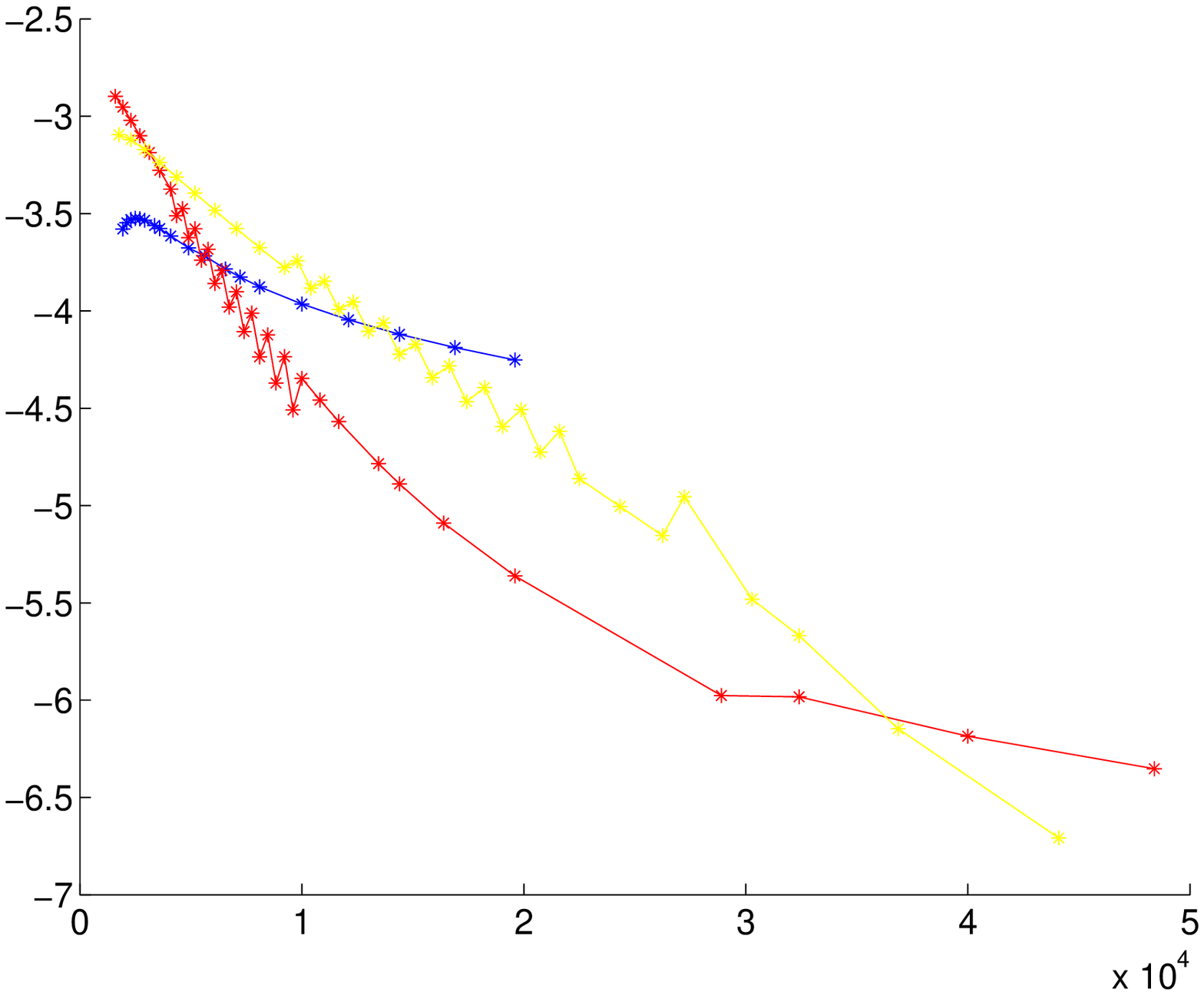}}
\hspace{0.3cm}
\subfigure[Errors according to the logarithm of the complexity]{
\label{FIG:MM}
\includegraphics[angle=0,width=6cm]{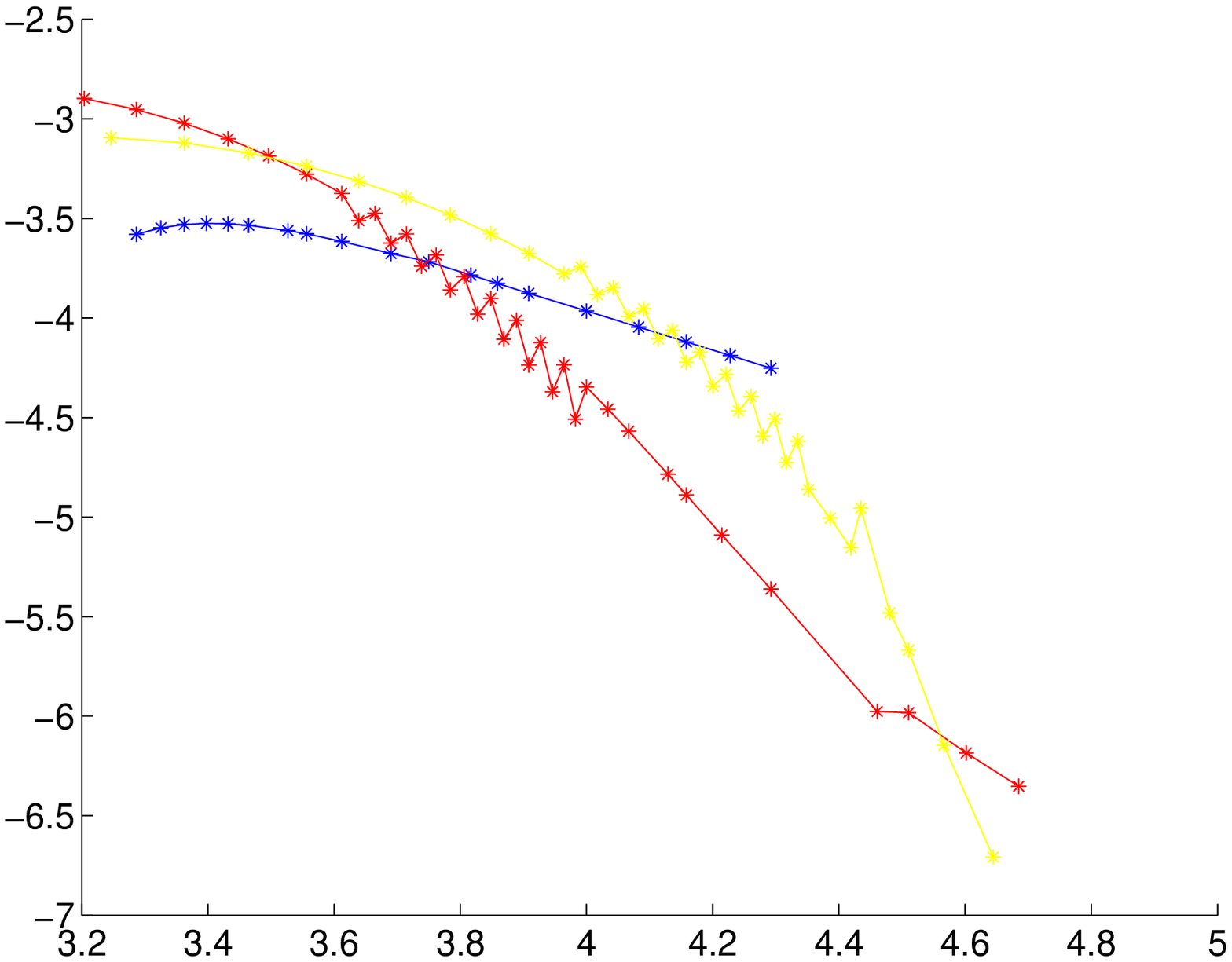}}
\caption{Errors on the energy according to the logarithm of the complexity}
\label{FIG:LLMM}
\end{figure}

In the two dimensional case, the results are drawn on Figure \ref{FIG:LLMM}.  The behaviour is similar.

The energy and the interface are essential physical quantities. From  a mathematical
point of view, it is also important to study the $\rm{L}^2$ error.

\begin{figure}[htp]
\centering
\subfigure[Versus the complexity.]{
\label{FIG:3L2ter.a}
\includegraphics[angle=0,width=6cm]{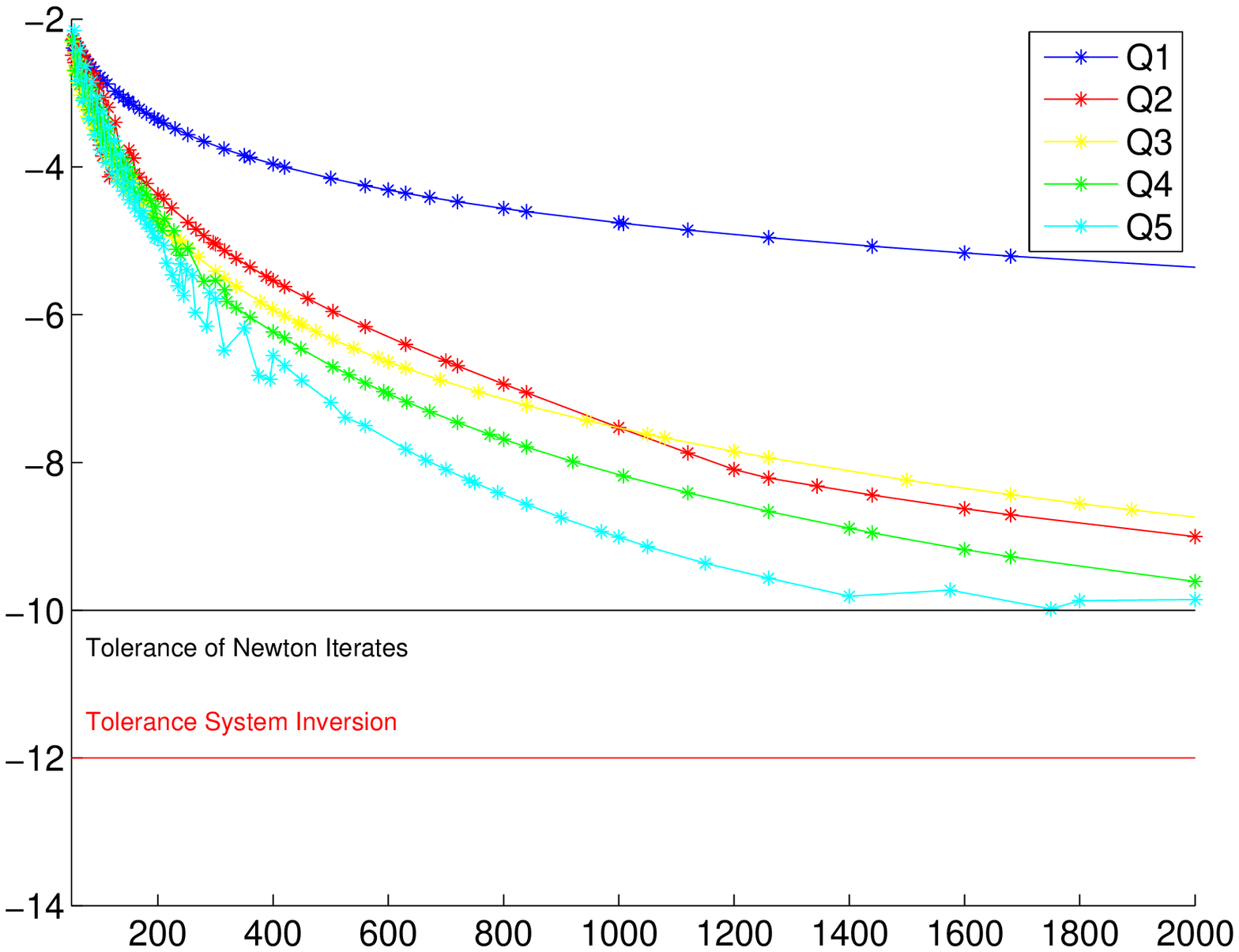}}
\hspace{1cm}
\subfigure[Versus the decimal logarithm of the complexity.]{
\label{FIG:3L2ter.b}
\includegraphics[angle=0,width=6cm]{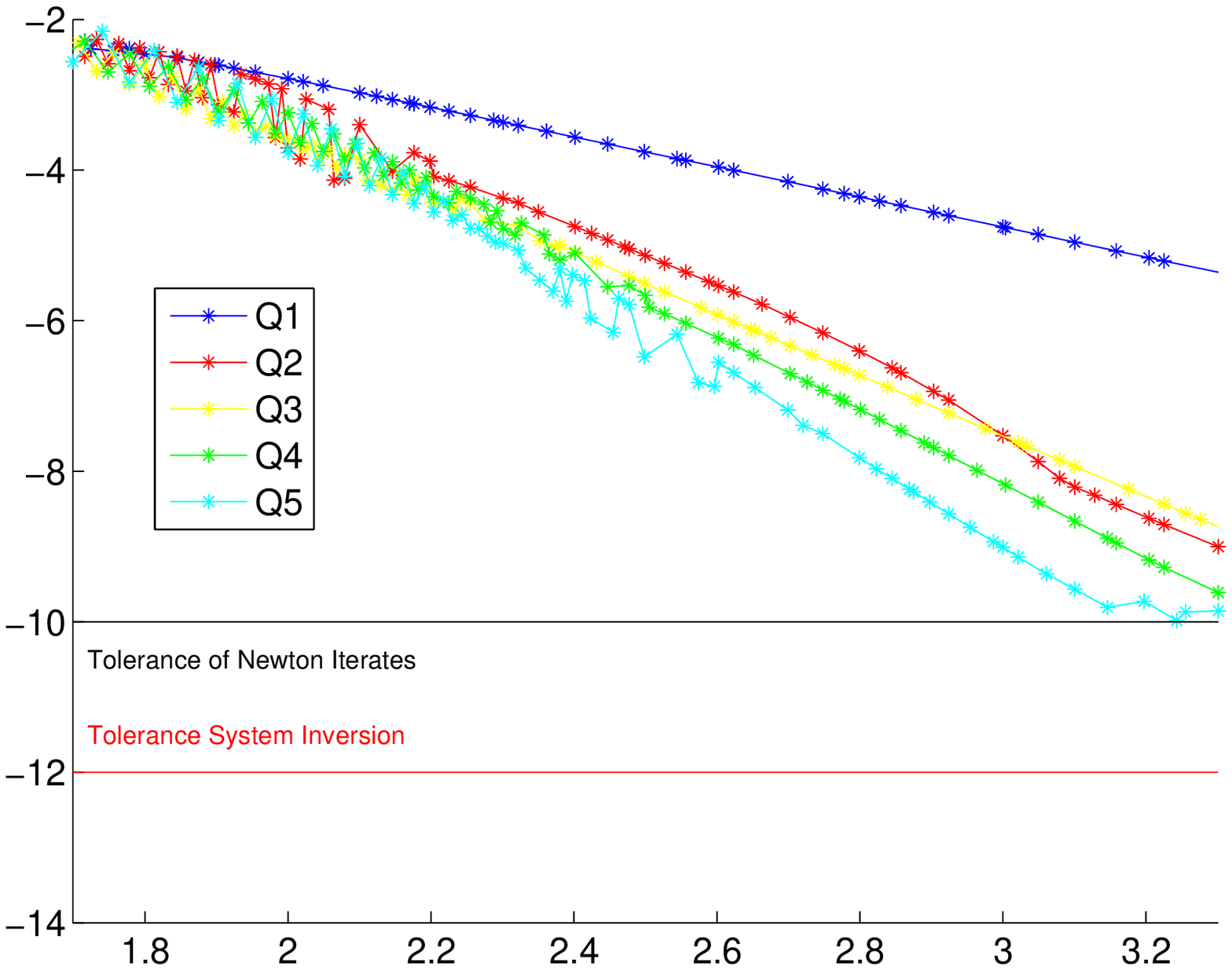}}
\caption{Comparison between the $\rm{L}^2$ errors under the same complexity.}
\label{FIG:3L2ter}
\end{figure}

Figures \ref{FIG:3L2ter.a} and \ref{FIG:3L2ter.b} show the $\rm{L}^2$ error according to the complexity for $Q_{1}$, $Q_{2}$, $Q_{3}$, $Q_{4}$ and $Q_{5}$ elements. On Figure \ref{FIG:3L2ter.b}, we have used a logarithmic scale in order to compare the convergence rate. 
We have computed the order of the speed of the convergence. If we extrapolate the lines, we can find the necessary complexity in order to reach the saturation.
\BD
\begin{array}{cccc}
\text{ Degrees }&\text{ Order }&\text{ Complexity for saturation }&\text{ Grid for saturation }\\
1&1.9960&423829&423829\\
2&3.9682&3558&1779\\
3&4.0216&4110&1370\\
4&4.9119&2364&591\\
5&5.9040&1475&295\\
\end{array}
\ED

Again, $Q_2$ and $Q_3$ elements give very good results for a reasonable computational cost. 
We have decided to prefer $Q_3$ elements because it seems that they provide better results on 
the interface length as shown on Figure \ref{FIG:4bis}.



\section{Stationary states}\label{S:4}
\begin{figure}[htp]
\centering
\includegraphics[angle=0,width=8cm]{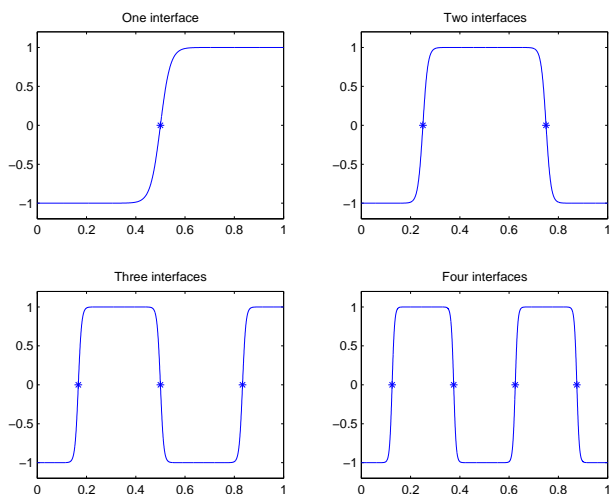}
\caption{Four numerically stable states.}
\label{FIG:9}
\end{figure}

The Cahn-Hilliard equation has a lot of asymptotic equilibria (see \cite{MR1331565}, \cite{MR1263907} and \cite{MR763473}). In the one dimensional case,  a state can be described by the number of interfaces and their positions. On Figure \ref{FIG:9}, we show four states which are numerically stable. It is possible to observe 
more than one interface only for small $\varepsilon$. Only when the interface is very thin - 
{\it i.e.} for small $\varepsilon$, the interfaces do not interact. Note that the energy increases
with the number of the interfaces. 

In fact, this is a bifurcation phenomenon. When $\varepsilon$ crosses critical values,  bifurcations
happen and more stationary solutions appear. 

In \cite{MR1950337}, the authors consider the stationary states of \eqref{Eq:0.6} on the square. They numerically study the solutions of the following semi-linear elliptic equation.
\BE\label{Eq:3.1}
\left\{\begin{array}{ll}
c = u-u^3+\varepsilon^2\Delta u ,&\text{ on } \Omega,\\
\\
\nabla u \cdot \nu = 0, &\text{ on } \partial\Omega,
\end{array}\right.
\EE
together with the mass constraint:
\BE\label{Eq:3.2}
\frac{1}{|\Omega|}\int_{\Omega} u(x) \dd x = m,
\EE
where $\Omega= [0,1]^2$ is the square, $c \in \R$ and $m\in \R$ are parameters. They study stationary solutions under the three-dimensional parameter space $\left(c,m,1/\varepsilon^2\right)$.
For this system and for all $\varepsilon$, a trivial solution is given by the constant solution $u\equiv m$ with $c = m -m^3$. The linearization around $u \equiv m$ of \eqref{Eq:3.1} under the mass constraint reads
\BE\label{Eq:3.3}
\left\{\begin{array}{ll}
0 = \left(1-3m^2\right)u+\varepsilon^2\Delta u ,&\text{ on } \Omega,\\
\\
\nabla u \cdot \nu = 0, &\text{ on } \partial\Omega.
\end{array}\right.
\EE
Let $v_{r}$ be an eigenfunction of the Laplace-Neumann operator in $\Omega$ defined in \eqref{Eq:4.1} with eigenvalue $r\in \R^+$, then $v_{r}$ is also an eigenfunction of \eqref{Eq:3.3} when
\BE\label{Eq:3.4}
\frac{1}{\varepsilon^2} = \frac{r}{1-3m^2} \text{ for } |m| < \frac{1}{\sqrt{3}}.
\EE
But $\sigma_{+} = 1/\sqrt{3}$ for the quartic double-well potential, so this equality shows that bifurcations may occur only for $m$ in the spinodal region. For the square domain $\Omega = [0,1]^2$, the eigenfunctions are:
\BE
v_{r}(x,y)=v_{k,l} (x,y) := \cos(\pi k x) \cos(\pi l y) \quad\text{ for } (x,y) \in [0,1]^2,
\EE
with $(k,l)\in \N^2$ such that $r = (k^2+l^2) \pi^2$. For the mode $v_{1,1}$ (i.e. $r=2\pi^2$), we obtain nontrivial solutions bifurcating at $u\equiv \pm m^*$ with $m^* = \sqrt{\left(1-\varepsilon^2r\right)/3}$. We fix $m=0$, such that the bifurcations occur as $1 =\varepsilon^2 r$.

\begin{figure}[htp]
\centering
\subfigure[Mode $v_{0,1}$]{
\label{FIG:V.a}
\includegraphics[angle=0,width=6cm]{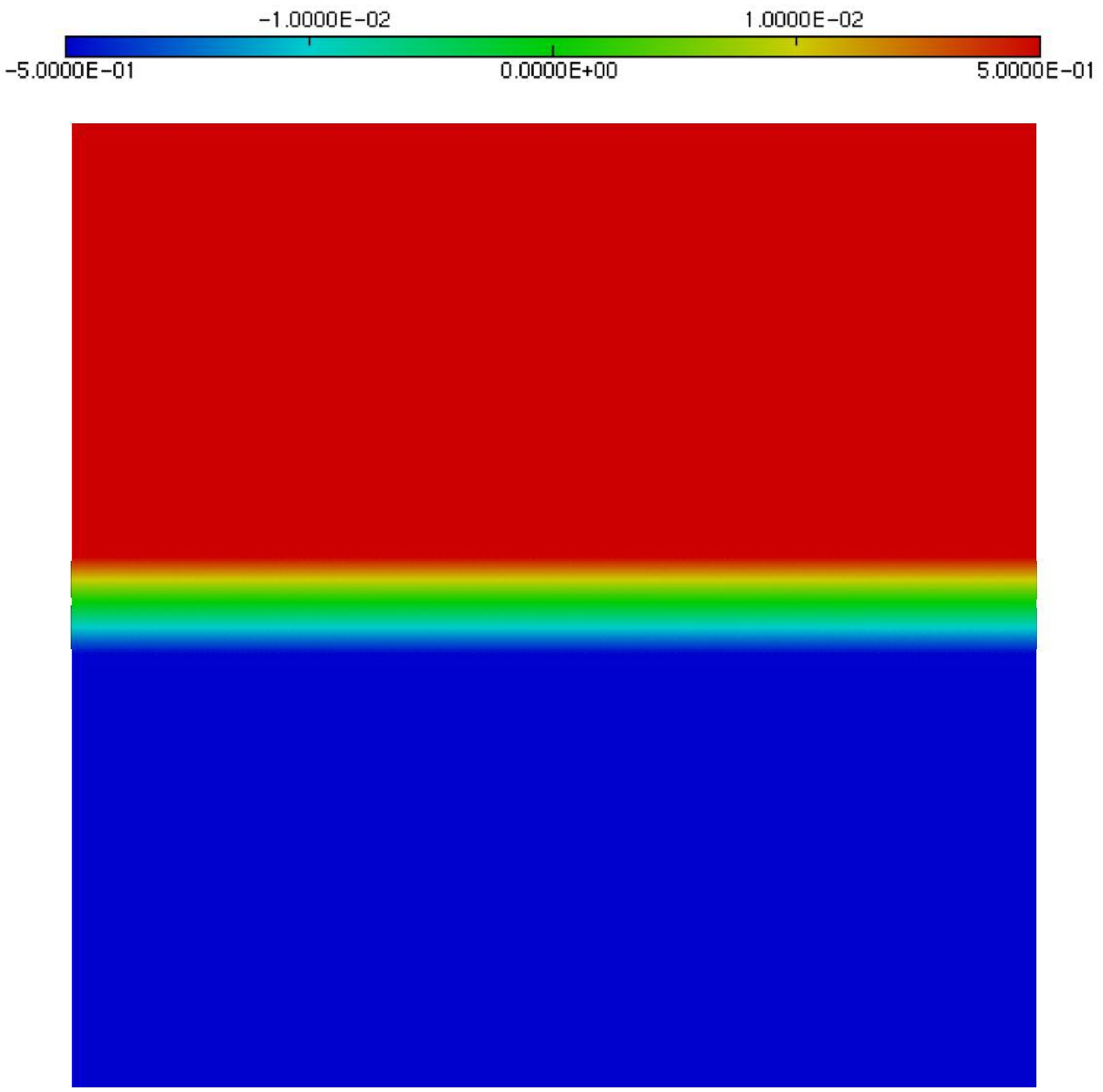}}
\hspace{0.3cm}
\subfigure[Mode $v_{0,1}$+$v_{1,0}$]{
\label{FIG:V.b}
\includegraphics[angle=0,width=6cm]{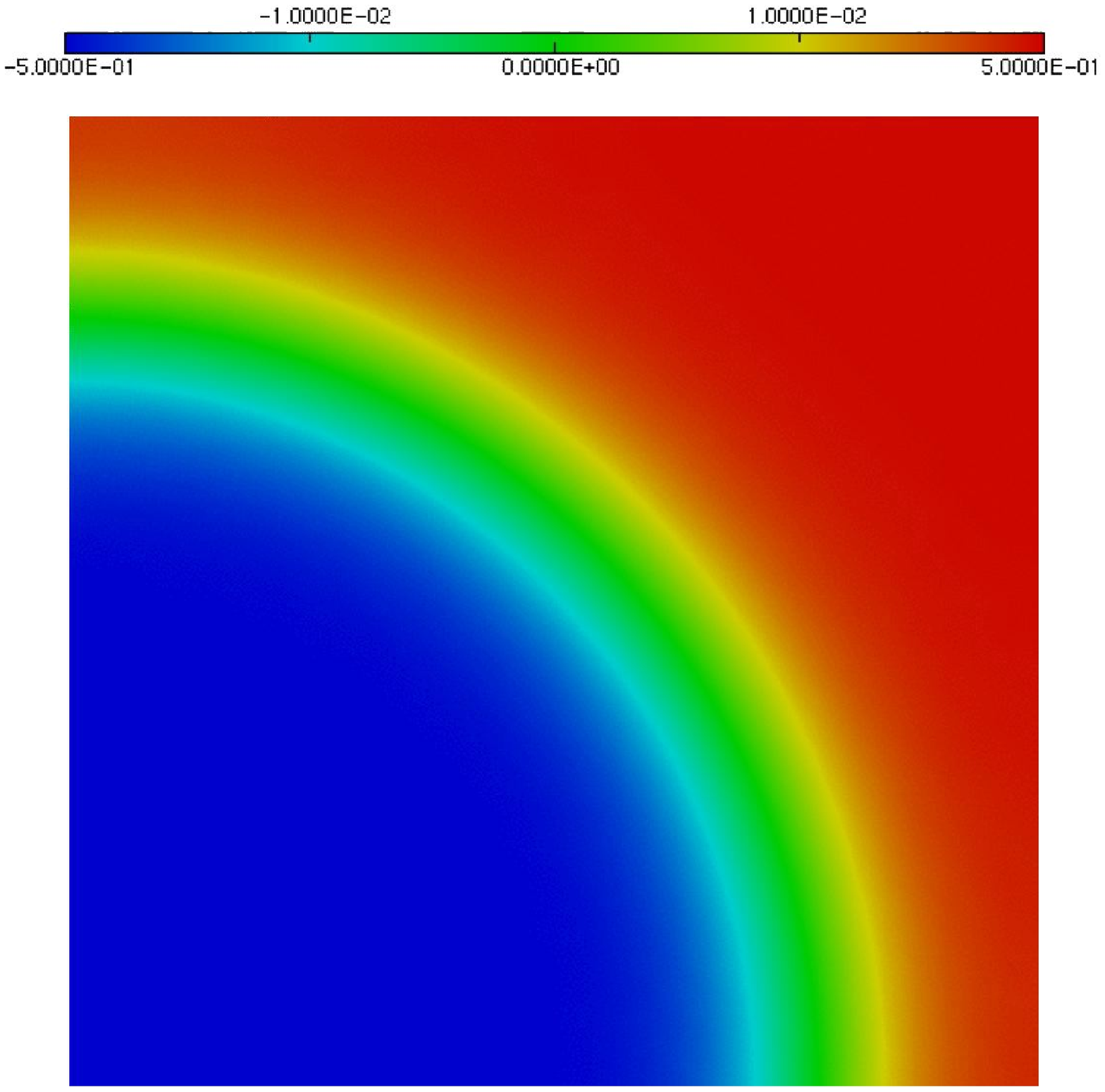}}
\\
\vspace{10pt}
\subfigure[Mode $v_{1,1}$]{
\label{FIG:V.c}
\includegraphics[angle=0,width=6cm]{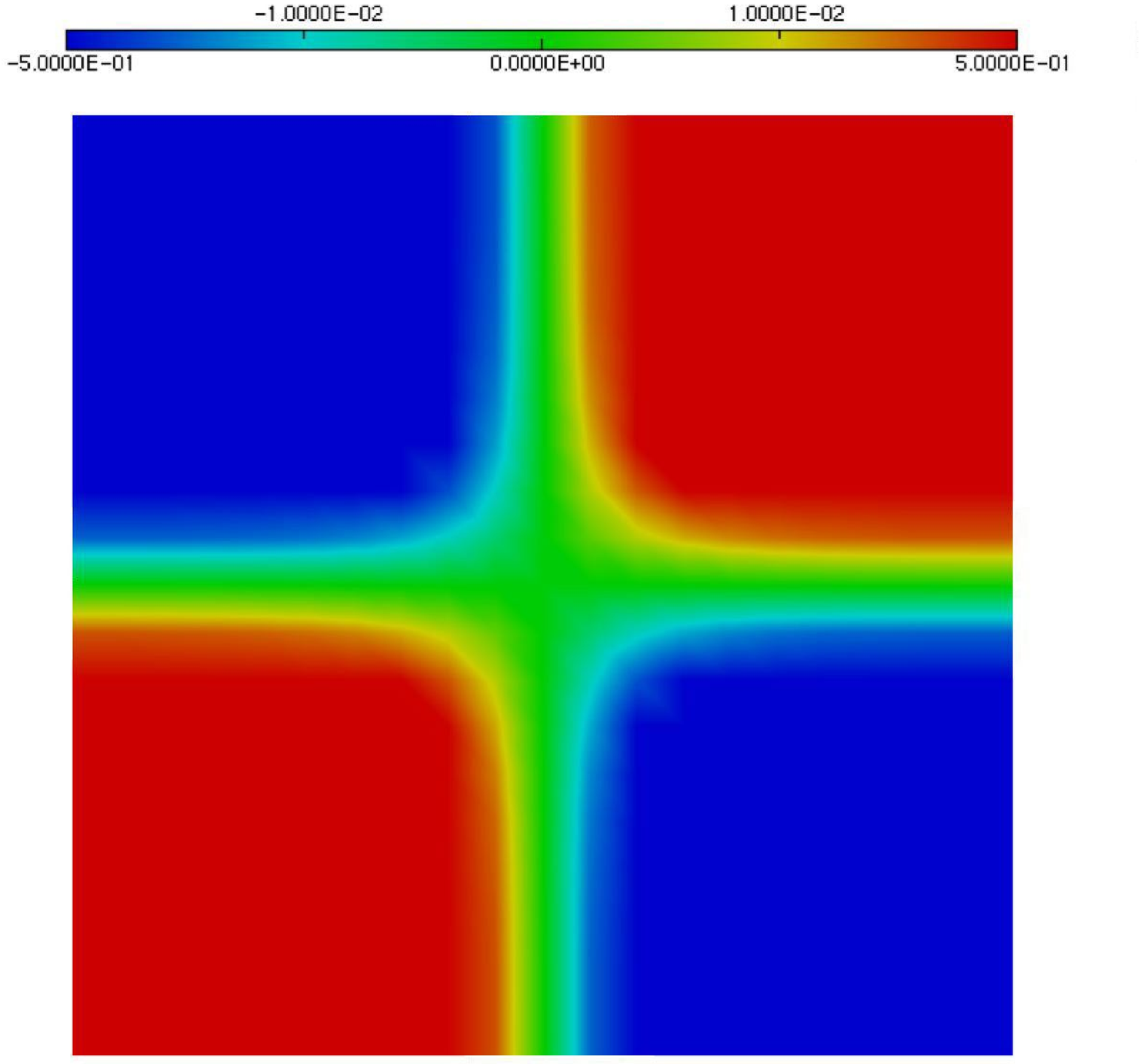}}
\hspace{0.3cm}
\subfigure[Mode $v_{1,4}$+$v_{4,1}$]{
\label{FIG:V.d}
\includegraphics[angle=0,width=6cm]{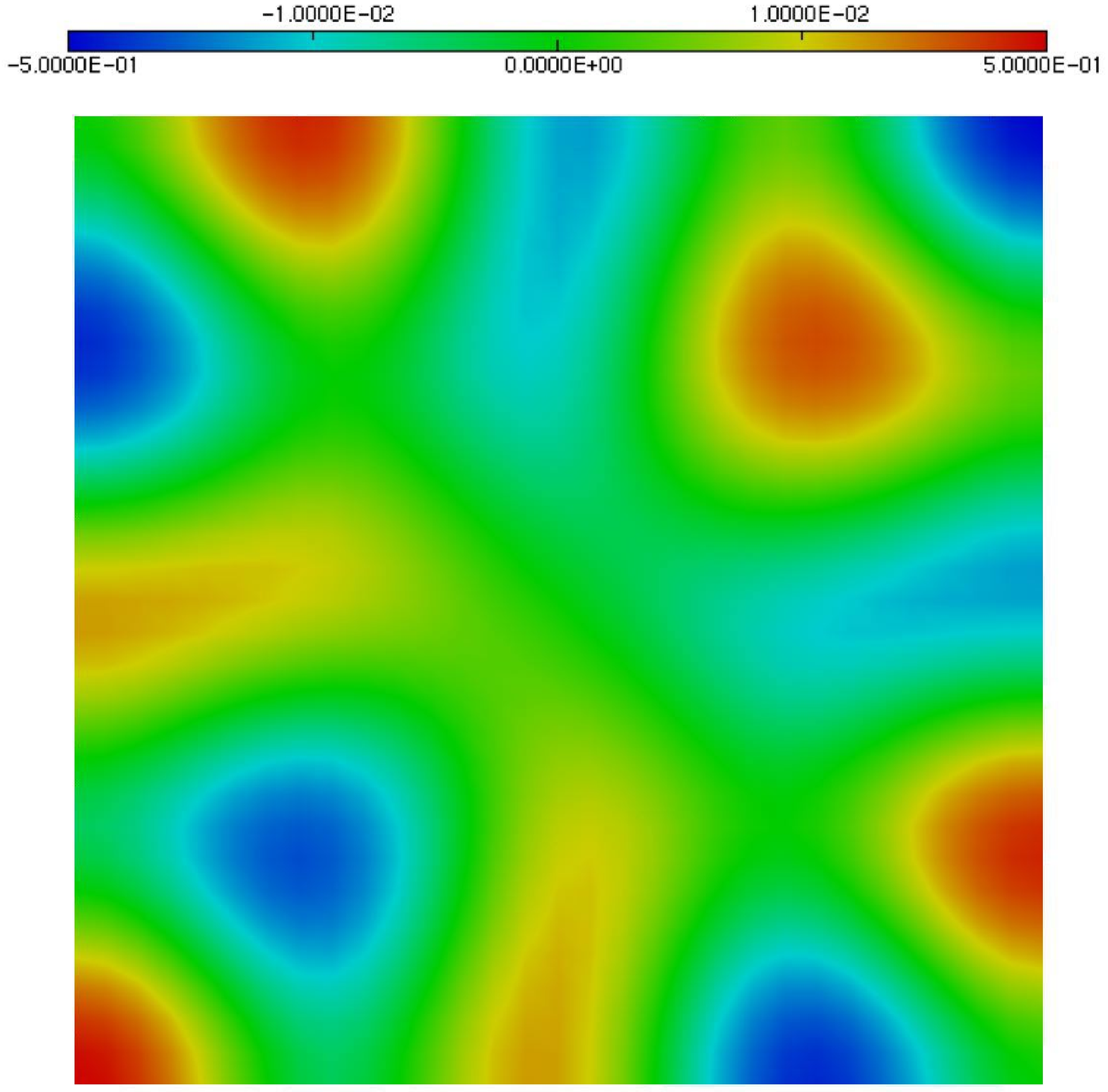}}
\caption{Asymptotic equilibria in the Maier-Paape-Miller nomenclature.}
\label{FIG:V}
\end{figure}

The previous asymptotic equilibria -- described in \cite{MR1950337} -- are asymptotic solutions of the dynamical evolution. For instance, we have obtained the $v_{1,1}$ mode as a stationary solution of a dynamical evolution (See Figure \ref{FIG:V.c}). A random start may lead to different modes, and actually we only see the most stable of them in long time. Figures \ref{FIG:V.a} and \ref{FIG:V.b} show the stable states that we see most of the time. 

All the symmetrical states are also stable. In \cite{MAWANG}, the authors have studied the global attractor on a square and they have proved 
that, after the first bifurcation, there exist 4 minimal attractors (see Theorem 4.2 in \cite{MAWANG}) obtained by symmetrization of Figure \ref{FIG:V.a}. The other stable states 
shown here appear after subsequent bifurcations. Starting the simulation with  well chosen initial data, we have been able to recover dynamically all the stable states described in \cite{MR1950337}. If we choose the mode $v_{4,1}+v_{1,4}$ (which is the last mode studied by Maier-Paape and Miller), we see on Figure~\ref{FIG:V.d} the asymptotic equilibria that we have obtained.

In general, the stable states are deeply dependent on the eigenvalues of the Laplace operator on $\Omega$. Let $\left(\rho_{k}\right)_{k\in\N}$ and $\left(v_{\rho_{k}}\right)_{k\in\N}$ be the eigenvalues and eigenvectors of the following problem:
\BE\label{Eq:4.1}
\left\{\begin{array}{ll}
-\Delta v_{\rho_{k}} = \rho_{k}v_{\rho_{k}} ,&\text{ on } \Omega\subset \R^n,\\
\\
\nabla v_{\rho_{k}} \cdot \nu = 0, &\text{ on } \partial\Omega,\\
\\
\int_{\Omega} v_{\rho_{k}}(\theta) \dd \theta = 0. 
\end{array}\right.
\EE
In \cite{MAWANG}, the authors have studied the bifurcations and the global attractors of the Cahn-Hilliard problem. In their nomenclature, they consider the following Cahn-Hilliard equation:
\BE\label{Eq:3.6}
\left\{\begin{array}{ll}
\partial_{t} v = \Delta w,&\text{ on } \Omega\subset \R^n,\\
\\
w = -\lambda v + \gamma_{2} v^2 + \gamma_{3}v^3-\Delta v ,&\text{ on } \Omega\subset \R^n,\\
\\
\nabla v \cdot \nu = 0 = \nabla w \cdot \nu, &\text{ on } \partial\Omega,\\
\end{array}\right.
\EE
where $\lambda$, $\gamma_{2}$ and $\gamma_{3}$ are parameters. If $u$ is a solution of the system \eqref{Eq:0.6} on $\Omega$, then $v$ is a solution on $\Omega/\sqrt{\varepsilon}$ of \eqref{Eq:3.6} if we define for all $t \in \R$ and $x \in \Omega$:
\BE\label{Eq:3.7}
v : (t,x) \mapsto u(t,x\sqrt{\varepsilon}).
\EE
and the correpondence is given by the following equalities.
\BE\label{Eq:3.8}
\lambda := \frac{1}{\varepsilon}, \quad\gamma_{2} := 0 \text{ and } \gamma_{3} := \frac{1}{\varepsilon}.
\EE
They prove that the first bifurcation occurs as their parameter $\lambda$ is greater than a particular value. For our problem, this bifurcation occurs as $\frac{1}{\varepsilon^2} > \rho_{1}$.

Below, we study this first bifurcation and illustrate theoretical results of \cite{MAWANG}.

\subsection{Asymptotic stable states on a rectangle}\label{S:4.1}

In the case of a rectangular domain $\Omega = [0,2]\times[0,1]$, the hypothesis of the Theorem 4.1 in \cite{MAWANG} holds. Accordingly, if $\frac{1}{\varepsilon^2}>\rho_{1} := \frac{\pi^2}{4}$ then there exist exactly two attractors $\pm u_{\varepsilon}$ which can be expressed as
\BE
\label{e4.1}
\pm u_{\varepsilon}(x,y) = \pm \frac{2\varepsilon}{\sqrt{3}}\sqrt{\frac{1}{\varepsilon^2}-\frac{\pi^2}{4}} \cos \left(\frac{\pi x}{2}\right) + \sqrt{\varepsilon}\ o\left(\left|\frac{1}{\varepsilon^2}-\frac{\pi^2}{4}\right|^{1/2}\right), \quad x \in [0,2], y \in [0,1].
\EE
We define the approximated attractors $\pm v_{\varepsilon}$ by
\BD
\pm v_{\varepsilon}(x) := \pm C(\varepsilon) \sqrt{\frac{1}{\varepsilon^2}-\frac{\pi^2}{4}} \cos \left(\frac{\pi x}{2}\right), \quad x \in [0,2], y \in [0,1],
\ED
where $C(\varepsilon)$ is a constant depending on $\varepsilon$. It is chosen in order to minimize
the $\rm{L}^2$ norm of $u'_\varepsilon - v_\varepsilon$.
 
For multiple values of the parameter $\varepsilon$ around the value $\frac{2}{\pi}$, we have obtained the corresponding numerical stationnary states $u'_{\varepsilon}$. 

\begin{figure}[htp]
\centering
\begin{psfrags}
\psfrag{Title}{}
\psfrag{Xlabel}{}
\psfrag{Ylabel}{}
\includegraphics[angle=0,width=8cm]{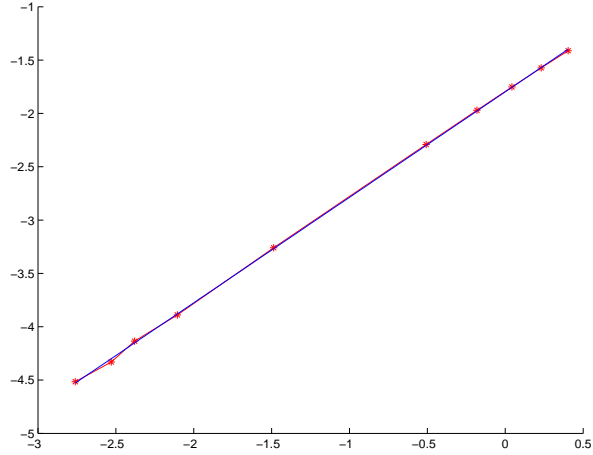}
\caption{Convergence of the ``bifurcationned'' solutions on a rectangular domain.}
\label{FIG:P}
\end{psfrags}
\end{figure}

\begin{figure}
\centering
\begin{psfrags}
\psfrag{Ordre Clambda VS EPS en Log.}{}
\includegraphics[angle=0,width=8cm]{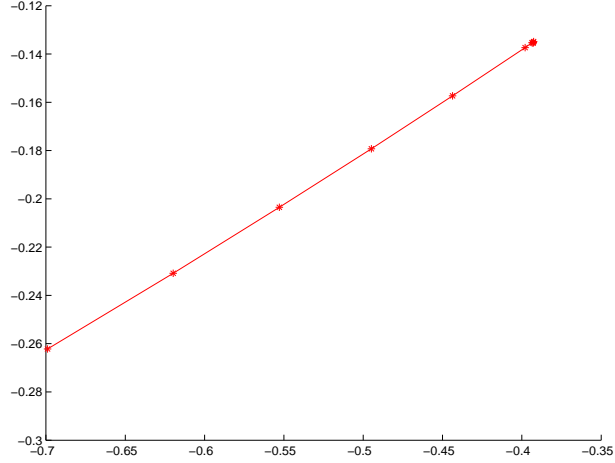}
\caption{Order of convergence of the minimal constant $C(\varepsilon)$.}\label{FIG:T}
\end{psfrags}
\end{figure}

We have checked numerically the validity of formula \eqref{e4.1}. 
We study the following quantity
\BD
\frac{\|u'_{\varepsilon}-v_{\varepsilon}\|_{2}}{\|v_{\varepsilon}\|_{2}} = \frac{\|u'_{\varepsilon}-v_{\varepsilon}\|_{2}}{C(\varepsilon)\sqrt{\left(\frac{1}{\varepsilon^2}-\frac{\pi^2}{4}\right)}}.
\ED
This relative $\rm{L}^2$ norm should converge to zero. On Figure \ref{FIG:P}, we have plotted the decimal logarithm of this relative $\rm{L}^2$ norm according to the decimal logarithm of $\frac{1}{\varepsilon^2} -\frac{\pi^2}{4}$.
Figure \ref{FIG:P} is in conformity with the expecting theoretical results. We can see that the relative $\rm{L}^2$ error converges to $0$ as $\frac{1}{\varepsilon^2}$ converges to $\frac{\pi^2}{4}$. We even can improve formula \eqref{e4.1} and find the exponent $\alpha_{rectangular}$ such that 
\BD
\| u'_{\varepsilon} - v_{\varepsilon} \|_{2} \sim \tilde{C} \left|\frac{1}{\varepsilon^2}-\frac{\pi^2}{4}\right|^{\alpha_{rectangular}}
\ED
where $\tilde{C}$ is an unknown constant. We find that the exponent $\alpha_{rectangular} = 1/2+0.98908$, almost $3/2$.
Moreover, since we know explicitly the attractor, we can verify that our minimal constant $C(\varepsilon)$ is near $\frac{2\varepsilon}{\sqrt{3}}$. 
On Figure \ref{FIG:T}, we have drawn the logarithm of our minimal constant $C(\varepsilon)$ according to the logarithm of $\varepsilon$.

We find
\BE
C(\varepsilon) \sim 1.0682\ *\ \varepsilon^{0.83603} 
\EE
this is in conformity with the fact that $C(\varepsilon)$ converges to $\lim_{\varepsilon\rightarrow\frac{2}{\pi}}\frac{2\varepsilon}{\sqrt{3}} = \frac{4}{\pi\sqrt{3}}$.

The segment may be seen as a degenerate rectangle. 
On the segment $[0,1]$, the eigenvalues of \eqref{Eq:4.1} are $\rho_{k}:=k^2\pi^2$ for all $k\in\N$. According to Theorem 4.2 in \cite{MAWANG}, there is a bifurcation at $\frac{1}{\varepsilon^2} > \pi^2$. Moreover, Remark 4.2 in \cite{MAWANG} states that there exist two minimal attractors $\pm u_{\varepsilon}$ which can be expressed as 
\BE\label{Eq:EquivalentSegment}
\pm u_{\varepsilon}(x) = \pm C(\varepsilon) \sqrt{\left(\frac{1}{\varepsilon^2}-\pi^2\right)} \cos \left(\pi x\right) + \sqrt{\varepsilon}\ o\left(\left|\frac{1}{\varepsilon^2}-\pi^2\right|^{1/2}\right), \quad x \in [0,1],
\EE
where $C(\varepsilon)$ is a constant which can depend on $\varepsilon$. Again, we define the approximated attractors $\pm v_{\varepsilon}$ by
\BD
\pm v_{\varepsilon}(x) := \pm C(\varepsilon) \sqrt{\left(\frac{1}{\varepsilon^2}-\pi^2\right)} \cos \left(\pi x\right), \quad x \in [0,1].
\ED
For multiple values of the parameter $\varepsilon$ around the value $\frac1\pi$, we have obtained the corresponding numerical stationary states $u'_{\varepsilon}$. We choose the constant $C(\varepsilon)$ in order to minimize the $\rm{L}^2$ norm of $u'_{\varepsilon}-v_{\varepsilon}$ and  study the convergence of $u'_{\varepsilon}$ to $v_{\varepsilon}$
using the quantity
\BD
\frac{\|u'_{\varepsilon}-v_{\varepsilon}\|_{2}}{\|v_{\varepsilon}\|_{2}} = \frac{2\|u'_{\varepsilon}-v_{\varepsilon}\|_{2}}{C(\varepsilon)\sqrt{\left(\frac{1}{\varepsilon^2}-\pi^2\right)}}.
\ED
\begin{figure}[htp]
\centering
\begin{psfrags}
\psfrag{Title}{}
\psfrag{Xlabel}{}
\psfrag{Ylabel}{}
\includegraphics[angle=0,width=8cm]{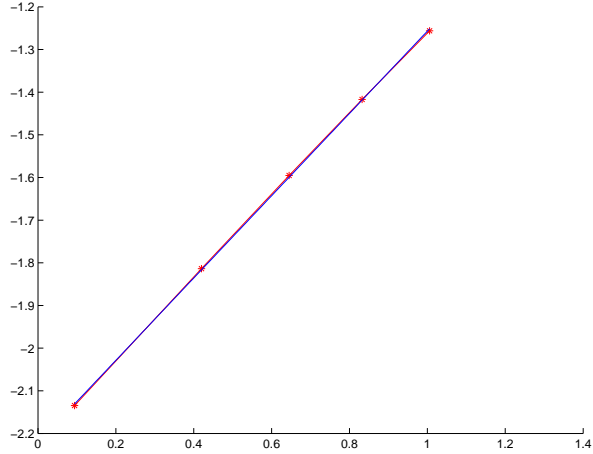}
\caption{Convergence of the ``bifurcationned'' solutions on a segment.}
\label{FIG:N}
\end{psfrags}
\end{figure}

Figure \ref{FIG:N} is in conformity with the expected theoretical results. We can see that the relative $\rm{L}^2$ error converges to $0$ as $\frac{1}{\varepsilon^2}$  converges to $\pi^2$
and find the exponent $\alpha_{segment}$ such that
\BD
\|u'_{\varepsilon} -  v_{\varepsilon} \|_{2} \sim \tilde{C} \left|\frac{1}{\varepsilon^2}-\pi^2\right|^{\alpha_{segment}}
\ED
where $\tilde{C}$ is an unknown constant. We have found $\alpha_{segment}=1/2 + 1.0254$, again almost
$3/2$.

\begin{figure}[htp]
\centering
\begin{psfrags}
\psfrag{Ordre Clambda VS sqrt(EPS) en Log.}{}
\includegraphics[angle=0,width=8cm]{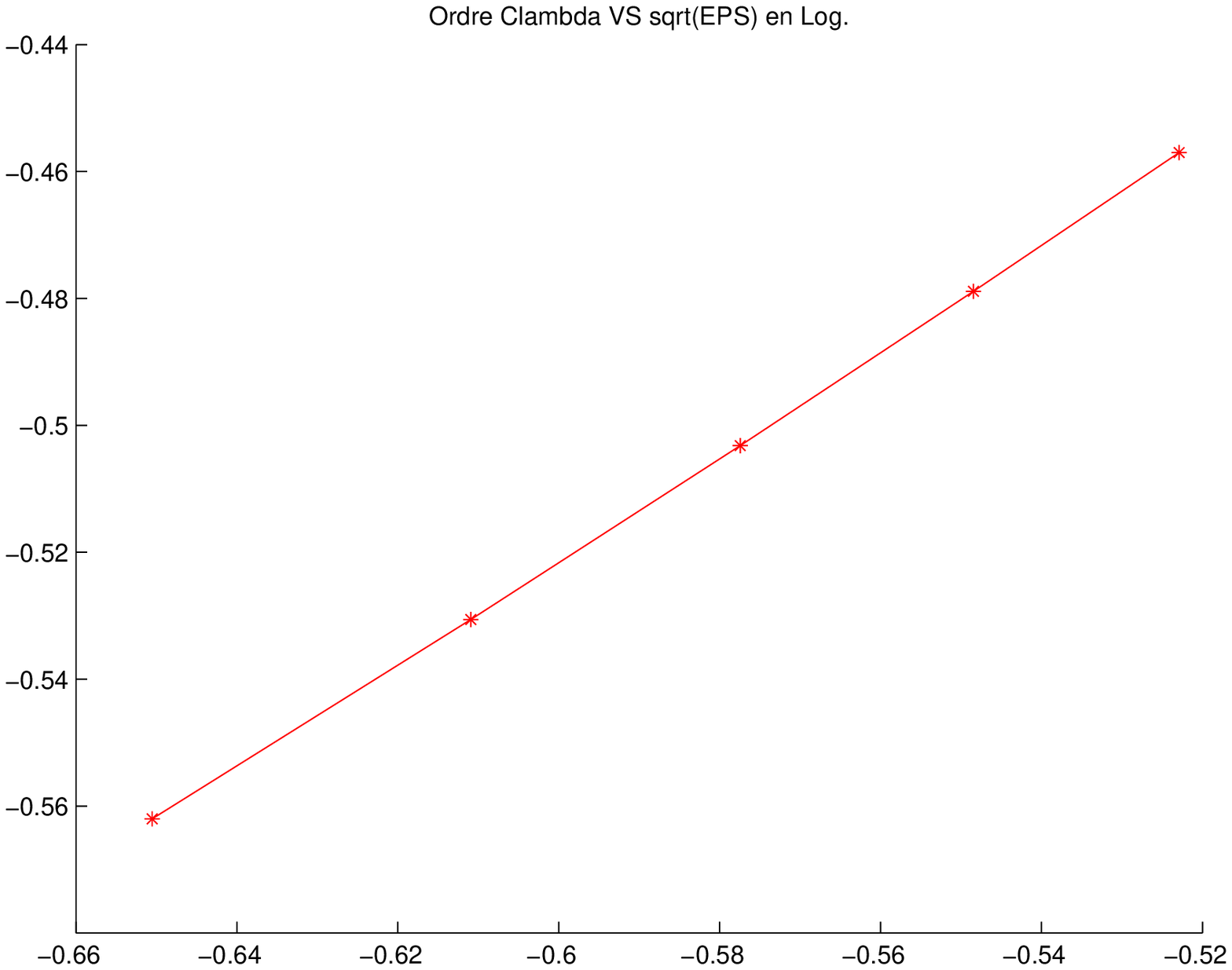}
\caption{Order of convergence of the minimal contant $C(\varepsilon)$.}
\label{FIG:NN}
\end{psfrags}
\end{figure}

We have said that the constants $C(\varepsilon)$ have been numerically chosen in order to minimize the $\rm{L}^2$ norm of $u'_{\varepsilon}-v_{\varepsilon}$. If we extend the results of \cite{MAWANG}, we expect that the constant $C(\varepsilon) \sim \frac{2}{\sqrt{3}}\varepsilon$. 
Thus, on Figure \ref{FIG:NN}, we have drawn the logarithm of our minimal constant $C(\varepsilon)$ according to the logarithm of $\varepsilon$.
If we study the slope, we find
\BE
C(\varepsilon) \sim 1.0844\ *\ \varepsilon^{0.94594}
\EE
which is again in conformity with the theoretical formula.


\subsection{Asymptotic stable states on smooth domains}\label{S:4.3}

\begin{figure}[htp]
\begin{center}
\includegraphics[angle=0,width=8cm]{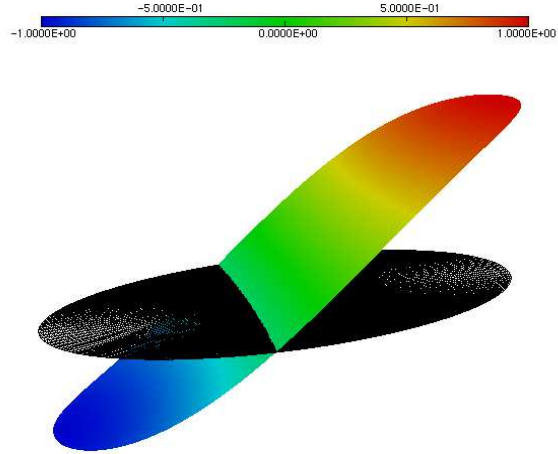}
\caption{First eigenvector on the ellipse.}\label{FIG:23}
\end{center}
\end{figure}

\begin{figure}[htp]
\begin{center}
\includegraphics[angle=0,width=8cm]{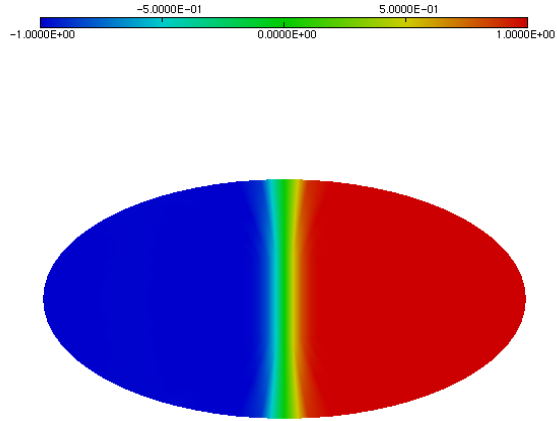}
\caption{Steady state on the ellipse.}\label{FIG:25}
\end{center}
\end{figure}

\begin{figure}[htp]
\centering
\begin{psfrags}
\psfrag{Title}{}
\psfrag{Xlabel}{}
\psfrag{Ylabel}{}
\includegraphics[angle=0,width=8cm]{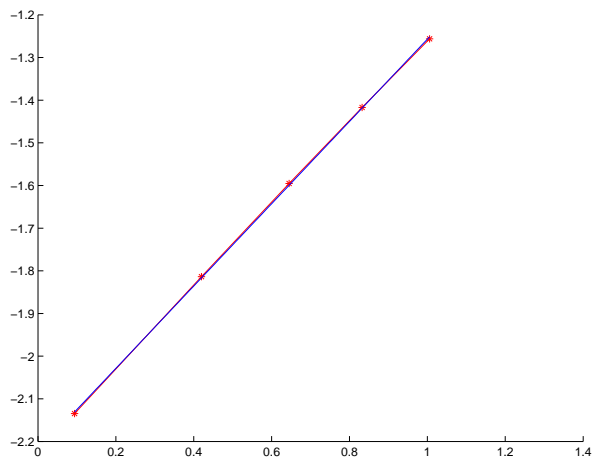}
\caption{Convergence of the ``bifurcationned'' solutions on the ellipse.}
\label{FIG:O}
\end{psfrags}
\end{figure}

For a smooth domain $\Omega$, the hypothesis of Theorem 3.1 in \cite{MAWANG} holds. 
We have considered an ellipse. On the ellipse, the first eigenvalue $\rho_{1}\simeq 0.8776$ is simple. On Figure \ref{FIG:23}, we have drawn the corresponding first eigenvector.
If $\frac{1}{\varepsilon^2} > \rho_{1}$, the problem \eqref{Eq:3.6} has two steady states $\pm u_{\varepsilon}$ which can be expressed as 
\BD
\pm u_{\varepsilon} = \pm C(\varepsilon) \sqrt{\left(\frac{1}{\varepsilon^2}-\rho_{1}\right)} v_{\rho_{1}} + \sqrt{\varepsilon}\ o\left(\left|\frac{1}{\varepsilon^2}-\rho_{1}\right|^{1/2}\right),
\ED
where $C(\varepsilon)$ is a constant which can depend on $\varepsilon$. We define the approximated attractors $\pm v_{\varepsilon}$ by
\BD
\pm v_{\varepsilon} := \pm C(\varepsilon) \sqrt{\left(\frac{1}{\varepsilon^2}-\rho_{1}\right)} v_{\rho_{1}},
\ED
where $v_{\rho_{1}}$ is a fixed eigenvector.
On the Figure \ref{FIG:25}, we have drawn a steady states.
For multiple values of the parameter $\varepsilon$ around the value $\frac{1}{\sqrt{\rho_{1}}}$, we have obtained the corresponding numerical stationnary states $u'_{\varepsilon}$. As in section \ref{S:4.1}, we choose the constant $C(\varepsilon)$ in order to minimize the $\rm{L}^2$ norm of $u'_{\varepsilon}-v_{\varepsilon}$ and study the convergence of $u'_{\varepsilon}$ to $v_{\varepsilon}$. We consider the quantity
\BD
\frac{\|u'_{\varepsilon}-v_{\varepsilon}\|_{2}}{\|v_{\varepsilon}\|_{2}} = \frac{\|u'_{\varepsilon}-v_{\varepsilon}\|_{2}}{C(\varepsilon)\sqrt{\left(\frac{1}{\varepsilon^2}-\rho_{1}\right)}\|v_{\rho_{1}}\|_{2}}.
\ED
According to theorem 3.1, this relative $\rm{L}^2$ norm should converge to zero. 
On Figure \ref{FIG:O}, we have drawn the decimal logarithm of this relative $\rm{L}^2$ norm according to the decimal logarithm of $\frac{1}{\varepsilon^2} - \rho_{1}$.

Figure \ref{FIG:O} corroborates the expected theoretical results. 
We can see that the relative $\rm{L}^2$ error converges to $0$ as $\frac{1}{\varepsilon^2}$  converges to $\rho_{1}$. Then we compute the exponent $\alpha_{ellipse}$ such that
\BD
\|u'_{\varepsilon} -  v_{\varepsilon} \|_{2} \sim \tilde{C} \left|\frac{1}{\varepsilon^2}-\rho_{1}\right|^{\alpha_{ellipse}}
\ED
where $\tilde{C}$ is an unknown constant. We find that exponent $\alpha_{ellipse}=1/2 + 0.9580$.

\subsection{Asymptotic stable states on a trapezoid}\label{S:4.4}

\begin{figure}
\begin{center}
\includegraphics[angle=0,width=8cm]{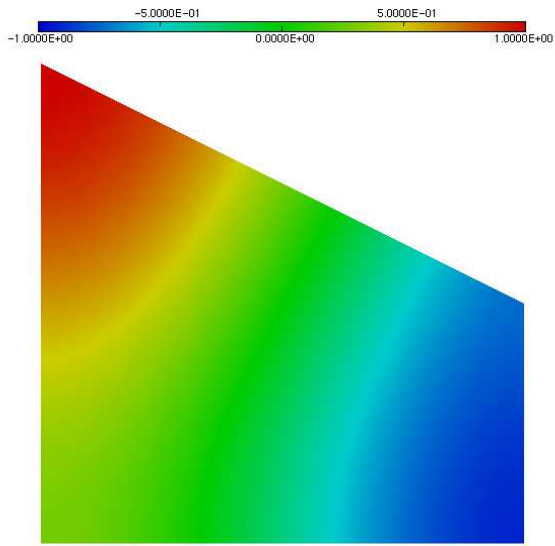}
\caption{Eigenvector of the first eigenvalue on the trapezoid.}\label{FIG:R}
\end{center}
\end{figure}

We try to see if the results of \cite{MAWANG} extend to non smooth domains. We have tested a trapezoid where the first eigenvalue $\rho_{1} \simeq 2.2417$ is simple. On Figure \ref{FIG:R}, we have drawn the corresponding first eigenvector.

\begin{figure}
\begin{center}
\includegraphics[angle=0,width=8cm]{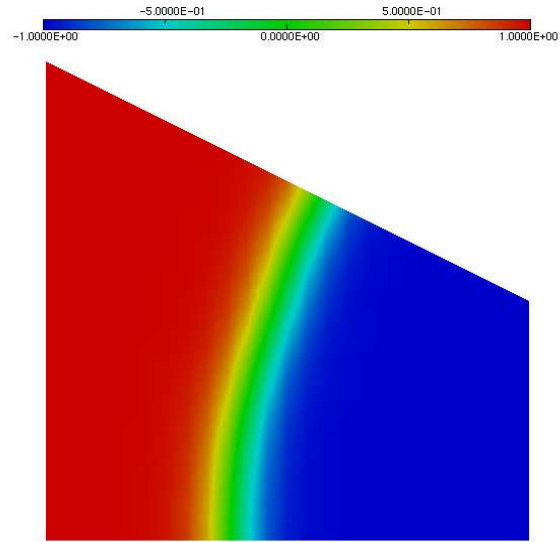}
\caption{Steady state on the trapezoid.}\label{FIG:S}
\end{center}
\end{figure}

The two steady states $\pm u_{\varepsilon}$ should be expressed as 
\BD
\pm u_{\varepsilon}(x) = \pm C(\varepsilon) \sqrt{\left(\frac{1}{\varepsilon^2}-\rho_{1}\right)} v_{\rho_{1}} + \sqrt{\varepsilon}\ o\left(\left|\frac{1}{\varepsilon^2}-\rho_{1}\right|^{1/2}\right), \quad x \in [0,1],
\ED
where $C(\varepsilon)$ is a constant which can depend on $\varepsilon$. 
We define the approximated attractors $\pm v_{\varepsilon}$ by
\BD
\pm v_{\varepsilon}(x) := \pm C(\varepsilon) \sqrt{\left(\frac{1}{\varepsilon^2}-\rho_{1}\right)} v_{\rho_{1}} , \quad x \in [0,1],
\ED
where $v_{\rho_{1}}$ is a fixed eigenvector.
On Figure \ref{FIG:S}, we have drawn a steady state.

\begin{figure}[htp]
\centering
\begin{psfrags}
\psfrag{Title}{}
\psfrag{Xlabel}{}
\psfrag{Ylabel}{}
\includegraphics[angle=0,width=8cm]{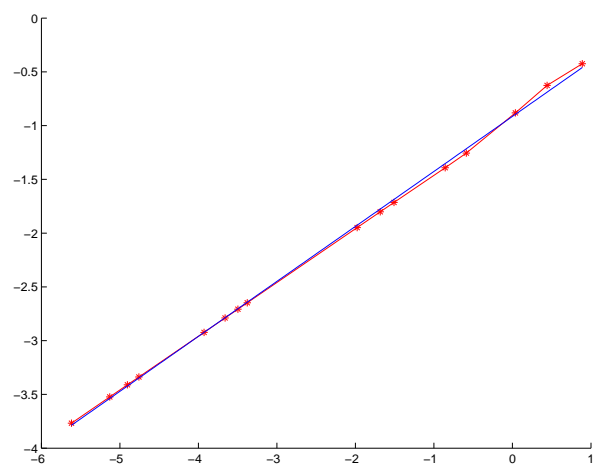}
\caption{Convergence of the ``bifurcationned'' solutions on the trapezoid domain.}
\label{FIG:Q}
\end{psfrags}
\end{figure}

For multiple values of the parameter $\varepsilon$ around the value $\frac{1}{\sqrt{\rho_{1}}}$, we have obtained the corresponding numerical stationnary states $u'_{\varepsilon}$. 
As in section \ref{S:4.1}, if we choose the constant $C(\varepsilon)$ in order to minimize the $\rm{L}^2$ norm of $u'_{\varepsilon}-v_{\varepsilon}$, we can study the convergence of $u'_{\varepsilon}$ to $v_{\varepsilon}$. We have found that
\BD
\frac{\|u'_{\varepsilon}-v_{\varepsilon}\|_{2}}{\|v_{\varepsilon}\|_{2}} = \frac{\|u'_{\varepsilon}-v_{\varepsilon}\|_{2}}{C(\varepsilon)\sqrt{\left(\frac{1}{\varepsilon^2}-\rho_{1}\right)}\|v_{\rho_{1}}\|_{2}}
\ED
does not converge to $0$. It seems that the bifurcation is different in this case. 
On Figure \ref{FIG:Q}, we have drawn the decimal logarithm of
${\|u'_{\varepsilon}-v_{\varepsilon}\|_{2}}$ according to the decimal logarithm of $\frac{1}{\varepsilon^2} - \rho_{1}$
We find that 
\BD
\|u'_{\varepsilon} -  v_{\varepsilon} \|_{2} \sim \tilde{C} \left|\frac{1}{\varepsilon^2}-\rho_{1}\right|^{\alpha_{trapezoid}},
\ED
with $\alpha_{trapezoid}=0.49937$. Thus this difference is of the same order as each term.





\begin{figure}[htp]
\begin{center}
\includegraphics[angle=0,width=6cm]{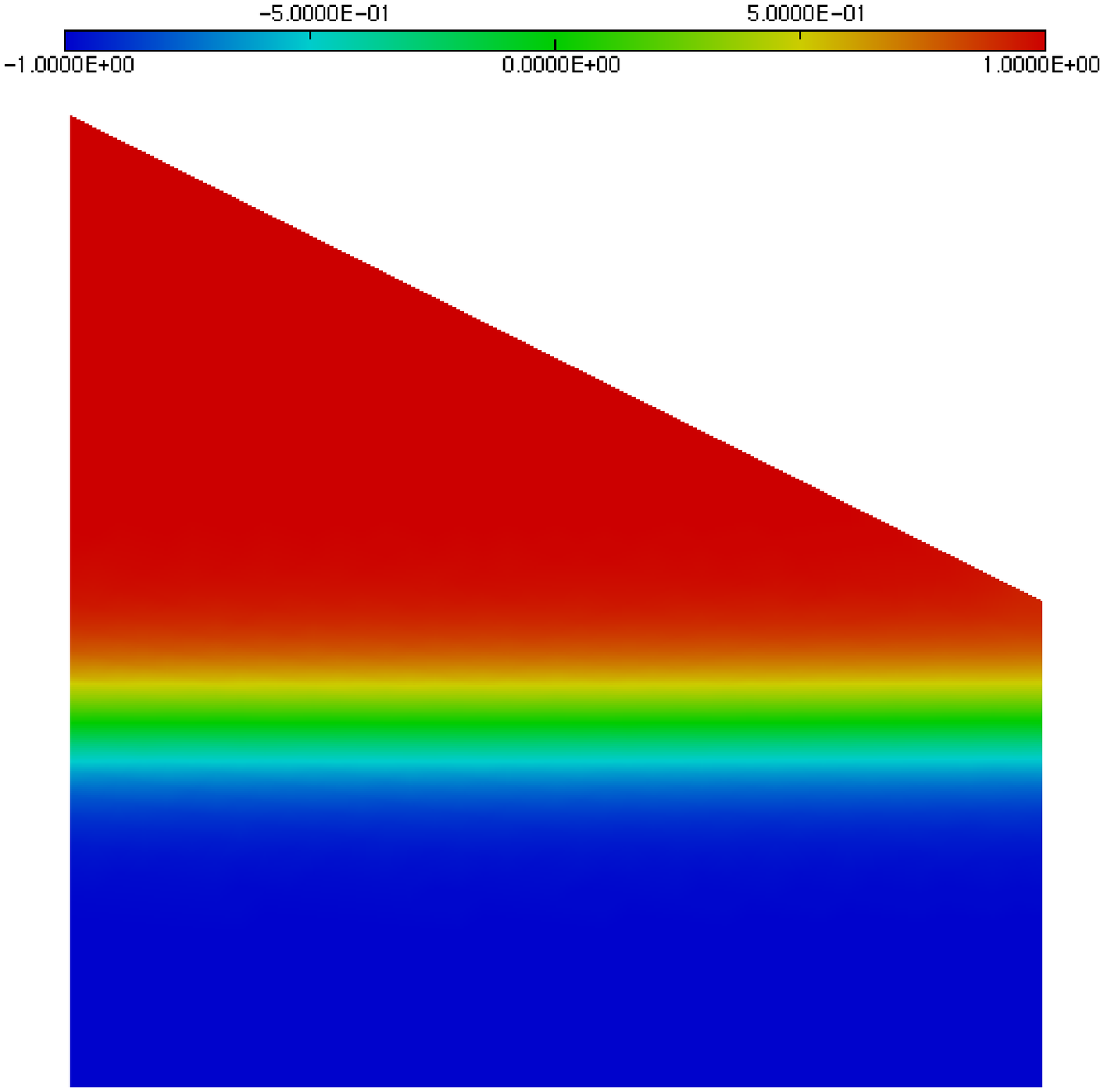}
\hspace{1cm}
\includegraphics[angle=0,width=6cm]{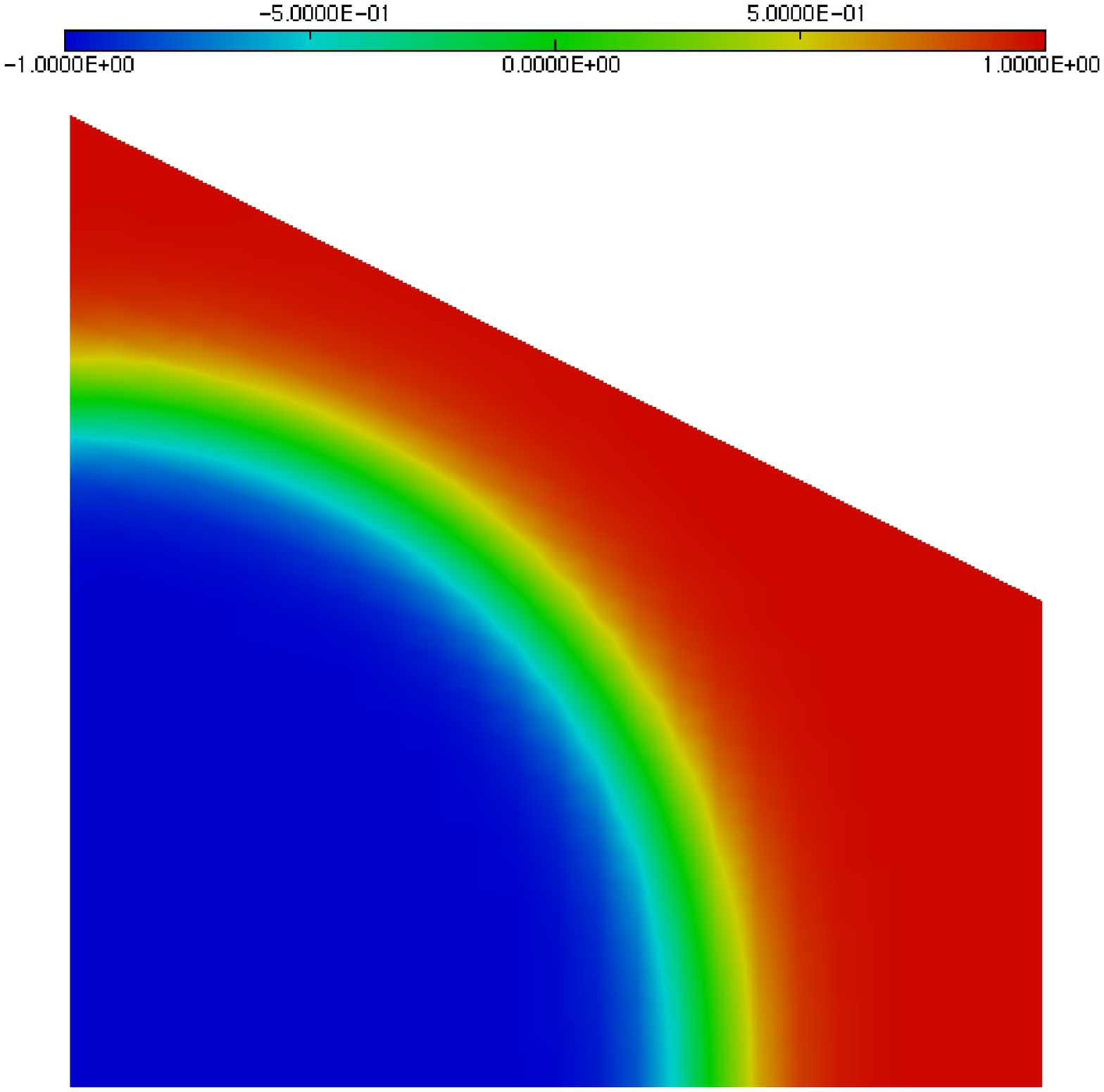}
\vspace{10pt}
\includegraphics[angle=0,width=6cm]{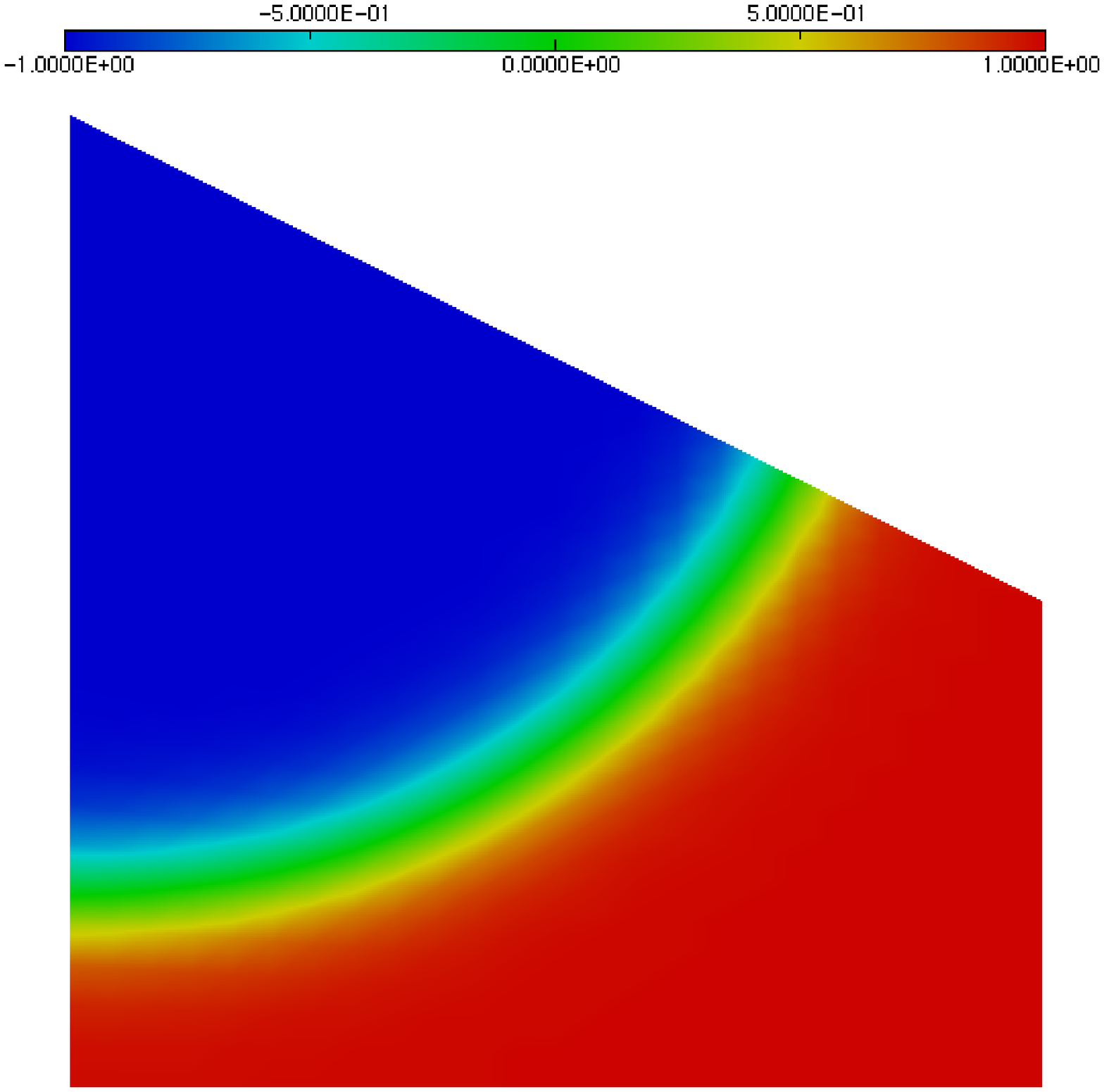}
\hspace{1cm}
\includegraphics[angle=0,width=6cm]{\CercleMD}
\caption{Numerically stable states on the trapezoid.}\label{FIG:29}
\end{center}
\end{figure}

\begin{figure}[htp]
\begin{center}
\includegraphics[angle=0,width=8cm]{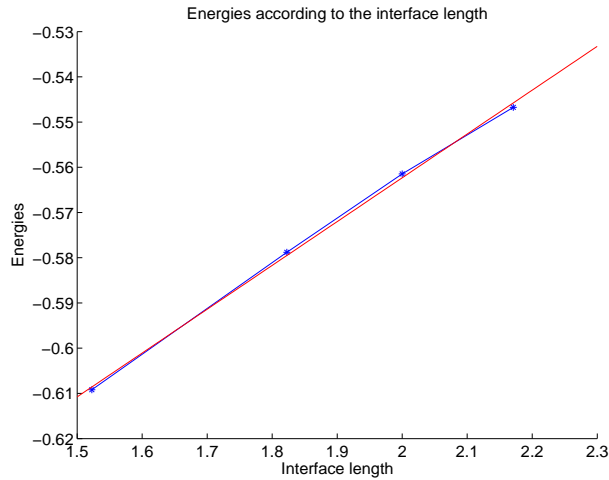}
\caption{Energies of the numerically stable states on the trapezoid according to the length of the interface.}\label{FIG:30}
\end{center}
\end{figure}

 As in the case of the square, we can find numerically stable states corresponding to the next modes in the nomenclature of Maier-Paape and Miller in \cite{MR1950337}. We have found 4 numerically stable states (see Figure \ref{FIG:29}), 
with energies that have been drawn on Figure \ref{FIG:30} against the length of the interface.
We can clear see the linear dependance between the two (the red line is the linear regression according to the least square method).

\clearpage

\noindent{\bf \large Acknowledgments:} We thank Pr. Arnaud Debussche for fruitful discussions about this subject.


\bibliographystyle{abbrv}
\bibliography{Bibliludo}
\end{document}